\documentstyle[amssymb,amsfonts,amscd]{amsart}

\def\G{{\hbox{\bf G}}}

\def\eps{\varepsilon}
\newenvironment{proof}{\noindent {\bf Proof} }{\endprf\par}
\def \endprf{\hfill  {\vrule height6pt width6pt depth0pt}\medskip}
\def\emph#1{{\it #1}}
\def\textbf#1{{\bf #1}}

\newcommand{\BL}{\operatorname{BL}}
\newcommand{\id}{\operatorname{id}}
\newcommand{\tr}{\operatorname{tr}}
\renewcommand{\div}{\operatorname{div}}
\newcommand{\diag}{\operatorname{diag}}
\newcommand{\linearspan}{\operatorname{span}}

\def\u{{\hbox{\bf u}}}
\def\g{{\hbox{\bf g}}}
\def\B{{\hbox{\bf B}}}
\def\K{{\hbox{\bf K}}}

\def\p{{\hbox{\bf p}}}
\def\R{{\hbox{\bf R}}}

\def\A{{\hbox{\bf A}}}
\def\f{{\hbox{\bf f}}}

\theoremstyle{plain}
  \newtheorem{theorem}[subsection]{Theorem}

  \newtheorem{example}[subsection]{Example}
  \newtheorem{proposition}[subsection]{Proposition}
  
  \newtheorem{lemma}[subsection]{Lemma}
  \newtheorem{corollary}[subsection]{Corollary}
  
\theoremstyle{remark}
  \newtheorem{remark}[subsection]{Remark}
  \newtheorem{remarks}[subsection]{Remarks}

\theoremstyle{definition}
  \newtheorem{definition}[subsection]{Definition}

\include{psfig}

\begin{document}

\title[Brascamp--Lieb inequalities]{The Brascamp--Lieb Inequalities: \\ 
Finiteness, Structure and Extremals}

\begin{abstract}
We consider the Brascamp--Lieb inequalities concerning multilinear
integrals of products of functions in several dimensions. We give
a complete treatment of the issues of finiteness of the constant,
and of the existence and uniqueness of centred gaussian extremals.
For arbitrary extremals we completely address the issue of
existence, and partly address the issue of uniqueness. We also
analyse the inequalities from a structural perspective. 
Our main tool is a monotonicity formula for positive solutions to
heat equations in linear and multilinear settings, which was first
used in this type of setting by Carlen, Lieb, and Loss \cite{CLL}.
In that paper, the heat flow method was used to obtain the rank one case
of Lieb's fundamental theorem concerning exhaustion by gaussians; we extend the
technique to the higher rank case, giving two new proofs of the general rank case
of Lieb's theorem.
\end{abstract}

\author{Jonathan Bennett}
\address{School of Mathematics, University of Birmingham, Birmingham, B15 2TT, UK}
\email{J.Bennett@@bham.ac.uk}

\author{Anthony Carbery}
\address{School of Mathematics, University of Edinburgh, Edinburgh, EH9 3JZ, UK}
\email{A.Carbery@@ed.ac.uk}

\author{Michael Christ}
\address{
        Department of Mathematics,
        University of California, 
        Berkeley, CA 94720-3840, USA}
\email{mchrist@@math.berkeley.edu}
\thanks{The first author was supported by EPSRC Postdoctoral Fellowship GR/
S27009/02, the second partially by a Leverhulme Study Abroad Fellowship and
third in part by NSF grant DMS-040126}

\author{Terence Tao}
\address{Department of Mathematics, UCLA, Los Angeles CA 90095-1555}
\email{tao@@math.ucla.edu}

\maketitle

\section{Introduction}

Important inequalities such as the multilinear H\"{o}lder inequality, the sharp
Young convolution inequality and the Loomis--Whitney inequality find their
natural generalisation in the \emph{Brascamp--Lieb inequalities}, which we now describe.

\begin{definition}[Brascamp--Lieb constant]\label{BL-def}  Define a \emph{Euclidean space}
to be a finite-dimensional
real Hilbert space\footnote{It is convenient to work with arbitrary finite-dimensional
Hilbert spaces instead of just copies of $\R^n$ in order to take advantage of invariance
under Hilbert space isometries, as well as such operations as restriction of a Hilbert
space to a subspace, or quotienting one Hilbert space by another.  In fact one could
dispense with the inner product structure altogether and work with finite-dimensional
vector spaces with a Haar measure $dx$, but as the notation in that setting is less
familiar, especially when regarding heat equations on such domains, we shall
retain the inner product structure for notational convenience.} $H$, endowed with the
usual Lebesgue measure $dx$; for instance, $\R^n$ is a Euclidean space for any $n$.
If $m \geq 0$ is an integer, we define an \emph{$m$-transformation} to be a triple
$$ \B := (H, (H_j)_{1 \leq j \leq m}, (B_j)_{1 \leq j \leq m})$$
where $H, H_1,\ldots,H_m$ are Euclidean spaces and for each $j$, $B_j: H \to H_j$ is a
linear transformation.  We refer to $H$ as the \emph{domain}
of the $m$-transformation $\B$.  We say that an $m$-transformation is \emph{non-degenerate}
if all the $B_j$ are surjective (thus $H_j = B_j H$) and the common kernel is trivial
(thus $\bigcap_{j=1}^m \ker(B_j) = \{0\}$). We define an \emph{$m$-exponent} to be an
$m$-tuple $\p = (p_j)_{1 \leq j \leq m} \in \R_+^m$ of non-negative real numbers.  We define
a \emph{Brascamp--Lieb datum} to be a pair $(\B,\p)$, where $\B$ is an
$m$-transformation and $\p$ is an $m$-exponent for some integer $m \geq 0$.
When we are in a situation which involves a Brascamp--Lieb datum $(\B,\p)$, it is always
understood that the objects $H$, $H_j$, $B_j$, $p_j$
denote the relevant components of this Brascamp--Lieb datum.  If $(\B,\p)$ is a Brascamp--Lieb
datum, we define an \emph{input} for $(\B,\p)$ to be an $m$-tuple
${\bf{f}} := (f_j)_{1 \leq j \leq m}$ of nonnegative measurable functions
$f_j: H_j \to \R^+$ such that $0 < \int_{H_j} f_j < \infty$,
and then define the quantity $0 \leq \BL(\B,\p;\f) \leq +\infty$
by the formula $$ \BL(\B,\p; \f) :=
\frac{\int_H \prod_{j=1}^m (f_j \circ B_j)^{p_j}}{\prod_{j=1}^m (\int_{H_j} f_j)^{p_j}}.$$
Note that if $(\B,\p)$ is non-degenerate and the $f_j$ are bounded with compact support, then $\BL(\B,\p;\f) < +\infty$.
We then define the \emph{Brascamp--Lieb constant} $\BL( \B, \p ) \in (0,+\infty]$ to be the
supremum of $\BL(\B,\p;\f)$ over all inputs $\f$.
Equivalently, $\BL(\B,\p)$ is the smallest constant for which the $m$-linear Brascamp--Lieb
inequality
\begin{eqnarray}\label{BL}
\int_{H} \prod_{j=1}^m (f_j \circ B_j)^{p_j}
\leq \BL( \B, \p ) \prod_{j=1}^m (\int_{H_j} f_j)^{p_j}
\end{eqnarray}
holds for nonnegative measurable functions $f_j: H_j \to \R^+$.
\end{definition}

\begin{remark}
By testing \eqref{BL} on functions which are strictly positive near the origin, one can
easily verify that the Brascamp--Lieb
constant must be strictly positive, though it
can of course be infinite.
We give this definition assuming only that the inputs $f_j$ are non-negative measurable, but
it is easy to see (using Fatou's lemma) that
one could just as easily work with strictly positive Schwartz functions with no change in
the Brascamp--Lieb constant.
One can of course define $\BL(\B,\p)$ when $\B$ is degenerate but it is easily seen that
this constant is infinite in that case (see also Lemma \ref{necc}).  Thus we shall often
restrict our attention to non-degenerate Brascamp--Lieb data.
\end{remark}

We now give some standard examples of Brascamp--Lieb data and their associated Brascamp--Lieb
constants.

\begin{example}[H\"older's inequality]\label{holder-ex} If $\B$ is the non-degenerate
$m$-transformation
$$ \B := (H, (H)_{1 \leq j \leq m}, (\id_H)_{1 \leq j \leq m})$$
for some Euclidean space $H$ and some $m \geq 1$, where $\id_H: H \to H$
denotes the identity on $H$,
then the multilinear H\"older inequality asserts that $\BL( \B, \p )$ is
equal to 1 when $p_1 + \ldots + p_m = 1$, and is equal to $+\infty$ otherwise.
\end{example}

\begin{example}[Loomis--Whitney inequality]\label{lw-ex} If $\B$ is the non-degenerate
$n$ - transformation
$$ B := (\R^n, (e_j^\perp)_{1 \leq j \leq n}, (P_j)_{1 \leq j \leq n})$$
where $e_1,\ldots,e_n$ is the standard basis of $\R^n$, $e_j^\perp \subset \R^n$ is the orthogonal complement of $e_j$,
and $P_j: \R^n \to e_j^\perp$ is the orthogonal projection onto $e_j$,
then the Loomis--Whitney inequality \cite{LW} can be interpreted as an assertion that $\BL( \B, \p ) = 1$
when $\p = (\frac{1}{n-1}, \ldots, \frac{1}{n-1})$, and is infinite for any other value of $\p$.
For instance, when $n=3$ this inequality asserts that
\begin{equation}\label{lw-eq}
 \int\int\int f(y,z)^{1/2} g(x,z)^{1/2} h(x,y)^{1/2}\ dx dy dz 
 \leq \|f\|_{L^1({\bf R}^2)}^{1/2} \|g\|_{L^1({\bf R}^2)}^{1/2} \|h\|_{L^1({\bf R}^2)}^{1/2}
 \end{equation}
whenever $f,g,h$ are non-negative measurable functions on $\R^2$.  More generally, Finner \cite{Finner} established
multilinear inequalities of Loomis-Whitney type involving orthogonal projections to co-ordinate subspaces.
\end{example}

\begin{example}[Sharp Young inequality]\label{bl-ex}
The sharp Young inequality (\cite{Be}, \cite{BL1}) can be viewed as an assertion that
if $\B$ is the non-degenerate $3$-transformation
$$ \B := (\R^d \times \R^d, (\R^d)_{1 \leq j \leq 3}, (B_j)_{1 \leq j \leq 3})$$
where $d \geq 1$ is an integer and the maps $B_j: \R^d \times \R^d \to \R^d$ are
defined for $j=1,2,3$ by
$$ B_1(x,y) := x; \quad B_2(x,y) := y; \quad B_3(x,y) := x-y$$
then we have
$$\BL( \B, (p_1,p_2,p_3) ) = \left(\prod_{j=1}^3 \frac{(1-p_j)^{1-p_j}}{p_j^{p_j}}\right)^{d/2}$$
if $p_1 + p_2 + p_3 = 2$ and $0 \leq p_1,p_2,p_3 \leq 1$,
with
$\BL(\B, (p_1,p_2,p_3)) = +\infty$ for any other
values of $(p_1,p_2,p_3)$.  (See Example \ref{bl-ex-2}.)
\end{example}

\begin{example}[Geometric Brascamp--Lieb inequality]\label{gbl}  Let $\B$ be the 
non-degenerate $m$-transformation
$$ \B := (H, (H_j)_{1 \leq j\leq m}, (B_j)_{1 \leq j \leq m})$$
where $H$ is a Euclidean space, $H_1,\ldots,H_m$ are subspaces of
$H$, and $B_j: H \to H_j$ are orthogonal projections, thus the
adjoint $B_j^*: H_j \to H$ is the inclusion map.  The geometric
Brascamp--Lieb inequality of Ball \cite{Ball} (and later
generalised by Barthe \cite{Barthe}) asserts that $\BL(\B, \p ) =
1$ whenever $\p = (p_1,\ldots,p_m) \in \R_+^m$ obeys the identity
$\sum_{j=1}^m p_j B_j^* B_j = 1$. Note that equality is obtained
when we have $f_j(x) = \exp(-\pi \|x\|_{H_j}^2)$ for all $j$. This
significantly generalises the Loomis--Whitney inequality (Example
\ref{lw-ex}) and the H\"older inequality (Example
\ref{holder-ex}).
\end{example}

\begin{example}[Rank-one Brascamp--Lieb inequality]\label{bbl} Let $\B$ be the non-degenerate
$m$-transformation
$$ \B := (H, (\R)_{1 \leq j\leq m}, (v_j^*)_{1 \leq j \leq m})$$
where $H$ is a Euclidean space, $v_1,\ldots,v_m$ are non-zero vectors in $H$ which span $H$,
and $v_j^*: H \to \R$ is the corresponding linear functional
$v_j^*(x) := \langle v_j, x \rangle_H$.  Then the work of Barthe \cite{Barthe} (see also \cite{CLL}) shows that $\BL(\B,\p)$ is finite if and only if $\p \in \R^m_+$ lies
in the convex polytope whose vertices are the points $(1_{j \in I})_{1 \leq j \leq m}$, where $I$ is a subset of $\{1,\ldots,m\}$ such that
the vectors $(v_j)_{j \in I}$ form a basis of $H$ (in particular, this forces $|I| = \dim(H)$).  Furthermore, if $\p$ lies in the $(m-1)$-dimensional
interior of this polytope, then equality in \eqref{BL} can be attained.  We shall reprove these statements as Theorem \ref{barthe-thm}.
\end{example}

\begin{remark}
The rank one case can differ dramatically from the general case. In particular, rearrangement
inequalities such as that of Brascamp, Lieb, and Luttinger \cite{BLL}
apply for rank one, that is, when all $H_j$ have dimensions one,
but for higher rank are only very rarely applicable.
\end{remark}

It is thus of interest to compute the Brascamp--Lieb constants $\BL(\B,\p)$ explicitly,
or at least to determine under what conditions these constants are finite.  A fundamental
theorem of Lieb \cite{L}
shows that these constants are exhausted by centred gaussians.
More precisely, given any positive definite transformation\footnote{By positive definite,
we mean that the transformation is self-adjoint and that the associated quadratic form
$\langle Ax, x \rangle$ is positive definite.  It is in fact slightly more natural
to view $A$ as a transformation from $H$ to $H^*$,
the dual of $H$, rather than $H$ itself, but of course since $H$ and $H^*$ are
canonically identifiable using the Hilbert space structure we will usually not
bother to enforce the distinction between $H$ and $H^*$.} $A: H \to H$, we consider
the associated
gaussian $\exp( - \pi \langle Ax, x \rangle_H )$.  As is well-known we have the formula
\begin{equation}\label{hilbert}
\int_{H} \exp( - \pi \langle Ax, x \rangle )\ dx =
(\mbox{$\det_{H}$} A)^{-1/2},
\end{equation}
where $\det_H$ is the determinant associated to transformations on the Euclidean space $H$.
Define a \emph{gaussian input} for $(\B,\p)$ to be any $m$-tuple $\A =
(A_j)_{1\leq j \leq m}$ of 
positive definite
linear transformations $A_j: H_j \to H_j$.
If we now test the inequality \eqref{BL} with the input
$(\exp( - \pi \langle A_j x, x \rangle_H))_{1 \leq j \leq m}$ for some arbitrary
positive-definite transformations $A_j: H_j \to H_j$, we conclude that
$$ \BL( \B,\p ) \geq \BL_\g(\B,\p;\A)$$
for arbitrary gaussian input $\A = (A_j)_{1 \leq j \leq m}$, where
$\BL_\g(\B,\p;\A) \in (0,+\infty)$ is the quantity
\begin{equation}\label{BLfunctional(not)}
\BL_\g( \B,\p; \A) := \left( \frac{\prod_{j=1}^m (\det_{H_j}
A_j)^{p_j}} {\det_H (\sum_{j=1}^m p_j B_j^* A_j B_j)}
\right)^{1/2},
\end{equation}
and $B_j^*: H_j \to H$ is the adjoint of $H$.
In particular, if we define
\begin{equation}\label{BLfunctional}
\BL_\g( \B,\p) = \sup \{ \BL_\g(\B,\p;\A): \A \hbox{ is a gaussian input for } (\B,\p)\}
\end{equation}
then we have
\begin{equation}\label{blg-trivial}
 \BL( \B,\p ) \geq \BL_\g( \B,\p )
 \end{equation}
for arbitrary Brascamp--Lieb data $(\B,\p)$.

Lieb \cite{L} showed that the above inequality is in fact an equality:

\begin{theorem}[Lieb's theorem]\label{lieb-thm}\cite{L}  For any Brascamp--Lieb datum $(\B,\p)$, we have
$\BL( \B, \p ) = \BL_\g( \B, \p )$.
\end{theorem}

\begin{remark}
Lieb's theorem was proved in the rank-one and certain intermediate cases by Brascamp
and Lieb \cite{BL1} in part using rearrangement techniques from \cite{BLL}. Lieb gave
the first full proof in \cite{L} using, {\em{inter alia}}, a clever $O(2)$-invariance
and arguments related to the central limit theorem. Barthe gave an alternative proof
using deep ideas of transportation of mass due to Brenier, McCann and Caffarelli;
see \cite{Barthe} and the references therein.  More relevantly to our own approach,
Carlen, Lieb and Loss \cite{CLL} gave a proof of Lieb's theorem in the rank one case
using heat flow methods.  Our own arguments, though rediscovered independently, can be viewed
as an extension of the arguments in \cite{CLL} to the higher rank case.
\end{remark}

The question of understanding the Brascamp--Lieb constants is solved by Theorem
\ref{lieb-thm}, in the sense that the task
is reduced to the simpler task of understanding the (in principle computable)
gaussian Brascamp--Lieb constants.  However, there are still
a number of issues that are not easily resolved just from this theorem alone; for instance,
it is not immediately obvious what the
necessary and sufficient conditions are for either $\BL( \B, \p )$ or $\BL_\g( \B, \p )$
to be finite.  Also, this theorem does not characterise
the extremals for \eqref{BL} or for \eqref{BLfunctional}, or even address whether
such extremals exist at all.  

More precisely, we formulate
\begin{definition}[Extremisability]\label{extremisability} 
A Brascamp--Lieb datum 
is said to be 
\emph{extremisable} 
if $\BL(\B,\p)$ is finite and there exists an input $\f$ for which
$\BL(\B,\p) = \BL(\B,\p; \f)$.

A Brascamp--Lieb datum 
is said to be \emph{gaussian-extremisable} if there exists a gaussian
input $\A$ for which $\BL_\g(\B,\p) = \BL_\g(\B,\p;\A)$ 
(so in particular, $\BL_\g(\B,\p)$ is finite). 
\end{definition}
Note that gaussian-extremisability does not formally imply extremisability.
Nonetheless,
Lieb's theorem asserts that every gaussian-extremisable datum is also
extremisable. However it does not give the converse (which turns out to be true,
see Proposition \ref{extremisor} below). 

In this paper we shall address these issues as follows. Firstly, we shall
use the multilinear heat flow monotonicity formula
technique to give two fully self-contained proofs of the general-rank case of Lieb's theorem (in
Section \ref{necsuf-sec} and Section \ref{regular-sec}
respectively).  This technique was first employed in this setting
\cite{CLL} to establish the analogues of most of the results here 
in the rank-one case $\dim(H_j) = 1$ (and also for more
general domains than Euclidean spaces), and our arguments are thus an extension of those in 
\cite{CLL}.  (See also \cite{Barthe-new} for another recent application of the heat flow method to
Brascamp--Lieb type inequalities.) Secondly, we establish a 
characterisation of when the Brascamp--Lieb constants are finite; in 
fact this will be achieved simultaneously with the first of our 
two new proofs of Lieb's theorem. Thirdly, we shall give necessary 
and sufficient conditions for the Brascamp--Lieb data $(\B,\p)$ to
be gaussian-extremisable. Fourthly we shall give necessary and sufficient 
conditions for the Brascamp--Lieb data $(\B,\p)$ to possess unique gaussian
extremisers (up to trivial symmetries). We shall achieve these results and 
others partly with the aid of a two-stage structural perspective on the 
Brascamp--Lieb inequalities \eqref{BL}.

In order to describe our criterion for finiteness we first need a 
definition.  This definition (or more precisely an equivalent formulation of
this definition) was first introduced in the context of the rank one problem in \cite{CLL}.

\begin{definition}[Critical subspace and simplicity]\label{crit}\cite{CLL}
Let  $(\B, \p)$ be a
Brascamp--Lieb datum. A \textit{critical subspace} $V$ for $(\B, \p)$ is a 
non-zero proper
subspace of $H$ such that
$$ \dim(V) = \sum_{j=1}^m p_j \dim(B_j V).$$
The datum $(\B, \p)$ is \textit{simple} if it has no critical subspaces.
\end{definition}

\begin{theorem}\label{sufblprecursor} Let $(\B, \p)$ be a Brascamp--Lieb datum. Then
$\BL_\g(\B,\p)$ is finite if and only if we have the scaling condition
\begin{equation}\label{scaling}
\dim(H) = \sum_j p_j \dim(H_j)
\end{equation}
and  the dimension condition
\begin{eqnarray}\label{dimension}
\dim(V) \leq \sum_{j=1}^m p_j \dim(B_j V) \hbox{ for all subspaces } V \subseteq H.
\end{eqnarray}
Furthermore, if $(\B,\p)$ is simple, then it is gaussian-extremisable.
\end{theorem}

\begin{remark}  In the rank-one case, a finiteness criterion for $\BL(\B,\p)$ and $\BL_\g(\B,\p)$ equivalent to
Theorem \ref{sufblprecursor} was given by
Barthe \cite{Barthe}, see Theorem \ref{barthe-thm} below.  
\end{remark}

By Lieb's theorem we have as an immediate corollary:

\begin{theorem}[Finiteness of Brascamp--Lieb constant]\label{sufbl}  Let $(\B, \p)$ be a
Brascamp--Lieb datum.  Then the following three statements
are equivalent.
\begin{itemize}
\item[(a)] $\BL(\B,\p)$ is finite.
\item[(b)] $\BL_\g(\B,\p)$ is finite.
\item[(c)] \eqref{scaling} and \eqref{dimension} hold.
\end{itemize}
Furthermore, if any of (a)-(c) hold, and $(\B,\p)$ is simple, then $(\B,\p)$ is
extremisable.
\end{theorem}

\begin{remarks}  The conditions \eqref{scaling}, \eqref{dimension} imply in particular
that $\B$ is non-degenerate, as can be seen by testing on $V:=H$
and on $V:=\bigcap_{j=1}^m \ker(B_j)$. The deduction of (c) from
(a) or (b) is very easy, see Lemma \ref{necc}. The more difficult
part of the theorem is to establish the reverse implication in
Theorem \ref{sufblprecursor}; we shall do so in Proposition
\ref{sufbl-meat}. One can also easily show that Theorem
\ref{sufbl} implies all the negative results in Examples
\ref{holder-ex}-\ref{gbl}. 
\end{remarks}

\begin{remark}
A key element of our analysis is 
a certain factorisation of the inequality through critical subspaces. 
This factorisation method played a similarly key role in the work of \cite{CLL} in the rank one case, and 
also in the work of Finner \cite{Finner}, who analysed the case of orthogonal projections to co-ordinate subspaces.
One can view this factorization method as generalizing the arguments of Loomis and Whitney \cite{LW}.
In the companion paper \cite{BCCT} we give an alternative proof 
of the equivalence (c) $\iff$ (a) in Theorem \ref{sufbl}, 
without recourse to Lieb's theorem,  heat flow deformation,
or related techniques, based on this factorisation together with
multilinear interpolation.  
Our first proof of Theorem \ref{lieb-thm} (given in Section
\ref{necsuf-sec}) combines heat flows with the factorisation, while
a second proof (in Section \ref{regular-sec}) is a pure heat flow argument.
The factorisation also plays a central role in our structural perspective. 
\end{remark}

\begin{remark} Theorem \ref{sufbl} implies in particular that for any fixed $\B$, the set of
all $\p$ for which $\BL(\B,\p)$ or $\BL_\g(\B,\p)$ is finite is a convex polytope, with the
faces determined by the natural numbers $n$, $n_1,\ldots,n_m$ for which an $n$-dimensional
subspace $V$ of $H$ exists with $n_j$-dimensional images $B_j V$.  
This is a purely geometric
condition, which can in principle be computed algebraically. However, the problem
of determining which $n,n_1,\ldots,n_m$ are attainable seems to be a problem in Schubert
calculus, and given the algebraic richness and complexity of this calculus (see e.g.
\cite{Belkale}), a fully explicit and easily computable description of these numbers
(and hence of the above polytope) may be too ambitious to hope for in general. However
in the rank-one case (Example \ref{bbl}) a concrete description of the polytope
has been given by Barthe \cite{Barthe}. As already observed in \cite{CLL}, Barthe's description
of the polytope is equivalent to that given by Theorem \ref{sufbl}; 
we recall this equivalence in Section \ref{necsuf-sec}.
\end{remark}

Our proofs of Lieb's theorem and Theorem \ref{sufbl} will occupy Sections
\ref{geom-sec}-\ref{necsuf-sec}.  However our methods provide 
a quite short proof of the geometric Brascamp--Lieb inequality (Example \ref{gbl}) 
which does not rely on Lieb's theorem; this is Proposition \ref{gbl-prop}.
Moreover, the geometric Brascamp--Lieb inequality, combined with
a simple linear change of variables argument, is already strong enough to obtain
many of the standard applications of the Brascamp--Lieb inequality, such as the
sharp Young inequality (Example \ref{bl-ex}), as well as the special case of
Theorem \ref{sufbl} and Lieb's theorem when $(\B,\p)$ is simple.  It will also be
one of two pillars of our first proof of Lieb's theorem, the other pillar being the
aforementioned factorisation argument.

In the work of Barthe \cite{Barthe} (Example \ref{bbl}), a structural analysis 
of the rank-one Brascamp--Lieb functional \eqref{BLfunctional(not)} was given, 
involving in particular a decomposition of the rank-one Brascamp--Lieb data 
into \emph{indecomposable} components. Barthe also considered the question of 
when extremals to \eqref{BLfunctional} exist, and when they are unique. The
uniqueness question was answered based upon this decomposition and
upon an explicit algebraic description of the gaussian
Brascamp--Lieb constant \eqref{BLfunctional} in the rank-one case.
Part of the purpose of this paper is to extend (albeit by different methods) Barthe's
programme to the higher rank case. Barthe's decomposition of
Brascamp--Lieb data is different from the factorisation method used here (and also in \cite{CLL}, \cite{Finner}); the
difference is analogous to the distinction between direct product
and semi-direct product in group theory. Barthe's decomposition
depends only on the data $\B$ while ours depends also upon
$\p$. We will combine a higher-rank analogue of Barthe's
decomposition of Brascamp--Lieb data into indecomposable
components with our notion of factorisation in order
to answer the question of when Brascamp--Lieb data is
gaussian-extremisable. For precise definitions see Section \ref{direct-sec}
below.

\begin{theorem} \label{gaussext} Let $(\B, \p)$ be a Brascamp--Lieb datum for 
which \eqref{scaling} and \eqref{dimension} hold. Then $(\B, \p)$
is gaussian-extremisable if and only if each indecomposable component for $\B$ is simple.
\end{theorem}

For the full story see Theorem \ref{extremiser} below. It is also shown in Proposition
\ref{extremisor} that extremisability is equivalent to gaussian-extremisability.
This latter equivalence uses Lieb's theorem and an iterated convolution
idea of Ball. Thus extremisability can be viewed as a kind of
\emph{semisimplicity} for the Brascamp--Lieb datum.

Theorem \ref{extremiser} gives a satisfactory description of the Brascamp--Lieb data for
which extremisers or gaussian extremisers exist. We can also show that
these gaussian extremisers are \emph{unique} (up to scaling) if and only if the data is simple:

\begin{theorem} \label{gaussuniq} Let $(\B, \p)$ be Brascamp--Lieb 
datum with $p_j >0$ for all $j$. Then gaussian extremizers for $(\B,\p)$ exist and are unique (up to scaling)
if and only if $(\B,\p)$ is simple.
\end{theorem}

See Corollary \ref{gauss-unique}.  As for the uniqueness (up to trivial symmetries) 
of general extremisers, this problem seems to be more difficult, and we have only 
a partial result (Theorem \ref{extreme-unique}).

As mentioned above, Carlen, Lieb and Loss in \cite{CLL} have
introduced the idea of using heat flow in the
Brascamp--Lieb context in the rank-one case (and also on the
sphere ${\bf S}^n$ instead of a Euclidean space). We remark also that prior to the work in \cite{CLL},
the known proofs of versions of
Lieb's theorem relied on methods such as rearrangement
inequalities \cite{BLL} and mass transfer inequalities
\cite{Barthe}.  While such methods are similar in spirit to that
of heat flows -- in that they are all ways of deforming
non-extremal solutions to extremal gaussians -- they are not
identical.  For instance the heat flow method will continuously
deform sums of gaussians to other sums of gaussians, whereas
continuous mass transfer does not achieve this; and it is not
obvious how to perform rearrangement in a continuous manner.

{\bf{Guide to the paper.}} In Section \ref{geom-sec} we introduce
the heat flow method and use it to prove the geometric
Brascamp--Lieb inequality. In Section 3 we consider the
gaussian-extremisable case and give a characterisation of gaussian
extremisers, recovering the sharp Young inequality Example
\ref{bl-ex} as an immediate application. Our first approach to the
structure of the Brascamp--Lieb inequalities, via 
factorisation, is taken up in Section 4 where we also establish the
necessity of \eqref{dimension} and \eqref{scaling} for
$\BL(\B,\p)$ or $\BL_{\g}(\B,\p)$ to be finite. In Section
\ref{necsuf-sec} we establish the sufficiency of these conditions
and prove Theorem \ref{lieb-thm} and Theorem \ref{sufbl} in the
general case. In Section 6 we prove that extremisability and
gaussian-extremisability are equivalent. Our second approach to
structural issues is made in Section \ref{direct-sec} where the
main characterisation of extremisability is also given. In Section
\ref{regular-sec}  we examine variants of the Brascamp--Lieb
inequalities for regularised inputs in a gaussian setting, and this leads to 
our second proof of Lieb's theorem via a purely heat-flow method. As a further consequence 
of our analysis, we give versions of Theorem \ref{lieb-thm} and Theorem \ref{sufbl} in a 
gaussian-localised setting in Corollary \ref{llt} and Theorem \ref{sufbl-local} respectively. 
Corollary \ref{llt} is very much in the spirit of Lieb's original paper \cite{L}. 
In Section 9 we discuss uniqueness of
extremals. Finally, in Section \ref{sliding-sec} we make some
remarks on the heat flow method in so far as it applies in
non-gaussian contexts. We give a general monotonicity formula for
log-concave kernels (the class of which of course includes
gaussians) in Lemma \ref{logconcave}. 

Our work in this paper was motivated by the forthcoming article \cite{BCT} in
which we obtain almost optimal results for the multilinear Kakeya and
restriction problems. It was in this context that we rediscovered the
applicability of heat flows in multilinear inequalities.

{\bf Acknowledgements.}  We are grateful to Assaf Naor for bringing \cite{CLL}
to our attention.  We also thank Rodrigo Ba\~nuelos for some corrections, and
Eric Carlen for help with the references.

\section{The geometric Brascamp--Lieb inequality}\label{geom-sec}

In this section we establish Lieb's theorem in the ``geometric'' case (see Proposition \ref{gbl-prop}).  We begin with a definition.

\begin{definition}[Geometric Brascamp--Lieb data]  A Brascamp--Lieb datum $(\B,\p)$ is said to be \emph{geometric}
if we have $B_j B_j^* = \id_{H_j}$ (thus $B_j^*$ is an isometry) and
\begin{equation}\label{pbb}
\sum_{j=1}^m p_j B_j^* B_j = \id_H.
\end{equation}
\end{definition}

\begin{remark}\label{blg-remark} The data in Example \ref{gbl} is of this type.  Conversely, if $(\B,\p)$ is geometric, and we let $E_j$ be the range of
the isometry $B_j^*$, then we may identify $H_j$ with $E_j$, and $B_j$ with the orthogonal projection from $H$ to $E_j$, thus
placing ourselves in the situation of Example \ref{gbl}.
The ``geometry'' here refers to Euclidean or Hilbert space geometry, since the inner product structure of $H$ is clearly
being used (for instance, to define the notion of orthogonal projection).  We remark 
that the condition $B_j B_j^* = \id_{H_j}$ forces $B_j$
to be surjective, and \eqref{pbb} implies that 
$\bigcap_{j=1}^m \ker(B_j) = \{0\}$, and thus geometric Brascamp--Lieb 
data is always non-degenerate.  Furthermore, by taking traces of \eqref{pbb} we conclude \eqref{scaling}.  More generally, if we let $V$ be any subspace of $H$, and let $\Pi: H \to H$ be the orthogonal projection onto $V$, then
by multiplying \eqref{pbb} by $\Pi$ and taking traces we conclude
$$ \sum_{j=1}^m p_j \tr(\Pi B_j^* B_j) = \dim(V).$$
Since $\Pi B_j^* B_j$ is a contraction and has range $B_j V$, we conclude \eqref{dimension} also.  
Thus the assertion that
geometric Brascamp--Lieb data have finite Brascamp--Lieb constants 
is a special case of
Theorem \ref{sufbl}; 
we establish this special case in Proposition \ref{gbl-prop} below.
\end{remark}

\begin{remark}  Let $(\B,\p)$ be a Brascamp--Lieb datum with $p_j > 0$ for all $j$,
and let $\bigoplus_{j=1}^m p_j H_j$ be the Cartesian product
$\prod_{j=1}^m H_j$ endowed with the inner product $\langle x, y
\rangle_{\bigoplus_{j=1}^m p_j H_j} := \sum_{j=1}^m p_j \langle
x_j, y_j \rangle_{H_j}$; this is thus a Euclidean space. Let $V$
be the image of $H$ in $\bigoplus_{j=1}^m p_j H_j$ under the map
$x \mapsto (B_j x)_{1 \leq j \leq m}$, and let $\pi_j: V \to H_j$
be the projection maps induced by restricting the projection maps
on $\bigoplus_{j=1}^m p_j H_j$ to $V$.  Then the statement that
$(\B,\p)$ is geometric is equivalent to the assertion that the map
$x \mapsto (B_j x)_{1 \leq j \leq m}$ is a Euclidean isomorphism 
from $H$ to $V$, and that the projection maps are co-isometries
(i.e. surjective partial isometries, or equivalently their
adjoints are isometries).  Thus geometric Brascamp--Lieb data can
be thought of as subspaces in $\bigoplus_{j=1}^m p_j H_j$ which
lie ``co-isometrically'' above each of the factor spaces $H_j$.
\end{remark}

Our approach to the geometric Brascamp--Lieb inequality is via monotonicity formulae, which as observed in \cite{CLL} seems to be especially
well suited for the geometric Brascamp--Lieb setting.
Abstractly speaking, if one wishes to prove an inequality
of the form $A \leq B$, one can do so by constructing a monotone non-decreasing quantity $Q(t)$ which equals $A$ when $t = 0$ (say)
and equals $B$ when $t = \infty$ (say).  It is thus of interest to have a general scheme for generating such monotone quantities, which we will also
use later in this paper to deduce some variants of the Brascamp--Lieb inequalities.  We begin 
with the following simple lemma.

\begin{lemma}[Monotonicity for transport equations]\label{transport-lemma}  Let $I \subseteq \R$ be a time interval, let $H$ be a Euclidean space,
let $u: I \times H \to \R^+$ be a smooth non-negative function, and $\vec v: I \times H \to H$ be a smooth vector field, such that $\vec v u$ is rapidly decreasing at spatial infinity $x \to \infty$ locally uniformly on $I$.  Let $\alpha \in \R$ be fixed.  Suppose that we have the transport inequality
\begin{equation}\label{transport-a}
\partial_t u(t,x) + \div( \vec v(t,x) u(t,x) ) \geq \alpha u / t
\end{equation}
for all $(t,x) \in I\times H$, where $\div$ of course is the
spatial divergence on the Euclidean space $H$.  Then the quantity
$$ Q(t) := t^{-\alpha} \int_{H} u(t,x)\ dx \in [0,+\infty]$$
is non-decreasing in time.  Furthermore, if \eqref{transport-a} holds with strict inequality for all $x,t$, then $Q$ is strictly
increasing in time.  Similarly, if the signs are reversed in \eqref{transport-a}, then $Q(t)$ is now non-increasing in time.
\end{lemma}

\begin{remark} We allow the velocity field $\vec v$ to vary in both space and time.  This is important as we will typically be dealing with functions $u$ which solve a heat equation such as $\partial_t u = \Delta u$, in which case the natural velocity vield $\vec v$ is given by
$\vec v = - \nabla \log u$.
\end{remark}

\begin{proof}  By multiplying \eqref{transport-a} by $t^\alpha$ and replacing $u$ by $\tilde u := t^{-\alpha} u$, we may reduce
to the case $\alpha = 0$.  Let $t_1 < t_2$ be two times in $I$.  From Stokes' theorem we have
\begin{eqnarray*}
\begin{aligned}
\int_{H} u(t_2,x) \psi(x)&\ dx - \int_{H} u(t_1,x) \psi(x)\ dx \\
&=\int_{t_1}^{t_2} \int_{H} (\partial_t u(t,x) \psi(x) + 
\div( \psi(x) \vec v(t,x) u(t,x))\ dx dt
\end{aligned}
\end{eqnarray*}
for any non-negative smooth cutoff function $\psi$.  Using the product rule and \eqref{transport-a} we conclude
$$ \int_{H} u(t_2,x) \psi(x)\ dx - \int_{H} u(t_1,x) \psi(x)\ dx \geq \int_{t_1}^{t_2} \int_{H} \langle \nabla \psi(x), \vec v(t,x) u(t,x) \rangle\ dx.$$
Letting $\psi$ approach the constant function $1$ and using the hypothesis that $\vec v u$ is rapidly decreasing at spatial infinity, uniformly in $[t_1,t_2]$, we obtain the claim.
\end{proof}

We now generalise this result substantially.

\begin{lemma}[Multilinear monotonicity for transport equations]\label{multitransport}
Let $p_1,\ldots,p_m > 0$ be real exponents, let $H$ be a Euclidean space, and for each $1 \leq j \leq m$
let $u_j: \R^+ \times H \to \R^+$ be a smooth strictly positive function, and 
$\vec v_j: \R^+ \times H \to H$ be a smooth vector field.  Let $\alpha \in \R$ be fixed.  
Suppose we also have an additional smooth vector field $\vec v: \R^+ \times H \to H$, such 
that $\vec v \prod_{j=1}^m u_j^{p_j}$ is rapidly decreasing in space locally uniformly on $I$, 
and we have the inequalities
\begin{align}
\partial_t u_j(t,x) + \div( \vec v_j u_j(t,x) ) &\geq 0  \hbox{ for all } 1 \leq j \leq m \label{transport-j}\\
\div( \vec v - \sum_{j=1}^m p_j \vec v_j ) &\geq \alpha/t \label{div-pos} \\
\sum_{j=1}^m p_j \langle \vec v - \vec v_j, \nabla \log u_j \rangle_H &\geq 0 \label{log-pos}.
\end{align}
Then the quantity
\begin{equation}\label{qt-def}
 Q(t) := t^{-\alpha} \int_{H} \prod_{j=1}^m u_j(t,x)^{p_j}\ dx
 \end{equation}
is non-decreasing in time.  Furthermore, if at least one of \eqref{transport-j}, \eqref{div-pos}, \eqref{log-pos} holds with strict inequality
for all $x,t$, then $Q$ is strictly increasing in time.  If all the signs in \eqref{transport-j}, \eqref{div-pos}, \eqref{log-pos} are reversed,
then $Q(t)$ is now non-increasing in time.
\end{lemma}

\begin{proof}  We begin with the first claim.  By Lemma \ref{transport-lemma}, it suffices to show that
$$ \partial_t \prod_{j=1}^m u_j^{p_j} + \div( \vec v \prod_{j=1}^m u_j^{p_j} ) \geq \alpha \prod_{j=1}^m u_j^{p_j}.$$
We divide both sides by the positive quantity $\prod_{j=1}^m u_j^{p_j}$ and then use the product rule, to reduce to showing that
$$ \sum_{j=1}^m p_j \frac{\partial_t u_j}{u_j} + \div(\vec v) + \langle \vec v, \sum_{j=1}^m p_j \frac{\nabla u_j}{u_j} \rangle_H \geq \alpha.$$
But the left-hand side can be rearranged as
\begin{align*}
&\sum_{j=1}^m \frac{p_j}{u_j} ( \partial_t u_j(t,x) + \div( \vec v_j u_j(t,x) ) ) \\
&+ \div( \vec v - \sum_{j=1}^m p_j \vec v_j )  \\
&+ \sum_{j=1}^m p_j \langle \vec v - \vec v_j, \nabla \log u_j \rangle_H
\end{align*}
and the claim follows from \eqref{transport-j}, \eqref{div-pos}, \eqref{log-pos}.  A similar argument gives the second claim.
\end{proof}

\begin{remark}
The above lemma shows that in order to show that the quantity \eqref{qt-def} synthesised from various spacetime functions $u_j$
is monotone, one needs to locate velocity fields $\vec v_j$, $\vec v$ obeying the pointwise inequalities \eqref{transport-j}, \eqref{div-pos},
\eqref{log-pos}.  In practice the $u_j$ will be chosen to obey a transport equation $\partial_t u_j + \div( \vec v_j u_j ) = 0$, so that
\eqref{transport-j} is automatically satisfied.  As for the other two inequalities, in the linear case $j=1$ one can obtain
\eqref{log-pos} automatically by setting $\vec v = \vec v_j$, leaving only \eqref{div-pos} to be verified.  In the multilinear case, one would
have to set $\vec v$ to be a suitable average of the $\vec v_j$.  This latter strategy seems to only work well in the case when the $u_j$
solve a heat equation $\partial_t u_j = \div( G_j \nabla u_j )$, since in this case the relevant velocity field $\vec v_j= - G_j \nabla \log u_j$
is related to $\nabla \log u_j$ by a positive definite matrix and one is more likely to ensure \eqref{log-pos} is positive.
\end{remark}

We now give the first major application of the above abstract machinery, by reproving K. Ball's
geometric Brascamp--Lieb inequality (Example \ref{gbl}).  We give some further applications in Section \ref{regular-sec} and Section \ref{sliding-sec}.

\begin{proposition}[Geometric Brascamp--Lieb inequality]\label{gbl-prop}\cite{Ball}, \cite{Barthe}  
Let $(\B,\p)$ be geometric Brascamp--Lieb data. Then
\begin{equation}\label{lieb-geom}
 \BL(\B,\p) = \BL_\g(\B,\p) = 1.
\end{equation}
Furthermore, $(\B,\p)$ is both extremisable and gaussian-extremisable.
\end{proposition}

\begin{proof} We may drop those exponents $p_j$ for which $p_j = 0$ as being irrelevant.  By Remark \ref{blg-remark} we may assume
without loss of generality that $H_j$ is a subspace of $H$ and that $B_j$ is the orthogonal projection from $H$ to $H_j$.
By considering \eqref{BLfunctional(not)} with the gaussian input $(\id_{H_j})_{1 \leq j \leq m}$
we observe that $\BL_\g(\B,\p) \geq 1$.  By \eqref{blg-trivial}
we see that it suffices to show that $\BL(\B,\p) \leq 1$; note that this will imply that $(\id_{H_j})_{1 \leq j \leq m}$ is a gaussian
extremiser and $(\exp(- \pi \|x\|_{H_j}^2))_{1 \leq j\leq m}$ is an extremiser.

By Definition \ref{BL-def}, it thus suffices to show that
for any non-negative measurable functions $f_j: H_j \to \R^+$ we have
\begin{equation}\label{fjpj}
\int_H \prod_{j=1}^m (f_j \circ B_j)^{p_j} \leq \prod_{j=1}^m (\int_{H_j} f_j)^{p_j}.
\end{equation}

By Fatou's lemma it suffices to verify this inequality when the $f_j$ are smooth, rapidly decreasing, and strictly positive.
Now let $u_j: \R^+ \times H \to \R^+$ be the solution to the heat equation Cauchy problem
\begin{align*}
\partial_t u_j(t,x) &= \Delta_H u_j(t,x) \\
u_j(0,x) &= f_j \circ B_j(x)
\end{align*}
where $\Delta_H := \div \nabla$ is the usual Laplacian on $H$.  More explicitly, we have
$$
u_j(t,x) = \frac{1}{(4\pi t)^{\dim(H)/2}} \int_H e^{-\| x-y \|_H^2/4t} f_j(B_j y)\ dy.$$
We can split $y$ into components in $H_j$ and in the orthogonal complement $H_j^\perp$, using Pythagoras' theorem to split $\| x-y\|_H^2$ correspondingly.  The contribution from the orthogonal complement can be evaluated explicitly by \eqref{hilbert}, to obtain
\begin{equation}\label{uj-form}
u_j(t,x) = \frac{1}{(4\pi t)^{\dim(H_j)/2}} \int_{H_j} e^{-\| B_j x - z \|_{H_j}^2/4t} f_j(z)\ dz.
\end{equation}
Alternatively, one can verify that $u_j$ also solves the above Cauchy problem.

In order to apply Lemma \ref{multitransport}, we rewrite the heat equation as a transport equation
$$ \partial_t u_j + \div( \vec v_j u_j ) = 0$$
where $\vec v_j := - \nabla \log u_j$; thus \eqref{transport-j} is trivially satisfied.  Next we set $\alpha := 0$ and
$$ \vec v := \sum_{j=1}^m p_j \vec v_j$$
so that \eqref{div-pos} is also trivially satisfied.

Next, we verify the technical condition that $\vec v \prod_{j=1}^m u_j^{p_j}$ is rapidly decreasing in space.  From \eqref{uj-form}
and the hypothesis that $f_j$ is smooth and rapidly decreasing, we see that for any $t$ in a compact interval in $(0,\infty)$,
$\vec v_j(t,x) = -\nabla u_j/u_j(t,x)$ grows at most polynomially in space.  On the other hand, $u_j$ is bounded and
rapidly decreasing in the directions $\| B_j x \|_{H_j} \to \infty$.  Since $\B$ is non-degenerate,
$\prod_{j=1}^m u_j^{p_j}$ decays rapidly in all spatial directions.  The claim follows.

Now we verify \eqref{log-pos}. Observe that $u_j(t,x)$ depends only on $B_j x$ and not on $x$ itself, which shows that $\vec v_j$ lies in
the range $H_j$ of the projection $B_j^* B_j$.  In particular we have
$$ \nabla \log u_j = - \vec v_j = -B_j^* B_j \vec v_j.$$
Hence we have
$$
\sum_{j=1}^m p_j \langle \vec v - \vec v_j, \nabla \log u_j \rangle_H
=
\sum_{j=1}^m p_j \langle B_j^* B_j(\vec v - \vec v_j), -\vec v_j \rangle_H.$$
Also, from \eqref{pbb} we have
$$ \sum_{j=1}^m p_j B_j^{*}B_{j} (\vec v - \vec v_j) = \vec v - \sum_{j=1}^m p_j B_j^* B_j \vec v_j = \vec v - \sum_{j=1}^m p_j \vec v_j = 0$$
and hence
\begin{equation}\label{vvj}
\sum_{j=1}^m p_j \langle \vec v - \vec v_j, \nabla \log u_j \rangle_H
=
\sum_{j=1}^m p_j \langle B_j^* B_j(\vec v - \vec v_j), (\vec v - \vec v_j) \rangle_H.
\end{equation}
Since the orthogonal projection $B_j^* B_j$ is positive semi-definite, we obtain \eqref{log-pos}.  We may now invoke Lemma \ref{multitransport}
and conclude that the quantity
$$ Q(t) := \int_{H} \prod_{j=1}^m u_j^{p_j}(t,x)\ dx$$
is non-decreasing for $0 < t < \infty$.  In particular we have
$$ \limsup_{t \to 0^+} Q(t) \leq \liminf_{t \to \infty} Q(t).$$
From Fatou's lemma we have
$$ \int_{H} \prod_{j=1}^m (f_j \circ B_j)^{p_j} \leq \limsup_{t \to 0^+} Q(t)$$
(in fact we have equality, but we will not need this), so it will suffice to show that
\begin{equation}\label{qinf}
 \liminf_{t \to \infty} Q(t) \leq \prod_{j=1}^m \left(\int_{H_j} f_j\right)^{p_j}.
\end{equation}
(Again, we will have equality, but we do not need this.)
Using \eqref{uj-form} we can write
$$ Q(t) = \frac{1}{(4\pi t)^{\sum_{j=1}^m p_j \dim(H_j)/2}}
\int_{H} \prod_{j=1}^m \left(\int_{H_j} e^{- \| B_j x - z\|_{H_j}^2/4t} f_j(z)\ dz\right)^{p_j}\ dx.$$
By taking traces of the hypothesis \eqref{pbb} we obtain \eqref{scaling}.
Thus by making the change of variables $x = t^{1/2} w$ we obtain
$$ Q(t) = \frac{1}{(4\pi)^{\dim(H)/2}}
\int_{H} \prod_{j=1}^m \left(\int_{H_j} e^{- \| B_j w - t^{-1/2} z\|_{H_j}^2/4} f_j(z)\ dz\right)^{p_j}\ dw.$$
Since the $f_j$ are rapidly decreasing and $\bigcap_{j=1}^m \ker(P_j) = \{0\}$, we may then use
dominated convergence to conclude
\begin{align*}
\liminf_{t \to \infty} Q(t) &=
\frac{1}{(4\pi)^{\dim(H)/2}}
\int_H \prod_{j=1}^m \left(\int_{H_j} e^{-\| B_j w \|_{H_j}^2/4} f_j(z)\ dz\right)^{p_j}\ dw \\
&=\frac{1}{(4\pi)^{\dim(H)/2}} \prod_{j=1}^m \left(\int_{H_j} f_j\right)^{p_j}
\int_{H} e^{-\langle \sum_{j=1}^m p_j B_j^* B_j w, w \rangle_H/4}\ dw.
\end{align*}
Using \eqref{pbb} and \eqref{hilbert}, the claim \eqref{qinf} follows.
\end{proof}

\begin{remark} It is clear from the above argument that the method extends to the case 
when $m$ is countably or even uncountably infinite; in the latter case the exponents 
$p_j$ will be replaced by some positive measure on the index set that $j$ ranges over.  
In this connection see \cite{Bar3}.
\end{remark}

\section{The gaussian-extremisable case}

We can now prove Theorem \ref{lieb-thm} and Theorem \ref{sufbl} in the gaussian-extremisable 
case.

We begin with a notion of equivalence between Brascamp--Lieb data.

\begin{definition}[Equivalence]\label{equiv-def}  Two $m$-transformations $\B = 
(H, (H_j)_{1 \leq j \leq m},$ $\hfill$ $(B_j)_{1 \leq j \leq m})$ and $\B' = (H', (H'_j)_{1 \leq j \leq m}, (B'_j)_{1 \leq j \leq m})$ are said to be \emph{equivalent} if there exist invertible linear
transformations $C: H' \to H$ and $C_j: H'_j \to H_j$ such that $B'_j = C_j^{-1} B_j C$ for all $j$; we refer to $C$ and $C_j$ as the \emph{intertwining transformations}.  Note in particular that this forces $\dim(H) = \dim(H')$ and $\dim(H_j) = \dim(H'_j)$.
We say that two Brascamp--Lieb data $(\B,\p)$ and $(\B',\p')$ are equivalent if $\B,\B'$ are equivalent and $\p = \p'$.
\end{definition}

\begin{remark} This is clearly an equivalence relation.  Up to equivalence, the only relevant features of a non-degenerate $m$-transformation
are the kernels $\ker(B_j)$ and how they are situated inside $H$.  More precisely, if one fixes the dimensions $n = \dim(H)$ and $n_j = \dim(H_j)$,
then the moduli space of non-degenerate $m$-transformations with these dimensions, quotiented out by equivalence, can be identified
with the moduli space of $m$-tuples of subspaces $V_j$ of $\R^n$, with $\dim(V_j) = n - n_j$ and $\bigcap_{j=1}^m V_j = \{0\}$, quotiented out
by the general linear group $GL(\R^n)$ of $\R^n$.  The problem of understanding this moduli space is part of the more general question of understanding quiver representations, which is a rich and complex subject (see \cite{DW} for a recent survey).
\end{remark}

The Brascamp--Lieb constants of two equivalent Brascamp--Lieb data are closely related:

\begin{lemma}[Equivalence of Brascamp--Lieb constants]\label{scaling-lemma}
Suppose that $(\B,\p)$, $(\B',\p')$ are two equivalent Brascamp--Lieb data, with
intertwining transformations $C: H' \to H$ and $C_j: H'_j \to H_j$.  Then we have
\begin{equation}\label{bl-scaling}
\BL( \B', \p' ) = \frac{ \prod_{j=1}^m |\det_{H'_j \to H_j} C_j|^{p_j} }{|\det_{H' \to H} C|} \BL ( \B,\p )
\end{equation}
and
\begin{equation}\label{blg-scaling}
 \BL_\g( \B', \p' ) = \frac{ \prod_{j=1}^m |\det_{H'_j \to H_j} C_j|^{p_j} }{|\det_{H' \to H} C|} \BL_\g ( \B,\p ).
 \end{equation}
Here of course $\det_{H' \to H} C$ denotes the determinant of the transformation $C: H' \to H$ with respect to the Lebesgue measures on $H', H$,
and similarly for $\det_{H'_j \to H_j} C_j$.  In particular $\BL(\B',\p')$ (or $\BL_\g(\B',\p')$) is finite if and only if
$\BL(\B,\p)$ (or $\BL_\g(\B,\p)$) is finite.

Furthermore, $(\B,\p)$ is extremisable if and only if $(\B',\p')$ is extremisable, and $(\B,\p)$ is gaussian-extremisable if and only if
$(\B',\p')$ are gaussian-extremisable.
\end{lemma}

\begin{proof}
To see \eqref{bl-scaling} we simply apply the change of variables $x = C^{-1} y$ on $\R^n$, and replace an input $(f_j)_{1 \leq j \leq m}$
for $(\B,\p)$ with the corresponding input $(f_j \circ C_j)_{1 \leq j \leq m}$ for $(\B',\p')$.
To see the second identity we similarly replace a gaussian input $(A_j)_{1 \leq j\leq m}$ for $(\B,\p)$ with
the gaussian input $(C_j^* A_j C_j)_{1 \leq j\leq m}$ for $(\B',\p')$.
The corresponding claims about extremisers are proven similarly.
\end{proof}

\begin{remark}  The transformation between $A_j$ and $A'_j$ can be described using the commutative diagram
$$
\begin{CD}
H        @>{B_j}>>  H_j        @>{A_j}>>  H_j^*       @>{B^*_j}>>  H^* \\
@AA{C}A            @AA{C_j}A             @VV{C_j^*}V            @VV{C^*}V\\
H'       @>{B'_j}>> H'_j       @>{A'_j}>>  (H'_j)^*       @>{(B'_j)^*}>>  (H')^*
\end{CD}
$$
where we have suppressed the identification between a Hilbert space $H$ and its
dual $H^*$ to emphasize the self-adjointness of the above diagram, and the horizontal rows do not denote exact sequences.
The above lemma shows that the inner product structure of $H$
is not truly relevant for the analysis of Brascamp--Lieb constants, however we retain this structure in order to take advantage of convenient
notions such as orthogonal complement or induced Lebesgue measure.
\end{remark}

We are now ready to give a satisfactory algebraic characterisation of gaussian 
extremisable Brascamp--Lieb data.

\begin{definition}[Ordering of self-adjoint transformations] If $A: H \to H$ and $B: H \to H$ are two self-adjoint linear transformations on
a Euclidean space $H$, we write $A \geq B$ if $A-B$ is positive semi-definite, and $A > B$ 
if $A-B$ is 
positive definite.  (We recall
that $A \geq B$ and $A \neq B$ do not together imply that $A > B$.)
\end{definition}

\begin{proposition}[Characterisation of gaussian-extremisers]\label{normal-form}
Let $(\B,\p)$ be a Bras- camp--Lieb datum with $p_j > 0$ for all $j$, and let $\A = (A_j)_{1 \leq j \leq m}$ 
be a gaussian input for $(\B,\p)$.
Let $M: H \to H$ be the positive semi-definite transformation $M := \sum_{j=1}^m p_j B_j^* A_j B_j$.
Then the following seven statements are equivalent.

\begin{enumerate}
\item  $\A$ is a global extremiser to \eqref{BLfunctional(not)} (in particular, $\BL_\g(\B,\p) = \BL_\g(\B,\p;\A)$ is finite and $(\B,\p)$ is gaussian-extremisable).
\item  $\A$ is a local extremiser to \eqref{BLfunctional(not)}.
\item  $M$ is invertible, and we have
\begin{equation}\label{pj-ident}
 A_j^{-1} - B_j M^{-1} B_j^* = 0 \hbox{ for all } 1 \leq j \leq m.
\end{equation}
\item  The scaling condition \eqref{scaling} holds, $\B$ is non-degenerate, and
\begin{equation}\label{pj-ident-alt}
 A_j^{-1} \geq B_j M^{-1} B_j^* \hbox{ for all } 1 \leq j \leq m.
\end{equation}
\item  The scaling condition \eqref{scaling} holds, $\B$ is non-degenerate, and
\begin{equation}\label{pj-ident2}
B_j^* A_j B_j \leq M \hbox{ for all } 1 \leq j \leq m.
\end{equation}
\item $(\B,\p)$ is equivalent to a geometric Brascamp--Lieb datum $(\B',\p')$ with intertwining operators $C :=
M^{-1/2}$ and $C_j := A_j^{-1/2}$, and
\begin{align*}
\BL(\B,\p) &=  \BL(\B',\p') \BL_\g(\B,\p;\A);\\
\BL_\g(\B,\p) &=  \BL_\g(\B',\p') \BL_\g(\B,\p;\A).
\end{align*}
\item  $(\B,\p)$ is equivalent to a geometric Brascamp--Lieb datum $(\B',\p')$ with intertwining operators $C :=
M^{-1/2}$ and $C_j := A_j^{-1/2}$, and $H_j$ equal to the range of $C_j$, and
$$ \BL(\B,\p) = \BL_\g(\B,\p) = \BL_\g(\B,\p;\A).$$
\end{enumerate}
\end{proposition}

\begin{remark}
One may wish to view \eqref{pj-ident} as a commutative diagram
$$\begin{CD}
H               @>{M^{-1}}>>   H  \\
@AA{B_j^*}A                @VV{B_j}V    \\
H_j             @>{A_j^{-1}}>> H_j
\end{CD}.
$$
We shall shortly expand this diagram in the proof below as
$$
\begin{CD}
H               @>{M^{-1/2}}>>   H    @>{id}>>                H  @>{M^{-1/2}}>>   H  \\
@AA{B_j^*}A                   @AA{T_j^*}A              @VV{T_j}V                   @VV{B_j}V    \\
H_j             @>{A_j^{-1/2}}>>  H_j @>{id}>>               H_j  @>{A_j^{-1/2}}>>  H_j
\end{CD}.
$$
\end{remark}

\begin{proof}  The implication (a) $\implies$ (b) is trivial.  Now we verify that (b) $\implies$ (c).  From Lemma \ref{necc} we see that $\bigcap_{j=1}^m \ker(B_j) = \{0\}$.
Since the $A_j$ are 
positive definite, this implies that $M$ is also.

Taking logarithms in \eqref{BLfunctional(not)}, we see that $\A > 0$ is a local maximiser for the quantity
$$ (\sum_{j=1}^m p_j \log \mbox{$\det_{H_j}$} A_j) - \log \mbox{$\det_H$} \sum_{j=1}^m p_j B_j^* A_j B_j.$$
Let us fix a $j$ in $\{1,\ldots m\}$, and let $Q_j: H_j \to H_j$
be an arbitrary self-adjoint transformation.  By perturbing $A_j$
by a small multiple of $Q_j$, we conclude that
$$ \frac{d}{d\eps}\left( p_j \log \mbox{$\det_{H_j}$} (A_j + \eps Q_j) - \log \mbox{$\det_H$}(M + \eps p_j B_j^* Q_j B_j)\right)\Bigl|_{\eps = 0} = 0.$$
We can thus subtract off the term $p_j \log \det_{H_j} A_j - \log
\det_H M$, which does not depend on $\eps$, from the above
equation and obtain
$$ \frac{d}{d\eps}\left( p_j \log \mbox{$\det_{H_j}$} (\id_{H_j} + \eps A_j^{-1} Q_j) -
\log \mbox{$\det_H$} (\id_H + \eps p_j M^{-1} B_j^* Q_j
B_j)\right)\Bigl|_{\eps = 0} = 0.$$ A simple Taylor expansion
shows
\begin{equation}\label{taylor}
\frac{d}{d\eps} \log\mbox{$\det_H$} (\id_H + A)|_{\eps = 0} =
\tr_H(A),
\end{equation}
hence
$$ p_j \tr_{H_j}(A_j^{-1} Q_j) - \tr_H( p_j M^{-1} B_j^* Q_j B_j ) = 0.$$
We rearrange this using the cyclic properties of the trace as
$$ \tr_{H_j}( (A_j^{-1} - B_j M^{-1} B_j^*) Q_j ) = 0.$$
Since $Q_j$ was an arbitrary self-adjoint transformation, and
$A_j^{-1} - B_j M^{-1} B_j^*$ is also self-adjoint, we conclude
\eqref{pj-ident}.

Now we verify that (c) $\implies$ (d). The implication of
\eqref{pj-ident-alt} from \eqref{pj-ident} is trivial.  From the
invertibility of $M$ we conclude $\bigcap_{j=1}^m \ker(B_j) =
\{0\}$, and from the surjectivity of $A_j^{-1}$ and
\eqref{pj-ident} we conclude that $B_j$ is surjective. Hence $\B$
is non-degenerate. Finally we compute
\begin{equation}\label{trid}
\tr_H(\id_H) = \sum_{j=1}^m p_j \tr_H(M^{-1} B_j^* A_j B_j ) = \sum_{j=1}^m p_j \tr_{H_j}( A_j B_j M^{-1} B_j^*) = \sum_{j=1}^m p_j \tr_{H_j}(\id_{H_j})
\end{equation}
which is \eqref{scaling}.

Now we verify (d) $\implies$ (c). From the non-degeneracy of $\B$
we see that $M$ is 
positive definite, hence invertible.
From \eqref{pj-ident-alt} we have $\tr_{H_j}(A_j B_j M^{-1} B_j^*)
\leq \tr_{H_j}(\id_{H_j})$.  Using \eqref{scaling} and reversing
the argument in \eqref{trid} we conclude
$$\tr_{H_j}(A_j B_j M^{-1} B_j^*) = \tr_{H_j}(\id_{H_j}),$$
 which together with \eqref{pj-ident-alt}
and the positive definiteness of $A_j$ yields \eqref{pj-ident}.

Now we verify that (d) $\iff$ (e).
If we introduce the operators $T_j: H \to H_j$ defined by
$T_j := A_j^{1/2} B_j M^{-1/2}$ then \eqref{pj-ident-alt} is equivalent to
$T_j T_j^* \leq \id_{H_j}$, while \eqref{pj-ident2} is equivalent to $T_j^* T_j \leq \id_H$.
Since $T_j T_j^*$ and $T_j^* T_j$ have the same operator norm, the claim follows.

Now we verify that (c) $\implies$ (f). Let $\B' := (H, (H_j)_{1
\leq j \leq m}, (T_j)_{1 \leq j \leq m})$ where $T_j$ was defined
earlier, thus $(\B',\p)$ is equivalent to $(\B,\p)$ with
intertwining maps $C = M^{-1/2}$ and $C_j := A_j^{-1/2}$. From the
hypothesis \eqref{pj-ident} we have $T_j T_j^* = \id_{H_j}$, and
from definition of $M$ we have
$$ \sum_{j=1}^m p_j T_j^* T_j = \sum_{j=1}^m M^{-1/2} p_j B_j^* A_j B_j M^{-1/2} = M^{-1/2} M M^{-1/2} = \id_H,$$
and hence $(\B',\p)$ is geometric. The remaining claims in (f)
then follow from Lemma \ref{scaling-lemma}.

Finally, the implication (f) $\implies$ (g) follows from Proposition \ref{gbl-prop}, and the implication (g) $\implies$ (a) is trivial.
\end{proof}

\begin{example}\label{bl-ex-2}  One can now deduce the sharp Young inequality (Example \ref{bl-ex}) from Proposition \ref{normal-form}
by setting $A_j := \frac{1}{p_j(1-p_j)} I_d$ for $j=1,2,3$.  In block matrix notation, one has
$$\sum_{i=1}^3 p_i B_i^* A_i B_i =
\left( \begin{array}{ll}
\frac{p_2}{(1-p_1)(1-p_3)} I_d & \frac{1}{1-p_3}I_d \\
\frac{1}{1-p_3} I_d & \frac{p_1}{(1-p_2)(1-p_3)} I_d
\end{array}\right)$$
and one easily verifies any of the conditions (c)-(e).
\end{example}

From Proposition \ref{normal-form} we see in particular that we
have proven Theorem \ref{lieb-thm} and Theorem \ref{sufbl} in the
case when $(\B,\p)$ is gaussian-extremisable. However, not every
Brascamp--Lieb datum is gaussian-extremisable, even if we assume
the Brascamp--Lieb constants to be finite, as the following example
shows.

\begin{example}\label{gauss-fail}  Consider the rank one case (Example \ref{bbl}) with $H = \R^2$, $m=3$, with any two of $v_1,v_2,v_3$
being linearly independent; this is the case for instance with the one-dimensional version of Young's inequality (Example \ref{bl-ex}).
Then one can verify that the quantity \eqref{BLfunctional} is finite if and only if $(p_1,p_2,p_3)$ lies in
the solid triangle with vertices $(1,1,0), (0,1,1), (1,0,1)$, but if $(p_1,p_2,p_3)$ lies on one of the open edges of this triangle then
no gaussian extremiser will exist.  However, gaussian extremisers do exist on the three vertices of the triangle and on the interior;
this corresponds to the well-known fact that extremisers to Young's inequality $\|f*g\|_r \leq \|f\|_p \|g\|_q$ with $1/r + 1 = 1/p + 1/q$
and $1 \leq p,q,r \leq \infty$ exist when $1 < p,q,r < \infty$, or if all of $p,q,r$ are equal to $1$ or $\infty$, but do not exist
in the remaining cases.  Observe also that this data is simple if and only if $(p_1,p_2,p_3)$ lies in the interior of this triangle.
\end{example}

Fortunately, it turns out that when a Brascamp--Lieb datum is not 
gaussian-extremisable, then it can be factored into lower-dimensional Brascamp--Lieb data, with the corresponding constants also factoring accordingly;
this was first observed by \cite{CLL} in the rank one case.  We will in fact have two means of factoring, which roughly speaking correspond to the notions of direct product $G = H \times K$ and semi-direct product $0 \to N \to G \to G/N \to 0$ in group theory; they also correspond to the notions of \emph{decomposability} and \emph{reducibility} respectively in quiver theory.  We begin with the analogue
of semi-direct product for Brascamp--Lieb data in the next section; the analogue of direct product will be studied in Section \ref{direct-sec}.

\section{Structural theory of Brascamp--Lieb data I : simplicity}

In this section we begin our analysis of the structure of general Brascamp--Lieb data, 
and in particular consider the questions of whether such data
can be placed into a ``normal form'' or decomposed into ``indecomposable'' components.  

Readers familiar with quivers (see e.g.\ \cite{DW} for an introduction) will recognise an $m$-transformation as a special example of a quiver representation, and indeed much
of the structural theory of $m$-transformations which we develop here can be viewed as a crude version of the basic representation theory for
quivers.  It is likely that the deeper theory of such representations is of relevance to this theory, but we do not pursue these connections
here.

We first make a trivial remark, that if one of the exponents $p_j$ in an $m$-exponent ${\bf{p}}$ is zero, then we may omit this exponent,
as well as the corresponding components of an $m$-transformation $\B$, without affecting the 
Brascamp--Lieb constants $\BL(\B,\p)$ or $\BL_\g(\B,\p)$ (or the conditions\eqref{scaling}, 
\eqref{dimension}).  Also, the omission of such exponents does not affect extremisability or gaussian-extremisability (though it does affect any issues regarding uniqueness of extremisers, albeit in a rather trivial way).
We shall use this remark from time to time to reduce to the cases where the exponents $p_j$ are strictly positive.

Next, we establish the (a) $\implies$ (c) and (b) $\implies$ (c) directions of Theorem \ref{sufbl}:

\begin{lemma}[Necessary conditions for finiteness]\label{necc} Let $(\B,\p)$ be Brascamp--Lieb data such that
$\BL(\B,\p)$ or $\BL_\g(\B,\p)$ is finite.  Then we have the scaling condition \eqref{scaling} and
the dimension inequalities \eqref{dimension}.  In particular this implies that $\B$ is non-degenerate.
\end{lemma}

\begin{proof}  By \eqref{blg-trivial} it suffices to verify the claim assuming that $\BL_\g(\B,\p)$ is finite.
Let $\lambda > 0$ be arbitrary.  By applying \eqref{BLfunctional} with the gaussian input $(\lambda \id_{H_j})_{1 \leq j \leq m}$
we see that
$$\BL_\g(\B,\p) \geq \lambda^{\frac{1}{2} [\sum_{j=1}^m p_j \dim(H_j) - \dim(H)]} / \det( \sum_{j=1}^m B_j^* B_j )^{1/2}.$$
Since $\lambda$ is arbitrary, we see that $\BL_\g(\B,\p)$ can only be finite if \eqref{scaling} holds.

Next, let $V$ be any subspace in $H$, and let $0 < \eps < 1$ be a small parameter.  Let $\A^{(\eps)} = (A_j)_{1 \leq j \leq m}$
be the gaussian input $A_j := \eps \id_{B_j V} \oplus \id_{(B_j V)^\perp}$.  Then $\det_{H_j}(A_j)$ decays like $\eps^{\dim(B_j V)}$ as $\eps \to 0$.
Also, we see that $\sum_{j=1}^m B_j^* A_j B_j$ is bounded uniformly in $\eps$, and when restricted to $V$ decays linearly in $\eps$.
Thus $\det_H(\sum_{j=1}^m B_j^* A_j B_j)$ decays at least as fast as $\eps^{\dim(V)}$ as $\eps \to 0$.  By \eqref{BLfunctional(not)}
we conclude that $\BL_\g(\B,\p; \A^{(\eps)})$ grows at least as fast as $\eps^{\frac{1}{2} [\sum_{j=1}^m p_j \dim(B_j V) - \dim(V)]}$.  Since
we are assuming that $\BL_\g(\B,\p)$ is finite, we obtain \eqref{dimension} as desired.

Comparing \eqref{dimension} for $V = H$ with \eqref{scaling} we conclude that the $B_j$ must all be surjective, and applying
\eqref{dimension} with $V = \bigcap_{j=1}^m \ker(B_j)$ we conclude that $\bigcap_{j=1}^m \ker(B_j)$ must equal $\{0\}$.  Thus $\B$ is necessarily
non-degenerate.
\end{proof}

We remark that by testing \eqref{dimension} on $V:= \ker B_j$ we conclude that in order for
$\BL(\B,\p)$ or $\BL_\g(\B,\p)$ to be finite it is necessary that each $p_j \leq 1$.

To proceed further, we need some notation.

\begin{definition}[Restriction and quotient of $m$-transformations] Let 
$\B = (H, (H_j)_{1 \leq j \leq m}, (B_j)_{1 \leq j \leq m})$ be
an $m$-transformation, and let $V$ be a subspace of $H $.  We define the 
\em{restriction} $\B_V$ of $\B$ to $V$ to be the $m$-transformation
$$ \B_V = (V, (B_j V)_{1 \leq j \leq m}, (B_{j,V})_{1 \leq j \leq m} )$$
where $B_{j,V}: V \to B_j V$ is the restriction of $B_j: H \to H_j$ to $V$, and we also define the \emph{quotient} $\B_{H/V}$ of $\B$ 
to be the $m$-transformation
$$ \B_{H/V} = (H/V, (H_j/(B_j V))_{1 \leq j \leq m}, (B_{j,H/V})_{1 \leq j \leq m} )$$
where $B_{j,H/V}: H/V \to H_j/(B_j V)$ is defined by setting $B_{j,H/V}(x+V)$ to be the coset $B_j x + B_j V$.  Equivalently, $\B_V$ and $\B_{H/V}$ are the unique $m$-transformations for which the diagram
$$
\begin{CD}
0       @>>>        V                           @>>>    H                   @>>>    H/V                         @>>>    0       \\
@.                      @VV{B_{j,V}}V               @VV{B_j}V               @VV{B_{j,H/V}}V         @.   \\
0       @>>>        B_j V                   @>>>    H_j             @>>>    H_j/(B_jV)          @>>>    0
\end{CD}
$$
commutes.
\end{definition}

\begin{remark} Using the Hilbert structure, one can of course identify $H/V$ 
with $V^\perp$, the orthogonal complement of $V$, and we will use this 
identification when convenient (particularly for computations).  However 
we prefer the notation $H/V$ as it reinforces the idea that the inner product structure 
is not truly essential to the Brascamp--Lieb analysis.  One can schematically 
represent the above commutative diagram by the ``short exact sequence'' 
$0 \to \B_V \to \B \to \B_{H/V} \to 0$.
\end{remark}

While one can restrict or quotient arbitrary $m$-transformations
using arbitrary subspaces $V$ of $H$, it turns out for the
purposes of analysing Brascamp--Lieb constants that the only
worthwhile subspaces $V$ to restrict to or quotient by are the
\emph{critical} subspaces. Recall that a subspace $V \subset H$ is a
\emph{critical subspace} for a Brascamp--Lieb datum $(\B,\p)$ if 
\begin{equation}\label{critical}
 \dim(V) = \sum_{j=1}^m p_j \dim(B_j V),
\end{equation}
and $V$ is non-zero and proper.
The trivial subspace $\{0\}$ satisfies \eqref{critical},
and the scaling condition \eqref{scaling} is simply the requirement
that $H$ itself satisfy \eqref{critical}.  
We remind the reader that we call a Brascamp--Lieb datum \emph{simple}
if there are no critical subspaces.  Also observe that
equivalent Brascamp--Lieb data have critical
subspaces in 1-1 correspondence; indeed if $V$ is a critical subspace for $(\B,\p)$ and
$(\B',\p)$ is an equivalent Brascamp--Lieb datum with intertwining
maps $C$ and $C_j$, then $C^{-1} V$ is a critical subspace for
$(\B',\p)$.

\begin{example}\label{crit-ex}
For H\"older's inequality (Example \ref{holder-ex}), every non-zero proper subspace of $\R^n$ is critical.
For the Loomis--Whitney inequality (Example \ref{lw-ex}), any co-ordinate plane (spanned by some subset of $\{e_1,\ldots,e_n\}$) is critical.
For Young's inequality (Example \ref{bl-ex}) with $d=1$, there are no critical subspaces when $p_1, p_2, p_3$ lie strictly between 0 and 1,
but if one of $p_1, p_2, p_3$ equals 1 then one of the lines $\{ (x,0): x \in \R\}$, $\{ (0,y): y \in \R \}$ or $\{ (x, -x): x \in  \R \}$ will be
critical.
\end{example}

\begin{remark} Critical subspaces play a role in Brascamp--Lieb data analogous to the 
role normal subgroups play in group theory, as will hopefully
become clearer from the other results in this section.
\end{remark}

A crucial observation for our analysis is
that the necessary conditions \eqref{scaling}, \eqref{dimension} factor through critical subspaces:

\begin{lemma}[Necessary conditions split]\label{critsum} 
Let $(\B,\p)$ be a 
Brascamp--Lieb datum, and let $V_c$ be a critical subspace.  Then
$(\B,\p)$ obeys the conditions \eqref{scaling}, \eqref{dimension} if and only if
both $(\B_{V_c},\p)$ and $(\B_{H/V_c},\p)$ obey the conditions \eqref{scaling}, \eqref{dimension}
for all subspaces $V$ of $V_c$, $H/V_c$ respectively.
\end{lemma}

\begin{proof}  The basic idea is to view $\B$ as the ``semi-direct product'' of $\B_{V_c}$ and $\B_{H/V_c}$, as the schematic diagram
$0 \to \B_{V_c} \to \B \to \B_{H/V_c} \to 0$ already suggests.  A similar idea will also underlie the proof of Lemma \ref{factor} below.

First suppose that $(\B,\p)$ obeys \eqref{scaling}, \eqref{dimension}.  Then it is clear
that $(\B_{V_c},\p)$ obeys \eqref{dimension}, simply by restricting the spaces $V \subseteq \R^n$ to be a subspace of
$V_c$.  Also, the scaling condition \eqref{scaling} for $(\B_{V_c},\p)$ is precisely \eqref{critical}.  To verify
that $(\B_{V_c},\p)$ obeys \eqref{dimension}, we let $V$ be a subspace of $H/V_c$.  Applying \eqref{dimension} to
the subspace $V + V_c$, followed by \eqref{critical}, we observe
\begin{align*}
\dim(V) &= \dim(V + V_c) - \dim(V_c) \\
&\leq \sum_{j=1}^m p_j \dim( B_j(V + V_c) ) - \dim(V_c) \\
&= \sum_{j=1}^m p_j \dim(B_j V + B_j V_c ) - \sum_{j=1}^m p_j \dim(B_j V_c) \\
&= \sum_{j=1}^m p_j \dim(B_{j,H/V_c} V)
\end{align*}
as desired.  A similar computation with $V = H/V_c$ shows that $(\B_{H/V_c},\p)$ also obeys \eqref{scaling}.  This proves the ``only if'' direction
of the lemma.

Now suppose conversely that $(\B_{V_c},\p)$ and $(\B_{H/V_c},\p)$ both obey \eqref{scaling}, 
\eqref{dimension}. Then by adding together the two instances of \eqref{scaling} we see that $(\B,\p)$ also obeys
\eqref{scaling}.  Also, for any $V \subseteq H$, we write $U := V \cap V_c$ and $W := (V + V_c) / V_c$ and compute using
both instances of \eqref{dimension},
\begin{align*}
\dim(V) &= \dim(V \cap V_c) + \dim(V + V_c) - \dim(V_c) \\
&= \dim(U) + \dim(W) \\
&\leq \sum_{j=1}^m p_j \dim(B_{j,V} U) + p_j \dim(B_{j,H/V} W) \\
&= \sum_{j=1}^m p_j \dim(B_j U) + p_j (\dim(B_j W + B_j V_c) - \dim(B_j V_c)) \\
&= \sum_{j=1}^m p_j [\dim(B_j U) + \dim(B_j V + B_j V_c) - \dim(B_j V_c)]
\end{align*}
since $W + V_c = V + V_c$.  Observe that
$$ \dim(B_j U) = \dim(B_j(V \cap V_c)) \leq \dim(B_j V \cap B_j V_c).$$
Since $\dim(B_j V \cap B_j V_c) + \dim(B_j V + B_j V_c) - \dim(B_j V_c) = \dim(B_j V)$,
we see that $(\B,\p)$ obeys \eqref{dimension} as desired.
\end{proof}

Another key fact is the submultiplicativity of  Brascamp--Lieb constants 
through critical subspaces; this was already observed in \cite{CLL} (in the rank one case)
and Finner (in the case of orthogonal projections to co-ordinate spaces).

\begin{lemma}[Submultiplicativity of Brascamp--Lieb constants]\label{subfactor}                
If $(\B,\p)$ is a Brascamp--Lieb datum and $V$ is a critical
subspace with respect to this datum, then
\begin{equation}\label{lblb}
\BL( \B,\p ) \leq \BL( \B_V, \p ) \BL( \B_{H/V}, \p ).
\end{equation}
\end{lemma}

\begin{proof}
We may assume that $\BL( \B_V, \p )$ and $ \BL( \B_{H/V}, \p )$ are finite.
Let $\f = (f_j)_{1 \leq j \leq m}$ be an arbitrary input for $(\B,\p)$.  We may normalise 
$\int_{H_j} f_j = 1$.  We have to show that
$$ \int_{H} \prod_{j=1}^m (f_j \circ B_j)^{p_j} \leq \BL( \B_V, \p ) \BL( \B_{H/V}, \p ).$$
By the Fubini-Tonelli theorem, the left-hand side can be rewritten as
$$ \int_{V^\perp} (\int_V \prod_{j=1}^m (f_j \circ B_j)^{p_j}(v+w)\ dv)\ dw.$$
But we can write $f_j \circ B_j(v+w) = f_{j,w} \circ B_{j,V}(v)$, where $f_{j,w}: B_j V \to \R^+$ is the function
$$ f_{j,w}(v_j) := f_j(v_j + B_j w) \hbox{ for all } v_j \in B_j V.$$
Applying \eqref{BL} we conclude that
$$ \int_V \prod_{j=1}^m (f_j \circ B_j)^{p_j}(v+w)\ dv \leq \BL( \B_V, \p )
\prod_{j=1}^m (\int_{B_j V} f_{j,w})^{p_j}$$
and we are reduced to showing that
\begin{equation}\label{bjv}
 \int_{V^\perp} \prod_{j=1}^m (\int_{B_j V} f_{j,w})^{p_j}\ dw \leq \BL( \B_{H/V}, \p ).
 \end{equation}
We then identify $V^\perp$ with $H/V$ and
observe that $\int_{B_j V} f_{j,w} = f_{j,H/V} \circ B_{j,H/V}(w)$, where $f_{j,H/V}: H_j/(B_j V) \to \R^+$
is the function
$$ f_{j,H/V}(w_j) = \int_{w_j + B_j V} f_j.$$
Applying \eqref{BL} again, we can bound the left-hand side of \eqref{bjv} by
$$ \BL( \B_{H/V}, \p ) \prod_{j=1}^m (\int_{H_j/(B_j V)} f_{j,H/V})^{p_j}$$
and the claim follows from the Fubini-Tonelli theorem and the normalisation
$\int_{H_j} f_j = 1$.  Note the same argument also shows that if
$(\B,\p)$ is extremisable, then $(\B_V,\p)$ and $(\B_{H/V},\p)$ are
also extremisable, because if equality is attained in \eqref{BL}
then all the inequalities above must in fact be equalities.
\end{proof}

This result will be used to prove the finiteness theorem ((c) $\implies$ (a) 
in Theorem \ref{sufbl}) in the next section. 
In fact, the Brascamp--Lieb constants and their gaussian versions are not just 
submultiplicative, but multiplicative:

\begin{lemma}[Brascamp--Lieb constants split]\label{factor}
If $(\B,\p)$ is a Brascamp--Lieb datum and $V$ is a critical
subspace with respect to this datum, then
$$ \BL( \B,\p ) = \BL( \B_V, \p ) \BL( \B_{H/V}, \p )$$
and
$$ \BL_\g( \B,\p ) = \BL_\g( \B_V, \p ) \BL_\g( \B_{H/V}, \p ).$$
Moreover, if $(\B,\p)$ is extremisable, then $(\B_V,\p)$ and $(\B_{H/V},\p)$ are also extremisable.
Similarly, if $(\B,\p)$ is gaussian-extremisable, then $(\B_V,\p)$ and $(\B_{H/V},\p)$ are 
also gaussian-extremisable.
\end{lemma}

\begin{proof}  If $(\B,\p)$ fails to obey \eqref{scaling} or \eqref{dimension}, then Lemma \ref{necc} and Lemma \ref{critsum} show that
all expressions in this lemma are infinite and so the conclusion trivially holds. Thus we may assume 
that $(\B,\p)$ obeys \eqref{scaling},
\eqref{dimension}.  By Lemma \ref{critsum} we see that $(\B_V, \p)$ and $(\B_{H/V}, \p)$ also obey 
\eqref{scaling}, \eqref{dimension}.
In particular all of these data are non-degenerate.

Having proved \eqref{lblb} above, 
let us first establish the reverse inequality
\begin{equation}\label{blbl}
\BL( \B,\p ) \geq \BL( \B_V, \p ) \BL( \B_{H/V}, \p ).
\end{equation}
Let $0 < C_V < \BL( \B_V, \p ) $ and $0 < C_{H/V} < \BL( \B_{H/V}, \p )$ be arbitrary
constants.  Then by definition of Brascamp--Lieb constant and homogeneity, we can find inputs $\f_V = (f_{j,V})_{1 \leq j \leq m}$
and $\f_{H/V} = (f_{j,H/V})_{1 \leq j \leq m}$ for $(\B_V,\p)$ and $(\B_{H/V},\p)$ respectively,
such that
$$ \BL(\B_V, \p; \f_V) > C_V; \quad \BL(\B_{H/V}, \p; f_{H/V} ) > C_{H/V}.$$
We may normalise
$$ \int_{B_j V} f_{j,V} = \int_{H_j/(B_j V)} f_{j,H/V} = 1 \hbox{ for all } j.$$
Let $\lambda > 1$ be a large parameter.  We define the input $\f^{(\lambda)} = (f_j)_{1 \leq j \leq m}$ for $(\B,\p)$ by the formula
$$ f_j( v + w ) := f_{j,V}( v_j ) f_{j,H/V}(\lambda w_j) \hbox{ whenever } v_j \in B_j V; w_j \in (B_j V)^\perp,$$
where we identify $H_j/(B_j V)$ with $(B_j V)^\perp$ in the usual manner.
Then by the Fubini-Tonelli theorem and the normalisation of $\f_V, \f_{H/V}$, we observe that
\begin{equation}\label{fj-int}
\int_{H_j} f_j = \lambda^{\dim(B_j V) - n_j}.
\end{equation}
Now we investigate the expression $\prod_{j = 1}^m f_j \circ B_j(v+w)$ where $v \in V$ and $w \in V^\perp$, where we identify $V^\perp$
with $H/V$ in the usual manner.
Observe for any $v \in V$, $w \in V^\perp$ that
$$ \prod_{j = 1}^m f_j \circ B_j(v + \lambda^{-1} w) =
\prod_{j = 1}^m f_{j,V}( B_{j,V} v + \lambda^{-1} \pi_{B_j V}(B_j w)) f_{j,H/V}(B_{j,H/V} w)$$
and thus by the Fubini-Tonelli theorem and rescaling
\begin{eqnarray*}
\begin{aligned}
& \hspace{25mm} \lambda^{\dim(V)-\dim(H)} \int_H \prod_{j = 1}^m f_j \circ B_j \\ 
& =\int_{V^\perp} \int_V \prod_{j = 1}^m f_{j,V}( B_{j,V} v 
+ \lambda^{-1} \pi_{B_j V}(B_j w))\ dv 
\prod_{j=1}^m f_{j,H/V}(B_{j,H/V} w)\ dw.
\end{aligned}
\end{eqnarray*}
From this, \eqref{fj-int}, \eqref{BL} and \eqref{critical} we conclude
$$ \BL( \B, \p ) \geq
\int_{V^\perp} (\int_V \prod_{j = 1}^m f_{j,V}( B_{j,V} v + \lambda^{-1} \pi_{B_j V}(B_j w))\ dv )
\prod_{j=1}^m f_{j,H/V}(B_{j,H/V} w)\ dw.$$
Since $\B_V$ and $\B_{H/V}$ are non-degenerate, we have
$$\bigcap_{j=1}^m \ker(B_{j,H/V}) = \bigcap_{j=1}^m \ker(B_{j,V}) = \{0\}.$$  
Thus in the above
integrals, $v$ and $w$ range over a compact set uniformly in $\lambda > 1$.
We may thus take limits as $\lambda \to \infty$ (using the smoothness of the $f_{j,V}$) to conclude
$$ \BL( \B,\p ) \geq
\int_{V^\perp}(\int_V \prod_{j = 1}^m f_{j,V}( B_{j,V} v )dv)
\prod_{j=1}^m f_{j,V^\perp}(B_{j,H/V} w)\ dw.$$ By construction of
the inputs $\f_V, \f_{H/V}$ we thus have
$$ \BL( \B,\p ) \geq C_V C_{H/V}$$
and upon taking suprema in $C_V$, $C_{H/V}$ we obtain \eqref{blbl} as desired.

We observe that the corresponding inequality for
the gaussian Brascamp--Lieb constants follows by an identical argument, with the $f_{j,V}$, $f_{j,H/V}$ (and hence $f_j$)
now being centred gaussians instead of test functions.  Note that one can still justify the limit 
as $\lambda \to \infty$ using dominated convergence. Thus one has 
\begin{eqnarray}\label{blblgaussian}
\BL_\g( \B,\p ) \geq \BL_\g( \B_V, \p ) \BL_\g( \B_{H/V}, \p ).
\end{eqnarray}

Turning now to the analogue of \eqref{lblb} for the gaussian case,
we first observe that the argument 
given above does not quite work in the 
gaussian case because if $f_j$ is a centred gaussian then $f_{j,w}$ is 
merely an uncentred gaussian except when $w=0$. For this reason it is 
pertinent to introduce the \emph{uncentred} gaussian Brascamp--Lieb 
constant $\BL_\u( \B,\p )$ defined as the best constant in \eqref{BL} when the input 
is restricted to consist of uncentred gaussians. 
The argument now shows that the analogue of
\eqref{lblb} holds for $\BL_\u( \B,\p )$. The proof of the inequality corresponding 
to \eqref{blbl} for $\BL_\u( \B,\p )$ follows by the identical argument, with the 
$f_{j,V}$, $f_{j,H/V}$ (and hence $f_j$) now being uncentred gaussians instead of test 
functions. Note that one can still justify the limit as $\lambda \to \infty$ using dominated convergence.
Thus the constants $\BL_\u( \B,\p )$ are multiplicative; that is 
\begin{eqnarray}\label{lmlm}
\BL_\u( \B,\p ) = \BL_\u( \B_V,\p ) \BL_\u( \B_{H/V},\p ).
\end{eqnarray}

It therefore suffices to show that the centred and uncentred gaussian constants coincide, that is 
$\BL_\g( \B,\p ) = \BL_\u( \B,\p )$. Clearly $\BL_\g( \B,\p ) \leq  \BL_\u( \B,\p )$ 
so it is enough to show $\BL_\g( \B,\p ) \geq \BL_\u( \B,\p ).$ While this is an obvious 
consequence of Lieb's theorem, an elementary argument is available. 

Consider an $m$-tuple of gaussians $f_{j}(x):=\exp(-\pi\langle 
A_{j}(x-\xi_{j}),(x-\xi_{j})\rangle)$
centred at $\xi_{j}\in\R^{n_{j}}$. Since $B_{j}$ is surjective we may take 
$\xi_{j}=B_{j}w_{j}$, 
for some $w_{j}\in\R^{n}$. Let $\overline{w}=
M^{-1}\sum_{j=1}^{m}p_{j}B_{j}^{*}A_{j}B_{j}w_{j}$, where
$M=\sum_{j=1}^{m}p_{j}B_{j}^{*}A_{j}B_{j}$, so that
$\sum_{j=1}^{m}p_{j}B_{j}^{*}A_{j}B_{j}(w_{j}-\overline{w})=0.$
Now by a simple change of variables,
\begin{eqnarray*}
\begin{aligned}
\int\prod f_{j}(B_{j}(x))^{p_{j}}\ dx&=
\int\exp\Bigl\{
-\pi\sum_{j=1}^{m}p_{j}\langle A_{j}B_{j}(x-w_{j}),B_{j}(x-w_{j})\rangle\Bigr\}
\ dx\\
=\int\exp\Bigl\{
-\pi&\sum_{j=1}^{m}p_{j}\langle A_{j}B_{j}(x-(w_{j}-\overline{w})),
B_{j}(x-(w_{j}-\overline{w}))\rangle\Bigr\}
\ dx\\
=\exp\Bigl\{-\pi\sum_{j=1}^{m}p_{j}\langle A_{j}&B_{j}(w_{j}-\overline{w}),
B_{j}(w_{j}-\overline{w})\rangle\Bigr\}\int\exp\Bigl\{
-\pi\sum_{j=1}^{m}p_{j}\langle A_{j}B_{j}x,
B_{j}x\rangle\Bigr\}
\ dx\\
&\leq\int\exp\Bigl\{
-\pi\sum_{j=1}^{m}p_{j}\langle A_{j}B_{j}x,
B_{j}x\rangle\Bigr\}
\ dx.
\end{aligned}
\end{eqnarray*}
Hence $\int\prod f_{j}(B_{j}(x))^{p_{j}}\ dx$ is maximised when all of the $\xi_{j}$'s are
zero, and thus $\BL_{\u}(\B,\p)\leq \BL_{\g}(\B,\p)$.

The assertion concerning extremisability of $(\B, \p)$ is made at the end of the proof 
of Lemma \ref{subfactor} above. Similarly, if $(\B, \p)$ is gaussian-extremisable, with 
$\A$ a gaussian input extremising \eqref{BLfunctional(not)}, then with $f_j(x) =
(\det A_j)^{1/2} \exp(-\pi \langle A_j x, x \rangle)$ we can follow the proof of \eqref{lblb},
to see that for all $w$ the uncentred gaussians $(f_{j,w})_{1 \leq j \leq m}$ attain $\BL_\u (\B_V, \p) 
= \BL_\g (\B_V, \p)$, and the centred gaussians $(f_{j,H/V})_{1 \leq j \leq m}$ 
attain $\BL_\g(\B_{H/V},\p)$. In particular, taking $w=0$, both $(\B_V, \p)$ and
$(\B_{H/V}, \p)$ are gaussian-extremisable.
\end{proof}

\begin{remark}\label{bar}
We note that the above proof of $\BL_\g( \B,\p ) \geq \BL_\u( \B,\p )$ 
shows that if a collection of uncentred gaussians are extremisers, then they are obtained 
from the corresponding centred gaussians by translation by $B_j \overline{w}$ for a fixed $
\overline{w} \in H$, and that any such translations also furnish extremisers.
This will complement Theorem \ref{extreme-unique}.
\end{remark}

\begin{remark}
Recall that a Brascamp--Lieb datum is \emph{simple} if it has no critical subspaces.  Lemma 
\ref{critsum} and Lemma \ref{factor} allow us
to reduce the problem of finiteness of Brascamp--Lieb constants (or of proving Lieb's theorem) 
to the case of simple Brascamp--Lieb data.  As we shall
see later, simple data obeying the necessary conditions \eqref{scaling}, \eqref{dimension}
are always gaussian-extremisable and hence equivalent to geometric Brascamp--Lieb data, and 
thus Lieb's theorem reduces entirely to checking the geometric Brascamp--Lieb case, which was already done in Proposition \ref{gbl-prop}.
\end{remark}

\begin{remark}
The question of whether gaussian extremisers exist in general is a little bit more complicated; the correct characterisation is not simplicity but \emph{semi-simplicity} (the direct sum of simple data).  We shall study this issue in depth in Section \ref{direct-sec}.
\end{remark}

\section{Sufficient conditions for finiteness, and Lieb's theorem}\label{necsuf-sec}

We are now ready to prove Theorem \ref{sufbl} and the general case on Theorem \ref{lieb-thm}.  We begin with a lemma.

\begin{lemma}\label{babygauss} Let $(\B,\p)$ be a Brascamp--Lieb datum such that
\begin{equation}\label{dimh}
 \dim(H/V) \geq \sum_{j=1}^m p_j \dim(H_j / (B_j V))
 \end{equation}
for all subspaces $V$ of $H$.  (In particular, \eqref{dimh} is equivalent to \eqref{dimension} if \eqref{scaling} holds.)
Then there exists a real number $c > 0$, such that for every orthonormal basis $e_1,\ldots,e_n$ of $H$ there exists a set $I_j \subseteq \{1,\ldots,n\}$ for each $1 \leq j \leq m$ with $|I_j| = \dim(H_j)$ such that
\begin{equation}\label{p-comb-alt-2}
\sum_{j=1}^m p_j |I_j \cap \{1,\ldots,k\}| \leq k \hbox{ for all } 0 \leq k \leq n
\end{equation}
and
\begin{equation}\label{bjei}
\| \bigwedge_{i \in I_j} B_j e_i \|_H \geq c \hbox{ for all } 1 \leq j \leq m.
\end{equation}
Here $\bigwedge_{i \in I_j} B_j e_i$ denotes the wedge product of the vectors $B_j e_i$, and the $\| \|_H$ norm is the usual norm on forms induced from the Hilbert space structure.
In particular, if \eqref{scaling} holds, we have
\begin{equation}\label{p-comb-alt}
\sum_{j=1}^m p_j |I_j \cap \{k+1,\ldots,n\}| \geq n-k \hbox{ for all } 0 \leq k \leq n.
\end{equation}
If furthermore there are no critical subspaces, then we can enforce strict inequality in \eqref{p-comb-alt} for all $0 < k < n$.
\end{lemma}

\begin{proof}  If \eqref{scaling} holds then $n = \sum_{j=1}^m p_j \dim(H_j)$, and \eqref{p-comb-alt} follows from \eqref{p-comb-alt-2} and
the fact that $|I_j| = \dim(H_j)$.  Thus we shall focus on proving
the claim \eqref{p-comb-alt-2}. The case $k=0$ of
\eqref{p-comb-alt-2} is trivial, and the $k=n$ case follows from
\eqref{dimh}, so we restrict attention to $0 < k < n$.  From the
hypotheses we know that the $B_j$ are surjective.

The space of all orthonormal bases is compact, and the number of possible $I_j$ is finite.  Thus by continuity and compactness we may replace
the conclusion \eqref{bjei} by the weaker statement
$$  \bigwedge_{i \in I_j} B_j e_i  \neq 0 \hbox{ for all } 1 \leq j \leq m.$$
In other words, we require the vectors $(B_j e_i)_{i \in I_j}$ to be linearly independent on $H_j$ for each $j$.

We shall select the $I_j$ by a backwards greedy algorithm.  Namely, we set $I_j$ equal to those indices $i$ for which
$B_j e_i$ is not in the linear span of $\{ B_j e_{i'}: i < i' \leq n\}$ (thus for instance $n$ will lie in $I_j$ as long as
$B_j e_n \neq 0$).  Since the $B_j$ are surjective, we see that $|I_j| = \dim(H_j)$.  To prove \eqref{p-comb-alt-2}, we apply the hypothesis
\eqref{dimh} with $V$ equal to the span of $\{e_{k+1},\ldots,e_n\}$, to obtain
$$ \sum_j p_j \dim(H_j / B_j V) \leq k.$$
But by construction of $I_j$ we see that $\dim(B_j V) = |I_j \cap \{k+1,\ldots,n\}|$
and hence $\dim(H_j / B_j V) = |I_j \cap \{1,\ldots,k\}|$.  The claim \eqref{p-comb-alt-2} follows.

If \eqref{scaling} holds and there are no critical spaces, then one always has strict inequality in \eqref{dimension} or \eqref{dimh},
hence in \eqref{p-comb-alt-2} or \eqref{p-comb-alt} when $0 < k < n$.
\end{proof}

We can now prove the sufficiency of \eqref{scaling} and \eqref{dimension} in Theorem \ref{sufbl} for the finiteness of the gaussian Brascamp--Lieb constant.  

\begin{proposition}\label{sufbl-meat} Let $(\B,\p)$ be a Brascamp--Lieb datum such 
that \eqref{scaling} and \eqref{dimension} hold.  Then
$\BL_\g(\B,\p)$ is finite.  Furthermore, if $(\B,\p)$ is simple 
(i.e., there is no critical subspace), then $(\B,\p)$ is gaussian-extremisable.
\end{proposition}

We remark that this proposition was established in the rank one case by \cite{CLL}.

\begin{proof} We may discard those $j$ for which $p_j = 0$ or for which $H_j = \{0\}$, as these factors clearly give no contribution
to the Brascamp--Lieb constants (or to \eqref{scaling} or \eqref{dimension}), nor does they affect whether $(\B,\p)$ is simple or not.
In particular $\B$ will still be non-degenerate after doing this.

Fix a gaussian input $\A = (A_j)_{1 \leq j \leq m}$, and let $M := \sum_j p_j B_j^* A_j B_j$.  This transformation is self-adjoint; since $\B$ is non-degenerate and $p_j > 0$, we also see
that it is 
positive definite. Thus by choosing an appropriate orthonormal basis
$\{e_1,\ldots,e_n\} \subset H$ we may assume that $M = \diag(\lambda_1, \ldots, \lambda_n)$ for some $\lambda_1 \geq \ldots \geq \lambda_n > 0$.

Applying Lemma \ref{babygauss}, we can find $I_j \subseteq \{1,\ldots,n\}$ for each $1 \leq j \leq m$
of cardinality $|I_j| = \dim(H_j)$
obeying \eqref{p-comb-alt} and \eqref{bjei}.  For each $i \in I_j$, we have
$$ \langle A_j B_j e_i, B_j e_i \rangle_{H_j} = \langle e_i, B_j^* A_j B_j e_i \rangle_{H}
\leq \frac{1}{p_j} \langle e_i, Me_i \rangle_{H} = \lambda_i / p_j.$$
On the other hand, from \eqref{bjei} we see that $(B_j e_i)_{i \in I_j}$ is a basis of $H_j$ with a lower bound on the degeneracy.
We thus conclude that
$$ \det(A_j) \leq C \prod_{i \in I_j} \lambda_i$$
for some constant $C > 0$ depending on the Brascamp--Lieb data.  Thus
$$ \prod_{j=1}^m (\det A_j)^{p_j} \leq C \prod_{i=1}^n \lambda_i^{\sum_{j=1}^m p_j |I_j \cap \{i\}|}.$$
We can telescope the right-hand side (using \eqref{scaling}) and obtain
$$ \prod_{j=1}^m (\det A_j)^{p_j} \leq C \lambda_1^n \prod_{0 \leq k \leq {n-1}} (\lambda_{k+1}/\lambda_k)^{\sum_{j=1}^m p_j |I_j \cap \{k+1,\ldots,n\}|}$$
where we adopt the convention $\lambda_0 = \lambda_1$.  Applying \eqref{p-comb-alt} we conclude
$$ \prod_{j=1}^m (\det A_j)^{p_j} \leq C \lambda_1^n \prod_{0 \leq k \leq {n-1}} (\lambda_{k+1}/\lambda_k)^{n-k}$$
which by reversing the telescoping becomes
$$ \prod_{j=1}^m (\det A_j)^{p_j} \leq C \lambda_1 \ldots \lambda_k = C \det(M).$$
Comparing this with \eqref{BLfunctional(not)} we conclude that $\BL_\g(\B,\p)$ is finite.

Now suppose that $(\B,\p)$ is simple.  Then we have strict inequality in \eqref{p-comb-alt}.  We may thus refine the above analysis
and conclude that
$$  \prod_{j=1}^m (\det A_j)^{p_j} \leq C \det(M) \prod_{1 \leq k \leq {n-1}} (\lambda_{k+1}/\lambda_k)^c$$
for some $c > 0$ depending on the Brascamp--Lieb data.  This shows that the expression in \eqref{BLfunctional(not)} goes to zero whenever
$\lambda_n / \lambda_1$ goes to zero.  Thus to evaluate the supremum it suffices to do so in the region $\lambda_1 \leq C \lambda_n$.
Also using the scaling hypothesis \eqref{scaling} we may normalise $\lambda_n = 1$.  This means that $M$ is now bounded above
and below, which by surjectivity of $B_j$ implies that $A_j$ is also bounded.  We may now also assume that $A_j$ is bounded from below since
otherwise the expression \eqref{BLfunctional(not)} 
in the supremum in \eqref{BLfunctional} will be small.  We have thus localised each the $A_j$ to a compact set, and hence by continuity we see that an extremiser exists.  Thus $(\B,\p)$ is gaussian-extremisable as desired.
\end{proof}

We can now prove Theorem \ref{lieb-thm} and Theorem \ref{sufbl} simultaneously and quickly.

\begin{proof}[of Theorem \ref{lieb-thm} and Theorem \ref{sufbl}]  As in \cite{CLL}, we induct on the 
dimension $\dim(H)$.  When $\dim(H) = 0$ the claim is trivial.
Now suppose inductively that $\dim(H) > 0$ and the claim has already been proven for smaller values of $\dim(H)$.  In light of
Lemma \ref{necc} we may assume that \eqref{scaling} and \eqref{dimension} hold.  From Proposition \ref{sufbl-meat} we know that $\BL_\g(\B,\p)$ is
finite; our task is to show that $\BL(\B,\p)$ is equal to $\BL_\g(\B,\p)$.

We divide into two cases.
First suppose that $(\B,\p)$ is simple.  Then
by Proposition \ref{sufbl-meat} the datum $(\B,\p)$ is gaussian-extremisable, and the claim then follows from
Proposition \ref{normal-form}.  Now consider the case when $(\B,\p)$ is not simple, i.e., there is a 
proper critical subspace $V$. Then we can split the Brascamp--Lieb data $(\B,\p)$ into $(\B_V, \p)$ and $(\B_{H/V}, \p)$.  By Lemma \ref{critsum}
and the induction hypothesis we see that
$$ \BL( \B_V,\p ) = \BL_\g( \B_V, \p ) < \infty$$
and
$$ \BL( \B_{H/V},\p ) = \BL_\g( \B_{H/V},\p ) < \infty$$
Applying Lemma \ref{factor} we conclude
$$ \BL( \B, \p ) = \BL_\g( \B, \p ) < \infty$$
thus closing the induction.
\end{proof}

Let us now fix the $m$-transformation $\B$, and let the $m$-exponent $\p = (p_1,\ldots,p_m)$ vary.  Let us define the following subsets of $\R^m$:
\begin{align*}
S(\B) &:= \{ \p \in \R^m: p_j \geq 0 \hbox{ for all } j; \quad \sum_{j=1}^m p_j \dim(H_j) = \dim(H) \}\\
\Pi(\B) &:= \{ \p \in \R^m: \sum_{j=1}^m p_j \dim(B_j V) \geq \dim(V) \hbox{ for all } \{0\} \subsetneq V \subseteq H \} \\
\Pi^\circ(\B) &:= \{ \p \in \R^m: \sum_{j=1}^m p_j \dim(B_j V) > \dim(V) \hbox{ for all } \{0\} \subsetneq V \subseteq H \}
\end{align*}

Observe that even though there are an infinite number of vector spaces $V$, there are only a finite number of possible values for the
dimensions $\dim(V), \dim(B_j V)$.  Thus $\Pi(\B)$ is a closed convex cone with only finitely many faces, and $\Pi^\circ(\B)$ is the $m$-dimensional interior of that
cone.

Theorem \ref{sufbl} thus asserts that the Brascamp--Lieb constant $\BL(\B,\p) = \BL_\g(\B,\p)$ is finite
if and only if $\p$ lies in $S(\B) \cap \Pi(\B)$, and furthermore if $\p$ lies in $S(\B) \cap \Pi(\B)^\circ$ then $(\B,\p)$ is also
gaussian-extremisable.  This of course implies that equality in \eqref{BL} can also be attained.

\begin{remark} It is perfectly possible for $S(\B) \cap \Pi(\B)$ or $S(\B) \cap \Pi^\circ(\B)$ to be empty.  For instance for the Loomis--Whitney
inequality (Example \ref{lw-ex}), $S(\B) \cap \Pi(\B)$ consists of a single point $(\frac{1}{n-1},\ldots, \frac{1}{n-1})$, and $S(\B) \cap \Pi(\B)^\circ$
is empty.  If $\B$ is degenerate, then $S(\B) \cap \Pi(\B)$ is empty.
\end{remark}

An alternative proof that 
the hypotheses \eqref{scaling} and \eqref{dimension} characterise
the Brascamp--Lieb constant $\BL(\B,\p)$
is given in the companion paper \cite{BCCT}.
That proof does not involve 
gaussians, extremisers, or monotonicity formulae,
but instead uses multilinear interpolation, H\"older's inequality,
and an induction argument on dimension based on factorisation through critical subspaces.
Such an induction and factorisation argument had previously appeared in \cite{CLL} (for the rank one case)
and \cite{Finner} (for the case of orthogonal projections to co-ordinate spaces).

We now connect the above results to the work of Barthe \cite{Barthe}, who considered 
the rank-one case (Example \ref{bbl}). This connection was also noted in \cite{CLL}, though our treatment
here differs slightly from that in \cite{CLL}.
We are assuming $\B$ to be non-degenerate, thus the $v_j$ are non-zero and
span $H$.  The condition \eqref{scaling} in this case simplifies to $\sum_{j=1}^m p_j = n$, while \eqref{dimension} in this case
simplifies to
$$ \dim(V) \leq \sum_{1 \leq j \leq m: v_j \not \in V^\perp} p_j.$$
Subtracting this from the scaling identity $\sum_{j=1}^m p_j = n$, we see that \eqref{dimension} has now become
$$ \sum_{1 \leq j \leq m: v_j \in V^\perp} p_j \leq \dim(V^\perp).$$
Since we can set $V^\perp$ to be any subspace, and in particular to be the span of any set of vectors $v_j$, we thus easily see that the constraints \eqref{dimension} are then (assuming \eqref{scaling}) equivalent to the assertion
\begin{equation}\label{isum}
\sum_{j \in I} p_j \leq d_I \hbox{ for all } I \subseteq \{1,\ldots,n\},
\end{equation}
where $d_I$ is the integer $d_I := \dim \linearspan((v_j)_{j \in I})$.  In particular we 
have $0 \leq p_j \leq 1$ for all $j$, (which, as we have noted above, is in general a
necessary condition). 
Thus in the rank one case, we have characterised the polytope where the Brascamp--Lieb
constant is finite as
$ S(\B) \cap \Pi(\B) $ $$=\{ (p_1,\ldots,p_m) \in \R_+^m: \sum_{1 \leq j \leq m} p_j = n; 
\sum_{j \in I} p_j \leq d_I \hbox{ for all } I \subseteq \{1,\ldots,m\}\}.$$
Let us now characterise the extreme points of this polytope.

\begin{lemma}\cite{CLL}  Suppose $H_j= \R$ for all $j$.
Let $(p_1,\ldots,p_m)$ be an extreme point of $S \cap \Pi$.  Then all of the $p_j$ are equal to 0 or 1.
\end{lemma}

\begin{proof}  We induct on $n+m$.  When $n+m=0$ the claim is trivial, so suppose that $n+m > 0$ and the claim has been proven for all smaller values of $n+m$.  If one of the $p_j$ is already equal to zero, say $p_m = 0$, then we can remove that index $m$ (and the associated vector $v_m$)
from $\B$ to form an $(m-1)$-transformation $\tilde \B$, and
observe that $(p_1,\ldots,p_{m-1})$ is an extreme point of 
$S(\tilde \B) \cap \Pi(\tilde \B)$.  Thus the claim follows
from the induction hypothesis.  Now suppose that one of the 
$p_j$ is equal to one, say $p_m = 1$.  Let $H' := H / \R v_m$ be the quotient of $H$ by
$v_m$, and let $v'_1,\ldots,v'_{m-1}$ be the image of 
$v_1,\ldots,v_m$ under this quotient map.  We let $\B'$ be the associated $(m-1)$-transformation.
If we then let $d'_I := \dim \linearspan( (v'_j)_{j \in I} )$ 
for all $I \subseteq \{1,\ldots,m-1\}$, then we have
\begin{align*}
S(\B') \cap \Pi(\B') &= \{ (q_1,\ldots,q_{m-1}) \in \R_+^{m-1}: 
\sum_{1 \leq j \leq m-1} q_j = n-1; \\
&\quad \sum_{j \in I} q_j \leq d'_I \hbox{ for all } I \subseteq \{1,\ldots,m\}\}.
\end{align*}
Now
$$ \sum_{j \in I} p_j = \sum_{j \in I \cup \{m\}} p_j - 1 \leq d_{I \cup \{m\}} - 1 = d'_I.$$
Thus $(p_1,\ldots,p_{m-1}) \in S(\B') \cap \Pi(\B')$.  Conversely, if $(q_1,\ldots,q_{m-1}) \in S(\B') \cap \Pi(\B')$ then reversing the above argument
shows that $(q_1,\ldots,q_{m-1},1) \in S(\B) \cap \Pi(\B)$.  Since $(p_1,\ldots,p_{m-1},1)$ was an extreme point of $S(\B) \cap \Pi(\B)$, this implies
that $(p_1,\ldots,p_{m-1})$ is an extreme point of $S(\B') \cap \Pi(\B')$.  Thus we can again apply the induction hypothesis to close the argument.

Finally, suppose that none of the $p_j$ are equal to 0 or 1.  Call a proper subset $\emptyset \subsetneq I \subseteq \{1,\ldots,m\}$
of $\{1,\ldots,m\}$ \emph{critical} if
$\sum_{j \in I} p_j = d_I$.  Since $d_I$ has to be an integer, we see that no singleton sets are critical.  On the other hand, the entire set $\{1,\ldots,m\}$ is critical.  Thus if we let $I_{min}$ be a non-empty critical set of minimal size, then $I_{min}$ has at least two elements.  Now suppose that $I$ is another critical set which intersects $I_{min}$.  Then
\begin{align*}
d_{I_{min}} + d_I &= \sum_{j \in I_{min}} p_j + \sum_{j \in I} p_j \\
&= \sum_{j \in I_{min} \cap I} p_j + \sum_{j \in I_{min} \cup I} p_j \\
&\leq d_{I_{min} \cap I} + d_{I_{min} \cup I} \\
&\leq \dim[\linearspan( (v_j)_{j \in I_{min}} ) \cap \linearspan( (v_j)_{j \in I} )] \\
&\quad\quad + \dim[\linearspan( (v_j)_{j \in I_{min}} ) + \linearspan( (v_j)_{j \in I} )] \\
&= d_{I_{min}} + d_I
\end{align*}
and hence all the above inequalities must in fact be equality.  In particular this implies that the non-empty set $I_{min} \cap I$
is critical, which by minimality of $I_{min}$ forces $I_{min} \subseteq I$.  Thus all the critical sets either contain $I_{min}$ or
are disjoint from it.  Now recall that $I_{min}$ has at least two elements, and that $p_j$ lies strictly between 0 and 1 for all $j$ in $I_{min}$.  Thus
one can reduce one of the $p_j$, $j \in I_{min}$ by an epsilon and increase another $p_j$, $j \in I_{min}$ by the same epsilon, and stay in
$S(\B) \cap \Pi(\B)$, or conversely.  This shows that $(p_1,\ldots,p_m)$ is not an extreme point of $S(\B) \cap \Pi(\B)$, a contradiction.
\end{proof}

Note that if $(p_1,\ldots,p_m) \in S(\B) \cap \Pi(\B)$ consists entirely of 
$0$'s and $1$'s, then the set $\{ v_j: p_j = 1 \}$ is a basis for $H$ (because
it has cardinality $n$ by \eqref{scaling}, and any subset of this set of cardinality $k$ must span a space of dimension at least $k$
by \eqref{isum}).  Conversely, if $I \subseteq \{1,\ldots,m\}$ is such that $\{v_j: j \in I\}$ is a basis for $H$, then
the $m$-tuple $(1_{j \in I})_{1 \leq j\leq m}$ is easily verified to lie in $S(\B) \cap \Pi(\B)$.  If one takes the convex hull of all these points, then any one of these points $(1_{j \in I})_{1 \leq j\leq m}$ will be vertices of these convex hull, as can be seen by maximising the linear functional $\sum_{j \in I} p_j$ on this convex hull.  Combining all these facts, we have reproved the following theorem of Barthe.

\begin{theorem}\label{barthe-thm}\cite[Section 2]{Barthe}  Let $v_1,\ldots,v_m$ be 
non-zero vectors which span a Euclidean space $\R^n$, and let $B_j: \R^n \to \R$ be 
the maps $B_j x := \langle x, v_j \rangle$.  Let $P \subset \R_+^m$ be the set of 
all exponents $\p$ for which the Brascamp--Lieb constant $\BL(\B,\p)$ is finite.  
Then $P$ is a convex polytope, whose vertices are precisely those points
$(1_{j \in I})_{1 \leq j \leq m})$ where $(v_j)_{j \in I}$ is a basis for $\R^n$.  
If $\p$ lies in the $(m-1)$-dimensional
interior of this polytope, $(\B,\p)$ is gaussian-extremisable.
\end{theorem}

\begin{remark} Barthe's proof of this theorem relies upon an analysis of the constant
\eqref{BLfunctional}. The $\sup$ in \eqref{BLfunctional} is now over $(0,\infty)^m$ and the
Cauchy-Binet theorem gives an explicit formula for the denominator in \eqref{BLfunctional(not)}.
This facilitates a direct analysis of \eqref{BLfunctional} in the rank-one case.
\end{remark}

\section{extremisers}

In this section we analyse extremisers to the Brascamp--Lieb inequality, and connect these extremisers to the heat flow approach used
previously.

We now recall an important observation of K. Ball (which can be
found for instance in \cite{Barthe-old}) which can be used to
motivate the monotonicity formula approach.  If $\f = (f_j)_{1
\leq j \leq m}$ and $\f' = (f'_j)_{1 \leq j \leq m}$ are inputs
for a Brascamp--Lieb datum $(\B,\p)$, we define their
\emph{convolution} $\f*\f'$ to be the $m$-tuple $\f * \f' := (f_j
* f'_j)_{1 \leq j \leq m}$, where of course $f_j*g_j(x) :=
\int_{H_j} f_j(x-y) g_j(y)\ dy$. Observe from the Fubini-Tonelli theorem
that $\f*\f'$ will also be an input for $(\B,\p)$.

\begin{lemma}[Convolution inequality]\label{conv}\cite{Barthe-old}
Let $(\B,\p)$ be a Brascamp--Lieb datum with $\BL(\B,\p)$ finite, and let $\f$, $\f'$ be inputs for $(\B,\p)$.  Then
$$ \BL(\B,\p; \f*\f') \geq \frac{\BL(\B,\p;\f) \BL(\B,\p; \f')}{\BL(\B,\p)}.$$
\end{lemma}

\begin{proof} Write $\f = (f_j)_{1 \leq j\leq m}$ and $\f' = (f'_j)_{1 \leq j \leq m}$
After normalisation we can then assume $\int_{H_j} f_j = \int_{H_j} f'_j = 1$.  Thus
$$ \int_H \prod_{j=1}^m (f_j \circ B_j)^{p_j} = \BL(\B,\p;\f); \quad \int_H \prod_{j=1}^m (f'_j \circ B_j)^{p_j} = \BL(\B,\p;\f').$$
Convolving these together we obtain
$$ \int_H \prod_{j=1}^m (f_j \circ B_j)^{p_j} * \prod_{j=1}^m (f'_j \circ B_j)^{p_j} = \BL(\B,\p;\f) \BL(\B,\p;\f').$$
The left-hand side can be rewritten using the Fubini-Tonelli theorem as
$$ \int_H (\int_H \prod_{j=1}^m (g_j^x \circ B_j)^{p_j})\ dx$$
where for each $x \in H$, $g_j^x: H_j \to \R^+$ is the function $g_j^x(y) := f_j(B_j x - y) f'_j(y)$.  Applying \eqref{BL} to the inner integral we conclude that
$$ \BL(\B,\p;\f) \BL(\B,\p;\f') \leq \BL(\B,\p) \int_H \prod_{j=1}^m (\int_{H_j} g_j^x)^{p_j}\ dx.$$
But clearly $\int_{H_j} g_j^x = f_j * f'_j(B_j x)$, hence
$$ \BL(\B,\p)^2 \leq \BL(\B,\p) \int_H \prod_{j=1}^m ((f_j*f'_j) \circ B_j)^{p_j}.$$
On the other hand, from \eqref{BL} we have
$$ \int_H \prod_{j=1}^m ((f_j*f'_j) \circ B_j)^{p_j} = \BL(\B,\p; \f*\f') \prod_{j=1}^m (\int_H f_j * f'_j)^{p_j} = \BL(\B,\p; \f*\f')$$
since by the normalisation we have $\int_H f_j * f'_j = (\int_H f_j) (\int_H f'_j) = 1$.  Combining these inequalities we obtain
the claim.
\end{proof}

\begin{remark} Lemma \ref{factor} can be viewed as a degenerate version of this inequality, in which $\f'$ is Lebesgue measure on
the subspace $V$.
\end{remark}

By applying this lemma we can conclude the following closure properties of extremisers.

\begin{lemma}[Closure properties of extremisers]\label{closures}  Let $(\B,\p)$ be a 
Brascamp--Lieb datum with $\BL(\B,\p)$ finite.
\begin{itemize}
\item (Scale invariance) If $\f = (f_j)_{1 \leq j \leq m}$ is an extremising input, 
then so is $(f_j(\lambda \cdot))_{1 \leq j \leq m}$ for any non-zero real number $\lambda$.
\item (Homogeneity) If  $\f = (f_j)_{1 \leq j \leq m}$ is an extremising input, 
then so is \\$(c_j f_j)_{1 \leq j \leq m}$ for any non-zero real numbers $c_1, \ldots, c_m$.
\item (Translation invariance) If $\f = (f_j)_{1 \leq j \leq m}$ is an extremising input, 
then so is $(f_j(\cdot-B_j x_0))_{1 \leq j \leq m}$
for any $x_0 \in H$.
\item (Closure in $L^1$)  If $\f^{(n)} = (f_j^{(n)})_{1 \leq j \leq m}$ is a sequence of extremising inputs, which converges in the product $L^1$ sense to another input $\f = (f_j)_{1 \leq j \leq m}$, so $\lim_{n \to \infty} \|f_j^{(n)} - f_j \|_{L^1(H_j)} = 0$ for all $1 \leq j \leq m$, then
$\f$ is also an extremising input.
\item (Closure under convolution) If $\f$ and $\f'$ are extremisers, then so is $\f * \f'$.
\item (Closure under multiplication)  If $\f$ and $\f'$ are extremisers, then the input
 $(f_j(\cdot-B_j x_0) f'_j(\cdot))_{1 \leq j \leq m}$ is an extremiser for almost 
every $x_0 \in H$ for which $\int_H f_j(x-B_j x_0) f'_j(x)\ dx > 0$ for all $j$.
If furthermore $\f$ and $\f'$ are bounded, then the ``almost'' in ``almost every'' can be removed.
\end{itemize}
\end{lemma}

\begin{proof} From Lemma \ref{necc} we see that \eqref{scaling} holds, which guarantees the scale invariance property.  The closure in $L^1$ follows
from standard arguments since one can use \eqref{BL} and the
hypothesis that $\BL(\B,\p)$ is finite to show that $\int_H
\prod_{j=1}^m (f_j^{(n)} \circ B_j)^{p_j}$ converges to $\int_H
\prod_{j=1}^m (f_j \circ B_j)^{p_j}$. The homogeneity and
translation invariance can be verified by direct computation.  The
closure under convolution follows from Lemma \ref{conv}.  Now
consider the closure under multiplication. Let us write $\tilde
f_j(x) := f_j(-x)$, and let us repeat the proof of Lemma
\ref{conv} with $\f = (f_j)_{1 \leq j \leq m}$ replaced by $\tilde
f = (\tilde f_j)_{1 \leq j \leq m}$. Note from scale invariance
that $\tilde f$ is still an extremiser.  We then argue as before
to obtain \begin{eqnarray*} \begin{aligned} \BL(\B,\p)^2 = \int_H
(\int_H \prod_{j=1}^m (g_j^x \circ B_j)^{p_j})\ dx &\leq
\BL(\B,\p)
\int_H \prod_{j=1}^m (\int_{H_j} g_j^x)^{p_j}\ dx\\
&= \BL(\B,\p) \BL(\B,\p; \tilde \f* \f'),
\end{aligned}
\end{eqnarray*}
where $g_j^x(y) := f_j(y-B_jx) f'_j(y)$. Since $\BL(\B,\p;\tilde
\f * \f') \leq \BL(\B,\p)$, the above inequality must in fact be
equality, and thus
$$ \int_H \prod_{j=1}^m (g_j^x \circ B_j)^{p_j}\ dx = \BL(\B,\p) \prod_{j=1}^m (\int_{H_j} g_j^x)^{p_j}$$
for almost every $x$.  Also from the Fubini-Tonelli theorem we know that $\int_{H_j} g_j^x$ is finite for almost every $x$.  This proves the first
part of the closure under multiplication.  Also, since $f_j$ and $f'_j$ are integrable and bounded, the convolution $g_j^x$ lies in $L^1(H_j)$ continuously in $x$.  Since we are assuming $\BL(\B,\p)$ to be finite, this means that the left-hand side of the above expression also depends continuously on $x$. Thus in fact we have equality for every $x$, not just almost every $x$, and the second part of the closure under
multiplication follows.
\end{proof}

If we specialise the above lemma to the case of centred gaussian extremisers and use Theorem \ref{lieb-thm}, observing that the convolution or
product of two gaussians is again a gaussian, we conclude

\begin{corollary}[Closure properties of gaussian extremisers]\label{closures-gauss}  Let 
$(\B,\p)$ be a Brascamp--Lieb datum with $\BL_\g(\B,\p)$ finite.
\begin{itemize}
\item (Scale invariance) If $\A$ is a gaussian extremiser, then so is $\lambda \A$ for any $\lambda > 0$.
\item (Topological closure)  If $\A^{(n)}$ is a sequence of extremising gaussian inputs which converge to a gaussian input $\A$, then $\A$ is also an extremiser.
\item (Closure under harmonic addition) If $\A = (A_j)_{1 \leq j \leq m}$ and $\A' = (A'_j)_{1 \leq j \leq m}$ are extremisers, then so is
$((A_j^{-1} + (A'_j)^{-1})^{-1})_{1 \leq j \leq m}$.
\item (Closure under addition)  If $\A$ and $\A'$ are gaussian extremisers, then so is $\A + \A'$.
\end{itemize}
\end{corollary}

Note that these properties can also be deduced from Proposition \ref{normal-form} (or Proposition \ref{extremiser-char} below), although the closure under harmonic addition is not
particularly obvious from these propositions.  By applying all of these closure properties we can now conclude

\begin{proposition}\label{extremisor}  Let $(\B,\p)$ be an extremisable Brascamp--Lieb datum $(\B,\p)$.  Then
$\BL(\B,\p)=\BL_\g(\B,\p)$ and $(\B,\p)$ is gaussian-extremisable.
\end{proposition}

Carlen, Lieb and Loss have proved this in the rank-one case; see 
Theorem 5.4 of \cite{CLL}.

\begin{proof}  Intuitively, the idea (which was first discussed in Barthe \cite{Barthe}) is to start with an extremising
input $\f$ and repeatedly convolve it using the closure under convolution property of extremisers, rescale these convolutions using the
scale invariance, and then take limits using the central limit theorem and the closure under $L^1$ to obtain a gaussian extremiser.  However
there is a technical difficulty because the central limit theorem requires some moment conditions on the input $\f$ which are not obviously available, and so we will instead proceed in stages, starting with an arbitrary extremiser $\f$ and successively replacing $\f$ with increasingly more regular and well-behaved extremisers, until we end up with a centred gaussian extremiser.

Let $\f =(f_j)_{1 \leq j \leq m}$ be an extremising input for a
Brascamp--Lieb datum $(\B,\p)$.  To begin with, all we know about
the $f_j$ are that they are non-negative, integrable, and that
$\int_{H_j} f_j > 0$.  However we can use the following trick to
create a more regular extremiser. Observe that the reflection
$\tilde \f(x) := (f_j(-x))_{1 \leq j \leq m}$ is also an
extremiser (by scale invariance).  Using closure under convolution
we conclude that $\f
* \tilde \f$ is also an extremiser, and is also symmetric around
the origin, and is strictly positive at zero.  In fact it must be strictly
positive on a small ball surrounding zero; this is because for each $j$ there must exist a set $E_j$ of positive
finite Lebesgue measure such that $f_j > c 1_{E_j}$ for some constant $c > 0$, and the convolution $c1_{E_j} * c1_{E_j}$ is continuous
and positive at zero. Thus,
replacing $\f$ with $\f * \tilde \f$ if necessary, we can (and
shall) assume that $\f$ is symmetric, and strictly
positive near the origin.

Now let $1 \leq \lambda \leq 2$ and $x_0 \in H$ be parameters with 
$\|x_0\|_H \leq \eps$ for some small $\eps$,
and consider the $m$-tuple $\f_{\lambda,x_0} = (f_j(x) 
f_j(\lambda x - B_j x_0))_{1 \leq j \leq m}$.
From the Fubini-Tonelli theorem one verifies that 
\begin{equation}\label{fub}
\int_{H_j} f_j(x) f_j(\lambda x - B_j x_0)\ dx < \infty
\end{equation}
for almost every $\lambda, x_0$, and by previous hypotheses we know that $f_j(x) f_j(\lambda x-B_j x_0)$ is strictly positive near the origin if $\eps$ is sufficiently small. In particular $\f_{\lambda,x_0}$ is an input for almost every $\lambda$, $x_0$ with $\|x_0\|_H$ small.  From closure under multiplication we thus
see that $\f_{\lambda,x_0}$ is an extremiser for almost every 
$\lambda$, $x_0$.  We now claim that $\f_{\lambda,x_0}$ obeys the moment condition
$$ \int_{H_j} f_j(x) f_j(\lambda x - B_j x_0) (1 + \|x\|_{H_j})\ dx < \infty$$
for each $j$.  To see this, we fix $j$.  By \eqref{fub} we only need to show that
$$ \int_{H_j: \|x\|_{H_j} > 1} f_j(x) f_j(\lambda x - B_j x_0) \|x\|_{H_j}\ dx < \infty.$$
By the Fubini-Tonelli theorem it then suffices
to show that
$$ \int_{x_0 \in H: \|x_0\|_H \leq \eps} \int_{1}^2 \int_{x \in H_j: \|x\|_{H_j} > 1} f_j(x) f_j(\lambda x - B_j x_0) \|x\|_{H_j}\ dx d\lambda dx_0
< \infty.$$
In fact we claim that
$$ \int_{x_0 \in H: \|x_0\|_H \leq \eps} \int_1^2 f_j(\lambda x - B_j x_0)\ d\lambda dx_0 \leq \frac{C}{\|x\|_{H_j}} \int_{H_j} f_j$$
for all $x \in H_j$ with $\|x\|_{H_j} > 1$, and some finite constant $C = C(\eps,B_j, H_j)$ depending only on $\eps$, $B_j$, and $H_j$.
To see this, we divide the interval $\{ 1 \leq \lambda \leq 2 \}$ into $O( 1/\|x\|_{H_j} )$ intervals $I$ of radius $O( \|x\|_{H_j} )$.
For each interval $I$, we observe that 
$$ \int_{x_0 \in H: \|x_0\|_H \leq \eps} f_j(\lambda x - B_j x_0)\ d\lambda dx_0 \leq C \int_{\|y-\lambda_I x\|_{H_j} \leq C(1+\eps)} f_j(y)\ dy$$
where $\lambda_I$ is the midpoint of $I$ and the constants $C$ can depend on $\eps,B_j,H_j$.  
The claim then follows by integrating in $\lambda \in I$, then summing in $I$, 
observing that the balls $\{ y \in H_j: \|y-\lambda_I x\|_{H_j} \leq C(1+\eps) \}$ have an overlap of at most $C$.

To conclude, we have located an extremiser $\f$ which is positive near the origin, and obeys the moment condition
$\int_{H_j} (1 + \|x\|_{H_j}) f_j(x)\ dx < \infty$.  We can then use closure under multiplication again, replacing the $f_j$ by $f_j(x) f_j(x - B_j x_0)$ for some small $x_0 \in H$ and arguing as before, to improve this moment condition to $\int_{H_j} (1 + \|x\|_{H_j}^2) f_j(x)\ dx < \infty$.
Indeed by iterating this we can ensure that $\int_{H_j} (1 + \|x\|_{H_j}^N) f_j(x)\ dx < \infty$ for any specified $N$ (e.g., $N = 100 \dim(H_j)$).
In particular, each $f_j$ is square integrable, which implies that the Fourier transforms $\hat f_j$ are also square-integrable.  It is also bounded,
since $f_j$ is integrable.  One can
then replace $|\hat f_j(\xi)|^2$ with $|\hat f_j(\xi)|^2 |\hat f_j(\lambda \xi)|^2$ for any $1 \leq \lambda \leq 2$, by convolving each $f_j$ with a rescaled version of itself; this may slightly reduce the amount of moment conditions available but this is of no concern since we can make $N$ arbitrary.  Note that the new input will still be an extremiser, thanks to Lemma \ref{closures}.  By arguing as before (but now in the Fourier domain) we see that we have the moment conditions $\int_{H_j} |\hat f_j(\xi)|^2 |\hat f_j(\lambda \xi)|^2 \|\xi\|_{H_j} \ d\xi$ for almost every $\lambda$.  Thus by replacing $\f$ if necessary  we can assume that the extremiser $\f$ obeys the regularity condition $\int_{H_j} |\hat f_j(\xi)|^2 (1 + \|\xi\|_{H_j})\ d\xi < \infty$ for each $j$.  Indeed one can iterate this argument and obtain an extremiser for which
$\int_{H_j} |\hat f_j(\xi)|^2 (1 + \|\xi\|_{H_j}^N)\ d\xi < \infty$.  To conclude, we have now obtained an extremiser which has any specified amount of Sobolev regularity and decay.

We can now convolve $\f$ with its reflection $\tilde \f$ as before to recover the symmetry of $\f$.  We may also normalise $\int_{H_j} f_j = 1$ for all $j$.

The input $\f$ now obeys enough regularity for the central limit
theorem.  If we set $\f^{(n)} = (f^{(n)}_j)_{1 \leq j \leq m}$ to
be the rescaled iterated convolution
$$ f^{(n)}_j(x) := n^{(\dim H_{j})/2}f_j * \ldots * f_j(\sqrt{n} x)$$
where $f_j * \ldots * f_j$ is the $n$-fold convolution of $f_j$,
then each of the $\f^{(n)}$ are extremisers thanks to Lemma
\ref{closures}. Also, the central limit theorem shows that
$f^{(n)}_j$ converges in the $L^1$ topology (for instance) to a
centred gaussian, normalised to have total mass one.  Applying
Lemma \ref{closures} again we conclude that we can find an
extremiser which consists entirely of centred gaussians, which
thus implies that every extremisable Brascamp--Lieb datum is gaussian
extremisable as desired.
\end{proof}
\begin{remark}
Proposition \ref{extremisor} implies Lieb's theorem (Theorem \ref{lieb-thm}) for extremisable data $(\B,\p)$, though this is not directly
interesting since it is difficult to verify that $(\B,\p)$ is extremisable without first using Lieb's theorem.
However, if one assumes Lieb's theorem as a ``black box'', then the above results do shed some light on why the monotonicity formula
approach worked in the gaussian-extremisable case, as follows.  Suppose that $(\B,\p)$ is gaussian-extremisable with some gaussian
extremiser $\A = (A_j)_{1 \leq j \leq t}$.   Using Lieb's theorem we see that $(\exp(-\pi \langle A_j x, x \rangle_{H_j}))_{1 \leq j \leq m}$
is then an extremiser. By Lemma \ref{closures}, for any $t > 0$ the heat kernel
$$ \K_\A(t)(x) := \left(\det( A_j / t )^{1/2} \exp( - \pi \langle A_j x, x \rangle_{H_j} / t )\right)_{1 \leq j \leq t}$$
is also an extremiser.  Define the heat operators $U_\A(t)$ on inputs $\f$ by $U_\A(t) \f := \f * \K_\A(t)$.
Applying Lemma \ref{conv} we see that
$$ \BL(\B,\p; U_\A(t) f) \geq \BL(\B,\p; f) \hbox{ for all } t > 0.$$
Indeed, thanks to the well-known semigroup law
$$U_\A(t+s) = U_\A(t) U_\A(s) \hbox{ for all } s,t > 0$$
for the heat equation, we conclude that $\BL(\B,\p; U_\A(t) f)$ is non-decreasing in time.  Also, one can compute (using \eqref{scaling})
that $\lim_{t \to +\infty} \BL(\B,\p; U_\A(t) \f) \leq \BL_\g(\B,\p; \A)$ for all sufficiently well-behaved $\f$ (e.g., rapidly decreasing $\f$ will suffice).  From this we can recover the type of heat flow monotonicity that was so crucial in the proof of Proposition \ref{gbl-prop}.
\end{remark}

\section{Structural theory II : semisimplicity and existence of extremisers}\label{direct-sec}

Let $(\B,\p)$ be a Brascamp--Lieb datum obeying the conditions
\eqref{scaling}, \eqref{dimension}, which we have shown to be
necessary and sufficient for the Brascamp--Lieb constant
$\BL(\B,\p) = \BL_\g(\B,\p)$ to be finite. We now explore further
the question of when $(\B,\p)$ is gaussian-extremisable.  We
already have Theorem \ref{sufbl}, which shows that when there are
no critical subspaces then $(\B,\p)$ is gaussian-extremisable.
However it is certainly possible to be gaussian-extremisable in
the presence of critical spaces; consider for instance H\"older's
inequality (Example \ref{holder-ex}), which has plenty of gaussian
extremisers but for which every proper subspace is critical.  To resolve
this question more satisfactorily we need the notion of an
\emph{indecomposable $m$-transformation}.

\begin{definition}[Direct sum]  
If $\B = (H, (H_j)_{1 \leq j \leq m}, (B_j)_{1 \leq j \leq m})$ and $\B'  = $ $ \hfill$
$(H', (H'_j)_{1 \leq j \leq m}, (B'_j)_{1 \leq j \leq m})$ are $m$-transformations, we define the \emph{direct sum} $\B \oplus \B'$ to be the $m$-transformation
$$ \B \oplus \B' := (H \oplus H', (H_j \oplus H'_j)_{1 \leq j \leq m}, (B_j \oplus B'_j)_{1 \leq j \leq m})$$
where $H \oplus H' := \{ (x,x'): x \in H, x' \in H' \}$ is the Hilbert space direct sum of $H$ and $H'$ with the usual direct sum inner product
$\langle (x,x'), (y,y') \rangle_{H \oplus H'} := \langle x, y \rangle_H + \langle x', y' \rangle_{H'}$, and $B_j \oplus B'_j: H \oplus H' \to H_j \oplus H'_j$ is the direct sum $(B_j \oplus B'_j)(x,x') := (B_j x, B'_j x')$.  We say that an $m$-transformation $\B$ is \emph{decomposable} if it
is equivalent to the direct sum $\B_1 \oplus \B_2$ of two $m$-transformations $\B_1,\B_2$ whose domains have strictly smaller dimension than that
of $\B$, in which case we refer to $\B_1$ and $\B_2$ as \emph{factors} of $\B$.  We say that $\B$ is \emph{indecomposable} if it is not decomposable.
\end{definition}

It is easy to verify that if $(\B_1,\p), (\B_2,\p)$ are Brascamp--Lieb data with domains $H_1, H_2$ which obey \eqref{scaling},
and $(\B,\p) = (\B_1 \oplus \B_2,\p)$ is the direct sum (with domain $H_1 \oplus H_2$), then the subspaces $H_1 \oplus \{0\}$ and $\{0\} \oplus H_2$
of $H_1 \oplus H_2$ are critical subspaces with respect to $(\B,\p)$.  Furthermore, the restriction of $\B$ to $H_1 \oplus \{0\}$ is equivalent
to $\B_1$, and the quotient of $\B$ by $H_1 \oplus \{0\}$ is equivalent to $\B_2$; similarly with the roles of 1 and 2 reversed.
In particular this implies (from Lemma \ref{factor} and Lemma \ref{scaling-lemma}) that
\begin{equation}\label{blp}
 \BL( \B_1 \oplus \B_2, \p ) = \BL( \B_1, \p ) \BL( \B_2, \p ).
\end{equation}
We of course have a similar assertion for the gaussian Brascamp--Lieb constants.

\begin{remark} The Krull--Schmidt theorem for quiver representations (see e.g., \cite{ARS}) implies (as a special case) that
every $m$-transformation has a factorisation (up to equivalence) as the direct sum
of indecomposable $m$-transformations, and furthermore that this factorisation is unique up to equivalence and permutations.  The question of classifying what the indecomposable factors are of a given $m$-transformation, however, is
quite difficult.
In the rank one case (Example \ref{bbl}) a satisfactory description was given by Barthe \cite[Proposition 1]{Barthe}.  In fact the indecomposable (and prime) components can be described explicitly in this case as follows.  Let us introduce the relation $\bowtie$ on $\{1,\ldots,m\}$ by requiring
$i \bowtie j$ whenever there exists a collection $I \subset \{1,\ldots,m\} \backslash \{i,j\}$ such that $(v_k)_{k \in I \cup \{i\}}$ and
$(v_k)_{k \in I \cup \{j\}}$ are bases for $H$.  Let $\sim$ be the transitive 
completion of $\bowtie$.  Then the indecomposable components of $H$ take
the form $\linearspan( v_j: j \in I )$, where $I$ is any equivalence class for $\sim$.  
In the higher rank case we cannot expect such a completely explicit factorisation; at a 
minimum, we must allow for factorisation to only be unique up to equivalence. For instance, 
for H\"older's inequality (Example \ref{holder-ex}), any decomposition of 
$H$ into $\dim(H)$ independent one-dimensional subspaces will induce a 
factorisation of the associated Brascamp--Lieb datum into indecomposables; 
these factorisations are only unique up to equivalence.
\end{remark}

We now give a geometric criterion for indecomposability.

\begin{definition}[Critical pair]  Let $\B$ be an $m$-transformation.
A pair $(V,W)$ of subspaces of $H$ is said to be a \emph{critical pair} for $\B$ if $V$ and $W$ are complementary in $H$
(thus $V \cap W = \{0\}$ and $V+W = H$), and for each $j$, $B_j V$ and $B_j W$ are complementary in $H_j$.  We say the critical pair $(V,W)$
is \emph{proper} if $\{0\} \subsetneq V,W \subsetneq H$.
\end{definition}

\begin{remarks}
Observe that the exponents $p_j$ play no role in this definition.  If the $B_j$ are surjective, then $(0,H)$ and $(H,0)$ are of course critical pairs, though they are not proper.  Furthermore, in this case one can use identities such as $\dim(V) = \dim(B_j V) + \dim(V \cap \ker(B_j))$ to
formulate an equivalent definition of critical pair by requiring that $V,W$ are complementary in $H$, and $V \cap \ker(B_j)$ and $W \cap \ker(B_j)$
are complementary in $\ker(B_j)$.  In the rank one case $\dim(H_j) \equiv 1$, this condition 
was identified by Barthe \cite{Barthe}, who used it in obtaining a factorisation of Brascamp--Lieb data into indecomposable components.  We shall extend that analysis here to the higher rank case.
\end{remarks}

\begin{example}\label{critpair-ex} For H\"older's inequality (Example \ref{holder-ex}), every 
complementary pair in $\R^n$ is a critical pair.  For the Loomis--Whitney
inequality (Example \ref{lw-ex}), any co-ordinate plane and its orthogonal complement will be a critical pair.
For Young's inequality (Example \ref{bl-ex}) with $d=1$, there are no proper critical pairs (compare with Example \ref{crit-ex}).
In the rank-one case (Example \ref{bbl}), $(V,W)$ is a critical pair if and only if $V$ and $W$ are complementary, and $\{v_1,\ldots,v_m\} \subseteq V \cup W$.  Also, if $(V, W)$ is a critical pair for $\B$, and $\B'$ is equivalent to $\B$ with intertwining transformations $C, C_j$, then $(C^{-1} V, C^{-1} W)$ is a critical pair for $\B'$.
\end{example}

\begin{example}  If $\B$ is an $m$-transformation with domain $H$, then the direct sum $\B \oplus \B$ on $H \oplus H$
has $(H \oplus \{0\}, \{0\} \oplus H)$ as an obvious critical pair.  But there is also a ``diagonal'' critical pair
$$ ( \{ (x,x): x \in H \}, \{ (y,-y): y \in H \} ).$$
By taking tensor products of inputs in the first critical pair and then projecting onto the second critical pair (in the spirit of
Lemma \ref{factor}) one can recover K. Ball's convolution inequality Lemma \ref{conv}, as well as several of the closure properties
in Lemma \ref{closures}.  We omit the details.
\end{example}

\begin{lemma}[Geometric criterion for indecomposability]\label{geom-crit}  An $m$-transformation $\B$ is indecomposable if and only if it has no proper critical pairs.
\end{lemma}

\begin{proof}  If $\B$ is decomposable, then it is equivalent to a direct sum $\B_1 \oplus \B_2$ on a domain $H_1 \oplus H_2$ where $H_1$ and $H_2$
are not $\{0\}$.  The pair
$(H_1 \oplus \{0\}, \{0\} \oplus H_2)$ can be easily seen to be a proper critical pair for $\B_1 \oplus \B_2$, and thus $\B$ also has a proper critial pair by the remarks in Example \ref{critpair-ex}.

Now suppose conversely that $\B$ has a proper critical pair $(V,W)$.  Since $V$ and $W$ are complementary, we can apply an invertible linear  transformation $C$ on $H$ to make $V$ and $W$ orthogonal complements, while replacing $\B$ with an equivalent $m$-transformation.  Thus we may assume without loss of generality that $V$ and $W$ are orthogonal complements.  Similarly we may assume that $B_j V$ and $B_j W$ are orthogonal complements.
But then $\B$ is canonically equivalent to the direct sum $\B_V \oplus \B_W$ of its restrictions to $V$ and $W$, and is thus decomposable.
\end{proof}

\begin{remark}\label{factor-remark} The above proof in fact shows that if $(V,W)$ is a critical pair for $\B$, then $\B$ is equivalent to $\B_V \oplus \B_W$.
\end{remark}

One feature of critical pairs is that they are ``universally critical'', in the sense 
that they are critical for all admissible exponents $\p$:

\begin{lemma}[Critical pairs are critical subspaces]\label{critpaircrit}  Let $(\B, \p)$ be a 
Brascamp--Lieb datum obeying \eqref{dimension}, \eqref{scaling}, and let $(V,W)$ be a proper 
critical pair for $\B$.  Then $V$ and $W$ are both critical subspaces for $(\B,\p)$.
\end{lemma}

\begin{remark} This lemma implies that all simple Brascamp--Lieb data are indecomposable.  The converse is not true however; consider for instance the indecomposable example in Example \ref{gauss-fail}.
\end{remark}

\begin{proof} From \eqref{dimension}, \eqref{scaling}, and the hypothesis that $(V,W)$ is a critical pair, we have
\begin{align*}
\dim(H) &= \dim(V) + \dim(W) \\
&\leq \sum_{j=1}^m p_j \dim(B_j V) + \sum_{j=1}^m p_j \dim(B_j W) \\
&= \sum_{j=1}^m p_j \dim(H_j) \\
&= \dim(H)
\end{align*}
and hence the inequality above must be equality.  This implies that $V$ and $W$ are critical, as claimed.
\end{proof}

Critical pairs are related to extremisability in two ways.  First of all, extremisability of a product is equivalent to extremisability
in the factors:

\begin{lemma}[Extremisability factors]\label{tensor}  Let $(\B, \p)$ be a Brascamp--Lieb datum
obeying the conditions \eqref{scaling}, \eqref{dimension}, and let $(V,W)$ be a critical pair for $\B$.  Then
$(\B,\p)$ is extremisable if and only if $(\B_V,\p)$ and $(\B_W,\p)$ are both extremisable.  Similarly,
$(\B,\p)$ is gaussian-extremisable if and only if $(\B_V,\p)$ and $(\B_W,\p)$ are both gaussian-extremisable.
\end{lemma}

\begin{proof} By Remark \ref{factor-remark} we know that $\B$ is equivalent to $\B_V \oplus \B_W$.  
From Lemma \ref{scaling-lemma} and Lemma \ref{factor} we thus see that if $(\B,\p)$ is extremisable then $(\B_V,\p)$ and $(\B_W,\p)$ are extremisable,
and if $(\B,\p)$ is gaussian-extremisable then $(\B_V,\p)$ and $(\B_W,\p)$ are gaussian-extremisable.

In the converse direction, if equality is attained in \eqref{BL}
for the input $(f_{j,V})_{1 \leq j \leq m}$ for $(\B_V,\p)$ and for the input $(f_{j,W})_{1 \leq j \leq m})$ for $(\B_W,\p)$,
then the direct sums $(f_{j,V} \oplus f_{j,W})_{1 \leq j \leq m}$ will attain equality in $(\B_V \oplus \B_W, \p)$, thanks to
\eqref{blp} (or Lemma \ref{factor}).  The claim for the extremisation problem \eqref{BLfunctional} is 
similar.
\end{proof}

Furthermore, in the presence of an extremiser one can obtain a converse to Lemma \ref{critpaircrit}.

\begin{lemma}[Critical complements for gaussian-extremisable data]\label{critpaircrit-converse}  
Let $(\B,\p)$ be a gaussian-extremisable Brascamp--Lieb datum.  Then for every critical 
subspace $V$ of $(\B,\p)$ there exists a complementary space $W$ to $V$ such that $(V,W)$ is 
a critical pair for $\B$. If furthermore $(\B,\p)$ is geometric, then $W$ is the orthogonal complement of $V$.
\end{lemma}

\begin{proof}  We may drop exponents for which $p_j=0$, and thus assume that $p_j > 0$ for all $j$.
By Proposition \ref{normal-form} we may assume without loss of generality that $(\B,\p)$ is a 
geometric Brascamp--Lieb datum, thus we may assume $H_j$ is a subspace of $H$, $B_j: H \to H_j$ is the orthogonal projection,
and \eqref{pbb} holds.

Now let $\Pi: H \to V$ be the orthogonal projection onto $V$.  Observe that $B_j \Pi$ is a contraction from $H$ to $B_j V$ and hence has
trace at most $\dim(B_j V)$.  Thus, since $V$ is a critical subspace,
\begin{align*}
\dim(V) &= \tr_H(V) \\
&= \tr_H(\sum_{j=1}^m p_j B_j V) \\
&\leq \sum_{j=1}^m p_j \dim(B_j V) \\
&= \dim(V)
\end{align*}
and hence we must in fact have $\tr_H(B_j \Pi) = \dim(B_j V)$ for all $j$.  Thus $B_j \Pi$ is 
in fact a co-isometry (the adjoint is an isometry), which means that $B_j V$
and $B_j(V^\perp)$ are orthogonal complements in $H_j$.  This implies that $(V, V^\perp)$ form 
a critical pair, and the claim follows.
\end{proof}

We can now give a characterisation of when extremisers exist.

\begin{theorem}[Necessary and sufficient conditions for extremisers]\label{extremiser}  Let 
$(\B,\p)$ be a Brascamp--Lieb datum
with $p_j > 0$ for all $j$.  Then the following nine statements are equivalent.
\begin{enumerate}
\item  $(\B,\p)$ is extremisable.
\item  $(\B,\p)$ is gaussian-extremisable.
\item  A local extremiser exists to \eqref{BLfunctional(not)}.
\item  There exists a gaussian input $\A = (A_j)_{1 \leq j \leq m}$ such that the matrix $M := \sum_{j=1}^m p_j B_j^* A_j B_j$ obeys $A_j^{-1} = B_j M^{-1} B_j^*$ for all $j$.
\item  The scaling condition \eqref{scaling} holds, $\B$ is non-degenerate, and
there exists a gaussian input $\A = (A_j)_{1 \leq j \leq m}$ such that the matrix $M := \sum_{j=1}^m p_j B_j^* A_j B_j$ obeys $A_j^{-1} \geq B_j M^{-1} B_j^*$ for all $j$.
\item  The scaling condition \eqref{scaling} holds, $\B$ is non-degenerate, and
there exists a gaussian input $\A = (A_j)_{1 \leq j \leq m}$ such that $B_j^* A_j B_j \leq \sum_{i=1}^m p_i B_i^* A_i B_i$.
\item  $(\B,\p)$ is equivalent to a geometric Brascamp--Lieb datum.
\item  The bounds \eqref{dimension}, \eqref{scaling} hold, and every critical subspace $V$ is part of a critical pair $(V,W)$.
\item  The bounds \eqref{dimension}, \eqref{scaling} hold, and every indecomposable factor of $(\B,\p)$ is simple.
\end{enumerate}
\end{theorem}

One can view the equivalence (a) $\iff$ (i) as a statement that extremisability is equivalent to being
``semisimple'' (the direct sum of simple factors).

\begin{proof} The implication (a) $\implies$ (b) is given by Proposition \ref{extremisor}, 
while the converse direction (b) $\implies$ (a)
is given by Theorem \ref{lieb-thm}.  The equivalence of (b)-(g) 
follows from Proposition \ref{normal-form}.
The implication (b) $\implies$ (h) is given by Theorem \ref{sufbl} 
and Lemma \ref{critpaircrit-converse}.  Now let us
prove the converse implication (h) $\implies$ (b).  
We induct on the dimension of $H$.  If there are no proper critical subspaces 
then the claim follows from
Theorem \ref{sufbl}.  Now suppose that there is a proper critical space 
$V$, which is then part of a critical pair $(V,W)$ of proper subspaces.
We then pass to the factors $(\B_V, \p)$ and $(\B_W, \p)$.  By Lemma \ref{critsum} these factors
also obey \eqref{dimension}, \eqref{scaling}.  Also, we observe that every critical subspace 
$V'$ of $V$ is part of a critical pair $(V', W')$ in
$V$.  To see this, observe from hypothesis that $V'$ is part of a critical pair $(V', W'')$ in $H$.  
Now set $W' := V \cap W''$; since $V'$ and $W''$ are complementary in $H$, and $V' \subset V$, we see 
that $V'$ and $W'$ are complementary in $V$.  In particular, $B_j V' + B_j W' = B_j V$.
On the other hand, since $B_j V'$ and $B_j W''$ are complementary in $H_j$, and $B_j W'$ is 
contained in $B_j W''$, we have $B_j V' \cap B_j W' = \{0\}$.  Thus $B_j V'$ and $B_j W'$ are 
complementary in $B_j V$, and hence $(V', W')$ is a critical pair in $V$ as claimed.

Applying the induction hypothesis, we conclude that the supremum in \eqref{BLfunctional} for $(\B_V,\p)$ is attained.
Similarly for $V$ replaced by $W$.  The claim now follows from Lemma \ref{critpaircrit-converse}.

Now we show the implication (b) $\implies$ (i).  Applying Lemma \ref{factor} we see that the supremum in \eqref{BLfunctional} is attained
for every indecomposable factor of $(\B,\p)$.  Applying the equivalence between (b) and (h) already proven, and noting that indecomposable factors
cannot contain proper critical pairs by definition, we thus see that every indecomposable factor contains no proper critical subspaces as claimed.

Finally, to show that (h) $\implies$ (b), we see by the equivalence of (b) and (h) that the supremum in \eqref{BLfunctional} is attained for every indecomposable factor of $(\B,\p)$.  Factoring $(\B,\p)$ in some arbitrary manner as a direct sum of indecomposable components and using Lemma \ref{tensor} we obtain the claim.
\end{proof}

It should be mentioned that the above equivalences only hold with the assumption $p_j > 0$.  Of course any exponent with $p_j = 0$
can be omitted without affecting most of the properties listed above (specifically, this does not affect (a)-(g)), but it does affect the
factorisation of data into indecomposable components.

\begin{example} Let us return to Example \ref{gauss-fail}.  This $3$-transformation $\B$ is indecomposable, and the datum $(\B,\p)$ is simple if $(p_1,p_2,p_3)$ lies in the interior of the triangle.  However if $(p_1,p_2,p_3)$ is on a vertex of this triangle, then one of the $p_j$ vanishes,
and on removing this exponent we are left with a $2$-transformation which is now decomposable (with the two components clearly being simple).  This is why
extremisers exist on the interior and vertices of the triangles but not on the open edges.
\end{example}

More generally, the above proposition yields an explicit test as to whether a given datum $(\B,\p)$ is extremisable.  First, one removes any
exponents for which $p_j = 0$.  Then one splits $\B$ into indecomposable components $(\B_i, \p)$.  If $\p$ lies in $S(\B_i) \cap \Pi(\B_i)^\circ$
for each component $\B_i$, then $(\B,\p)$ is extremisable, but if $\p$ lies on the boundary of $\Pi(\B_i)$ for any $i$ then $(\B,\p)$ is not
extremisable.

\section{Regularised and localised Brascamp--Lieb inequalities}\label{regular-sec}

We have seen a number of connections between Brascamp--Lieb constants, gaussians, and the heat equation.  We now pursue these connections
further by introducing a generalisation of the Brascamp--Lieb constants.  As a byproduct of this analysis we give an alternate proof of
Lieb's theorem (Theorem \ref{lieb-thm}) that does not rely on factorisation through critical subspaces.

We first examine the situation when the data are regularised at a certain scale; morally this 
corresponds to a discrete setting. This allows us to then examine the situation when the data 
are gaussian-localised, effected by appending a fixed gaussian multiplier to the input and then 
scaling. The gaussian-localised situation falls in principle under the framework of
\cite{L}. 

Our first definition formalises our notion of regularity.

\begin{definition}[Type $G$ functions]  Let $H$ be a Euclidean space, and let $G: H \to H$ be 
positive definite.  We define a \emph{type $G$ function} to be any function $u: H \to \R^+$ of the form
$$ u(x) = \mbox{$\det_H$}(G)^{1/2} \int_H \exp(-\pi \langle G (x-y), (x-y) \rangle_H)\ d\mu(y)$$
where $\mu$ is a positive finite Borel measure on $H$ with non-zero total mass.
If $\mu$ is a point mass, we say that $u$ is of \emph{extreme type $G$}; thus the functions
of extreme type $G$ are simply the translates and positive scalar multiples of the function $\exp(-\pi \langle G x, x \rangle_H)$.
\end{definition}

\begin{remark}
Observe that type $G$ functions are smooth and strictly positive.  Also, if $G_1 \geq G_2 > 0$ then every function of type $G_2$ is also of type $G_1$, since a gaussian such as $\exp(-\pi \langle G_2 x, x\rangle_H)$ can be expressed as a convolution of $\exp(-\pi \langle G_1 x, x \rangle_H)$ with
a positive finite measure (this can be seen for instance using the Fourier transform).  Positive finite measures themselves can be informally viewed
as ``functions of type $+\infty$''.
\end{remark}

\begin{remark}
Type $G$ functions arise naturally in the study of heat equations.  More precisely, if $u: \R^+ \times H \to \R^+$ is a solution to the heat
equation
$$ u_t = \frac{1}{4\pi} \div( G^{-1} \nabla u )$$
and $u(0) = \mu$ is a positive finite measure, then from the explicit solution
$$  u(t,x) = \mbox{$\det_H$}(G/t)^{1/2} \int_H \exp(-\pi \langle G (x-y), (x-y) \rangle_H / t)\ d\mu(y)$$
we see that $u(t)$ is of type $G/t$ for all $t > 0$.  More
generally, if $u(s)$ is of type $A$ for some $s \geq 0$ and $A >
0$, then $u(t)$ is of type $(A^{-1} + (t-s)G^{-1})^{-1}$ for all
$t > s$, as can be seen for instance by using the Fourier
transform, or alternatively by convolving one fundamental solution
to a heat equation with another using \eqref{hilbert}.  Note that
this operation of harmonic addition has already appeared in
Corollary \ref{closures-gauss}.
\end{remark}

\begin{definition}[generalised Brascamp--Lieb constant]  If $(\B,\p)$ is a Brascamp--Lieb datum,
we define a \emph{generalised Brascamp--Lieb datum} to be a triple
$(\B,\p,\G)$, where $(\B,\p)$ is a Brascamp--Lieb datum and $\G$ is
a gaussian input for $(\B,\p)$.  We shall refer to $\G$ as an \emph{$m$-type}. We say that an input $\f =
(f_j)_{1 \leq j \leq m}$ is of \emph{type $\G = (G_j)_{1 \leq j
\leq m}$} if each $f_j$ is of type $G_j$. We define the
\emph{generalised Brascamp--Lieb constant} $\BL(\B,\p,\G)$ to be
\begin{eqnarray}\label{BL-gen}
\BL( \B, \p, \G )
=\sup_{\f \text{ of type } \G }
\BL(\B,\p; \f).  
\end{eqnarray}
We also define the generalised gaussian Brascamp--Lieb constant
$\BL_\g(\B,\p,\G)$ by the formula
\begin{eqnarray}\label{BLfunctional-gen}
\BL_\g( \B,\p,\G ) := \sup_{\A \leq \G} \BL_\g(\B,\p; \A),
\end{eqnarray}
where the supremum extends over all gaussian inputs $\A = (A_j)_{1 \leq j\leq m}$ such that $A_j \leq G_j$ for all $j$.
We say that $(\B,\p,\G)$ is \emph{extremisable} if $\BL(\B,\p,\G)$ is finite and equality can be attained in \eqref{BL-gen} for some input $\f$ of
type $\G$, and is \emph{gaussian-extremisable} if $\BL_\g(\B,\p,\G)$ is finite and the supremum in \eqref{BLfunctional-gen} can be attained for
some gaussian input $\A \leq \G$.
\end{definition}

\begin{remark} One may restrict attention to non-degenerate $\B$ for the following reasons.  If $\bigcap_{j=1}^m \ker(B_j) \neq \{0\}$ then it is easy to see that both constants will be infinite.  If one or more of the $B_j$ is not surjective, one can simply restrict $H_j$ to $B_j V$ (and restrict
the quadratic form associated to $G_j$ to $B_j V$ also) to obtain an equivalent problem.
\end{remark}

Clearly we have
\begin{equation}\label{blg-trivial-gen}
\BL(\B,\p,\G) \geq \BL_\g(\B,\p,\G);
\end{equation}
we shall later show in Corollary \ref{qw} that this is in fact an equality, in analogy with Theorem \ref{lieb-thm}.
Clearly $\BL(\B,\p,\G) \leq \BL(\B,\p)$ and $\BL_\g(\B,\p,\G) \leq \BL_\g(\B,\p)$, but it is certainly possible for strict inequality to hold
for either inequality.  On the other hand, a simple regularisation argument shows that
\begin{equation}\label{tlim}
 \BL(\B,\p) = \lim_{\lambda \to \infty} \BL(\B,\p,\lambda\G); \quad \BL_\g(\B,\p) = \lim_{\lambda \to \infty} \BL_\g(\B,\p,\lambda\G)
 \end{equation}
for any $m$-type $\G$.  Furthermore, it turns out that there are analogues not only of Lieb's theorem 
but also of Theorem \ref{sufbl} in this generalised setting. We shall state these results presently, 
but let us first develop some preliminary lemmas.  We begin with a basic log-convexity estimate.

\begin{lemma}[Log-convexity lemma]\label{log-convex}  Let $G: H \to H$ be 
positive definite, let $u$ be of type $G$, and let $v$ of be of
extreme type $G$.  Then $u/v$ is log-convex (i.e., $\nabla^2 \log (u/v)$ is positive semi-definite).  If $u$ is not of extreme type $G$, then $u/v$ is
strictly log-convex (i.e., $\nabla^2\log (u/v)$ is 
positive definite).
\end{lemma}

\begin{proof} By applying an invertible linear transformation if necessary one can reduce to the case $H = \R^n$, $G = I_n$.  By translation we may then
assume that $v$ is the extreme type $I_n$ function
$$ v(x) = \exp(-\pi |x|^2).$$
We can write the type $I_n$ function $u(x)$ as
$$ u(x) = \int_{\R^n} \exp(-\pi |x-y|^2)\ d\mu(y)$$
and hence
\begin{equation}\label{uv-def}
\frac{u(x)}{v(x)} = \int_{\R^n} e^{2\pi \langle x, y \rangle / 2}\ e^{-\pi |y|^2} d\mu(y).
\end{equation}
At this point we could argue that $u/v$ is log-convex by observing that the exponential functions $e^{2\pi \langle x,y \rangle}$ are log-convex,
and that the superposition of log-convex functions is again log-convex.  However we shall give a more explicit argument which also yields the second claim.  Taking gradients of the above equation we obtain
$$ \frac{u(x)}{v(x)} \nabla \log \frac{u(x)}{v(x)} = \int_{\R^n}  2\pi y e^{2\pi \langle x, y \rangle}\ e^{-\pi |y|^2} d\mu(y)$$
and hence by \eqref{uv-def} again
$$ \int_{\R^n} (\nabla \log \frac{u(x)}{v(x)} -  2\pi y)  e^{2\pi \langle x, y \rangle}\ e^{-\pi |y|^2} d\mu(y) = 0.$$
We multiply this by an arbitrary transformation $A: \R^n \to \R^n$ and take divergences, to conclude
$$ \int_{\R^n} (\div A \nabla \log \frac{u(x)}{v(x)} + \langle 2\pi A y, 
\nabla \log \frac{u(x)}{v(x)} -  2\pi y \rangle )  e^{2\pi \langle x, y \rangle}\ e^{-\pi |y|^2} d\mu(y) = 0.$$
Rearranging this using \eqref{uv-def} again, we conclude
$$ \div A \nabla \log \frac{u(x)}{v(x)} = \frac{ \int_{\R^n} \langle A (\nabla \log \frac{u(x)}{v(x)} - 2\pi y), \nabla \log \frac{u(x)}{v(x)} -  2\pi y \rangle e^{2\pi\langle x,y \rangle}\ e^{-\pi |y|^2} d\mu(y) }
{ \int_{\R^n} e^{2\pi \langle x,y \rangle}\ e^{-\pi |y|^2} d\mu(y) }.$$
Observe that if $A$ is positive semi-definite, then the right-hand side is non-negative.  This shows that $\nabla^2 \log \frac{u}{v}$ is
log-convex as claimed.  Also, if $A$ is 
positive definite, then we see that for any fixed $x$ there is at most one $y$ for which
the integrand in the numerator is zero.  Thus if $\mu$ is not a point mass then the numerator is always strictly positive, and hence $\frac{u}{v}$ is strictly convex in this case.  This concludes the proof.
\end{proof}

\begin{corollary}\label{divag}  Let $G: H \to H$ be 
positive definite, let $u: H \to \R^+$ be of type $G$, and let $A: H \to H$ be positive semi-definite.
Then
$$ \div( A \nabla \log u ) \geq -2\pi \tr_H(AG).$$
If $A$ is not identically zero, then equality occurs if and only if $u$ is extreme type $G$.
\end{corollary}

\begin{remark} Intuitively speaking, this corollary asserts that positive solutions to the heat equation cannot diverge too rapidly, and can be viewed as an explanation for the monotonicity for heat flows.
\end{remark}

\begin{proof}
Let $v: H \to \R^+$ be the extreme type $G$ function $v(x) := \exp( - \pi \langle Gx, x \rangle_H \rangle)$.  Observe
that $\div(A \nabla \log v ) = -2\pi \tr_H(AG)$, and from Lemma \ref{log-convex} we have $\hfill$
$\div(A \nabla \log \frac{u}{v} ) \geq 0$,
with strict inequality if $A$ is not identically zero and $u$ is not extreme type $G$.  The claim follows.
\end{proof}

We can now give the monotonicity formula for multiple heat flows.  We begin with a rather general statement.

\begin{proposition}[Multilinear heat flow monotonicity formula]\label{flowform} Let $(\B,\p,\G)$ be a
generalised Brascamp--Lieb datum with all the $B_j$ surjective, and suppose that there is a gaussian 
input $\A = (A_j)_{1 \leq j \leq m}$ with $\A \leq \G$, such that the transformation
$M: H \to H$ defined by $M := \sum_{j=1}^m p_j B_j^* A_j B_j$ is invertible and obeys the inequality
\begin{equation}\label{bam}
 A_j^{-1} - B_j M^{-1} B_j^* \geq  0
 \end{equation}
and the identity
\begin{equation}\label{caj}
(A_j^{-1} - B_j M^{-1} B_j^*) (G_j - A_j) = 0
\end{equation}
for each $1 \leq j \leq m$.  Also, for each $1 \leq j \leq m$, let $\tilde u_j(1): H_j \to \R^+$ be of type $G_j$,
and let $\tilde u_j: [1,\infty) \times H_j \to \R^+$ be the solution to the heat equation
$$ \partial_t \tilde u_j = \frac{1}{4\pi} \div( A_j^{-1} \nabla \tilde u_j )$$
with initial data $\tilde u_j(1)$ at $t=1$.  Then the quantity
$$ Q(t) := t^{(\sum_{j=1}^m p_j \dim(H_j) - \dim(H))/2} \int_H \prod_{j=1}^m (\tilde u_j \circ B_j)^{p_j}$$
is monotone non-decreasing in time for $t > 1$.

Furthermore, if there is a $j$ for which equality in \eqref{bam} does not hold, and $\tilde u_j(1)$ is not of extreme type $G_j$, then
$Q(t)$ is strictly increasing for $t > 1$.
\end{proposition}

\begin{proof}  We have
\begin{equation}\label{tuj1}
\tilde u_j(1,x) = \mbox{$\det_{H_j}$}(G_j)^{1/2} \int_{H_j}
\exp(-\pi \langle G_j (x-y), (x-y) \rangle_{H_j} )\ d\mu_j(y)
\end{equation}
for some positive finite Borel measure $\mu_j$. A simple computation using the Fourier 
transform then shows that
\begin{eqnarray}\label{tuj}
\begin{aligned}
\tilde u_j(t,x) &= \mbox{$\det_{H_j}$}(G^{-1}_j + (t-1)
A_j^{-1})^{-1/2} \\
&\times \int_{H_j} \exp(-\pi \langle (G^{-1}_j + (t-1)
A_j^{-1})^{-1} (x-y), (x-y) \rangle_{H_j} )\ d\mu_j(y).
\end{aligned}\end{eqnarray}
By a Fatou lemma argument we may assume without loss of generality that the measures $\mu_j$ are compactly supported.  We may also assume $p_j > 0$ for all $j$ as the degenerate exponents $p_j = 0$ can be easily omitted.

Let $D_j: H_j \to H$ be any linear right-inverse to $B_j$, thus $B_j D_j = \id_{H_j}$.  Such an inverse exists since $B_j$ is assumed
to be surjective.  The exact choice of $D_j$ will not be relevant; for instance
one could take $D_j = B_j^* (B_j B_j^*)^{-1}$.  We will apply Lemma \ref{multitransport} with
\begin{align*}
I &:= (0,+\infty) \\
u_j &:= \tilde u_j \circ B_j \\
\vec v_j &:= -\frac{1}{4\pi} D_j A_j^{-1} (\nabla \log \tilde u_j) \circ B_j \\
\vec v &:= M^{-1} \sum_{j=1}^m p_j B_j^* A_j B_j \vec v_j \\
\alpha &:= \frac{1}{2}( \dim(H) - \sum_{j=1}^m p_j \dim(H_j)).
\end{align*}
Let us first verify the technical condition that $\vec v
\prod_{j=1}^m u_j^{p_j}$ is rapidly decreasing in space. Because
the $\mu_j$ are compactly supported, one can verify that $\nabla
\log \tilde u_j$ grows at most polynomially in space, locally
uniformly on $I$. Hence $\vec v$ also grows at most polynomially
in space.  Since $M$ is invertible, we have $\bigcap_{j=1}^m
\ker(B_j) = 0$.  Since each of the $\tilde u_j$ are rapidly
decreasing in space, locally uniformly on $I$, we see that
$\prod_{j=1}^m u_j^{p_j}$ is also, and the claim follows.

By Lemma \ref{multitransport}, we will now be done as soon as we verify the inequalities \eqref{transport-j}, \eqref{div-pos}, \eqref{log-pos}.
We begin with \eqref{transport-j}.  On the one hand, we have
$$ \partial_t u_j = (\partial_t \tilde u_j) \circ B_j = \frac{1}{4\pi} \div( A_j^{-1} \nabla \tilde u_j ) \circ B_j.$$
On the other hand, from the chain rule and the choice of $D_j$ we have
\begin{align*}
\div(\vec v_j u_j) &= -\div([ \frac{1}{4\pi} D_j A_j^{-1} (\nabla \log \tilde u_j) \tilde u_j ] \circ B_j) \\
&= -\frac{1}{4\pi} \div( D_j A_j^{-1} (\nabla \tilde u_j) \circ B_j ) \\
&= -\frac{1}{4\pi} \div( B_j D_j A_j^{-1} \nabla \tilde u_j) \circ B_j \\
&= -\frac{1}{4\pi} \div( A_j^{-1} \nabla \tilde u_j) \circ B_j
\end{align*}
and so the left-hand side of \eqref{transport-j} is zero.

Next we verify \eqref{div-pos}.  We expand the left-hand side as
\begin{align*}
\div( \vec v - \sum_{j=1}^m p_j \vec v_j ) &= \sum_{j=1}^m p_j \div( M^{-1} B_j^* A_j B_j \vec v_j - \vec v_j ) \\
&= -\frac{1}{4\pi} \sum_{j=1}^m p_j \div( (M^{-1} B_j^* A_j B_j - \id_{H_j}) D_j A_j^{-1} (\nabla \log \tilde u_j) \circ B_j ) \\
&= -\frac{1}{4\pi} \sum_{j=1}^m p_j \div( B_j (M^{-1} B_j^* A_j B_j - \id_{H_j}) D_j A_j^{-1} (\nabla \log \tilde u_j) ) \circ B_j \\
&= \frac{1}{4\pi} \sum_{j=1}^m p_j \div( (A_j^{-1} - B_j M^{-1} B_j^*) \nabla \log \tilde u_j ) \circ B_j.
\end{align*}
From \eqref{tuj} we see that $\tilde u_j(t)$ is of type $(G^{-1}_j + (t-1) A_j^{-1})^{-1}$.
Since $G_j \geq A_j > 0$ and $t > 1$,
we have $(G^{-1}_j + (t-1) A_j^{-1})^{-1} \leq G_j / t$.  Thus $\tilde u_j(t)$ is also of type
$G_j/t$.  Applying \eqref{bam} and
Corollary \ref{divag} we have
$$ \div( (A_j^{-1} - B_j M^{-1} B_j^*) \nabla \log \tilde u_j ) \geq
-2\pi \tr_{H_j}( (A_j^{-1} - B_j M^{-1} B_j^*) G_j / t ).$$
Using \eqref{caj} we conclude that
$$ \div( (A_j^{-1} - B_j M^{-1} B_j^*) \nabla \log \tilde u_j ) \geq
-\frac{2\pi}{t} \tr_{H_j}( (A_j^{-1} - B_j M^{-1} B_j^*) A_j ).$$
Inserting this into the preceding computation, we conclude
\begin{align*}
\div( \vec v - \sum_{j=1}^m p_j \vec v_j ) &\geq
-\frac{1}{2t} \sum_{j=1}^m p_j \tr_{H_j}( (A_j^{-1} - B_j M^{-1} B_j^*) A_j ) \\
&= \frac{1}{2t} \sum_{j=1}^m p_j (\tr_H(M^{-1} B_j^* A_j B_j) - \dim(H_j)) \\
&= \frac{1}{2t} (\tr_H( M^{-1} M ) - \sum_{j=1}^m p_j \dim(H_j)) \\
&= \frac{1}{2t} ( \dim(H) - \sum_{j=1}^m p_j \dim(H_j))\\
&= \frac{\alpha}{t}
\end{align*}
as desired.

Finally, we verify \eqref{log-pos}.  From the chain rule we have
$$ \nabla \log u_j = B_j^* (\nabla \log \tilde u_j) \circ B_j$$
and hence by definition of $\vec v_j$ and $D_j$
$$ \nabla \log u_j = -B_j^* A_j B_j \vec v_j.$$
Thus we can write the left-hand side of \eqref{log-pos} as
$$
\sum_{j=1}^m p_j \langle B_j^* A_j B_j (\vec v - \vec v_j), -\vec v_j \rangle_H.$$
On the other hand, by definition of $\vec v$ and $M$ we have
$$ \sum_{j=1}^m p_j B_j^* A_j B_j (\vec v - \vec v_j)  = M \vec v - \sum_{j=1}^m p_j B_j^* A_j B_j\vec v_j = 0.$$
Thus we can write the left-hand side of \eqref{log-pos} as
$$
\sum_{j=1}^m p_j \langle B_j^* A_j B_j (\vec v - \vec v_j), (\vec v - \vec v_j) \rangle_H.$$
Since $A_j$ is positive definite, we see that this expression is non-negative as desired.
This completes the proof of monotonicity.  The proof of strict monotonicity when equality does not hold in \eqref{bam} for some $j$, and
$\tilde u_j(1)$ is not of extreme type $G_j$, follows by a minor variation of the above argument which we omit.
\end{proof}

Let us compute the limiting behaviour of $Q(t)$ as $t \to \infty$, under the assumption that the $\mu_j$ are compactly supported.
After a change of variables we have
$$ Q(t) = t^{\sum_{j=1}^m p_j \dim(H_j)/2} \int_H \prod_{j=1}^m \tilde u_j(t,t^{1/2} B_j x)^{p_j}\ dx$$
and then after applying \eqref{tuj} we can write this as
\begin{align*}
\int_H &\prod_{j=1}^m \biggl[t^{\dim(H_j)/2} \mbox{$\det_{H_j}$}(G^{-1}_j + (t-1) A_j^{-1})^{-1/2} \\
&\int_{H_j} \exp(-\pi \langle (G^{-1}_j + (t-1) A_j^{-1})^{-1} (t^{1/2} B_j x-y), (t^{1/2} B_j x-y) \rangle_{H_j} )\ d\mu_j(y)\biggr]^{p_j}\ dx.
\end{align*}
Via another change of variables we may rewrite this as
\begin{align*} \int_H &\prod_{j=1}^m \biggl[\mbox{$\det_{H_j}$}(t^{-1} G^{-1}_j + (1-t^{-1}) A_j^{-1})^{-1/2} \\
&\int_{H_j} \exp(-\pi \langle (t^{-1} G^{-1}_j + (1-t^{-1}) A_j^{-1})^{-1} (B_j x-t^{-1/2} y), (B_j x-t^{-1/2} y) \rangle_{H_j} )\ d\mu_j(y)\biggr]^{p_j}\ dx.
\end{align*}

Taking limits as $t \to \infty$ using dominated convergence (which can be justified since $\mu_j$ is compactly supported) we conclude
$$ \lim_{t \to \infty} Q(t)
= \int_H \prod_{j=1}^m \left[\mbox{$\det_{H_j}$}(A_j^{-1})^{-1/2}
\int_{H_j} \exp(-\pi \langle A_j B_j x, B_j x \rangle_{H_j} )\
d\mu_j(y)\right]^{p_j}\ dx,$$ 
which simplifies using \eqref{hilbert} to
\begin{eqnarray*}  
\begin{aligned}
\lim_{t \to \infty} Q(t) 
&= \frac{\prod_{j=1}^m \det_{H_j}(A_j)^{p_j/2}}{\det_H(M)^{1/2}} 
\prod_{j=1}^m \mu_j(H_j)^{p_j} \\ 
&= \frac{\prod_{j=1}^m \det_{H_j}(A_j)^{p_j/2}}{\det_H(M)^{1/2}} 
\prod_{j=1}^m \left(\int_{H_j} \tilde u_j(1)\right)^{p_j}.
\end{aligned}
\end{eqnarray*}

In particular we have
$$ Q(1) \leq \frac{\prod_{j=1}^m \det_{H_j}(A_j)^{p_j/2}}{\det_H(M)^{1/2}} \prod_{j=1}^m \left(\int_{H_j} \tilde u_j(1)\right)^{p_j}.$$
We thus conclude (using a Fatou lemma argument to eliminate the hypothesis that $\mu_j$ is compactly supported)

\begin{corollary}[Towards a generalised Brascamp--Lieb inequality]\label{flowform-bl} Let $(\B,\p,\G)$ be 
a generalised Brascamp--Lieb datum with all the $B_j$ surjective,
and let $\A = (A_j)_{1 \leq j \leq m}$ be a gaussian input with $\A \leq \G$, such that the transformation
$M: H \to H$ defined by $M := \sum_{j=1}^m p_j B_j^* A_j B_j$ is invertible and obeys the inequality \eqref{bam} and the
identity \eqref{caj} for all $1 \leq j \leq m$.  Then
$$ \BL(\B,\p,\G) = \BL_\g(\B,\p;\A).$$
\end{corollary}

One special case of this Corollary is when $\G=\A$, as in this case the condition \eqref{caj} is automatic.  Another
special case is the limiting case $\G \to +\infty$, which recovers portions of Proposition
\ref{normal-form}.  For instance:

\begin{corollary}[Lieb's theorem generalised, gaussian-extremisable case]\label{lieb-nondeg-gen}  Let $\hfill$
$(\B,\p,\G)$ be a generalised Brascamp--Lieb datum
which is non-degenerate and gauss- ian-extremisable.  Then $\BL(\B,\p,\G) = \BL_\g(\B,\p,\G)$.  In particular $\BL(\B,\p,\G)$ is finite in this case.
\end{corollary}

\begin{proof}  We may assume that $p_j > 0$ for all $j$ since otherwise we can simply omit any exponents for which $p_j = 0$.  Let
$A_1,\ldots, A_m > 0$ be an extremiser to \eqref{BLfunctional-gen}.
Since $\B$ is non-degenerate, we have $\bigcap_{j=1}^\infty \ker(B_j) = \{0\}$.  In particular if we set $M: H \to H$ to be the
transformation $M := \sum_{j=1}^m p_j B_j^* A_j B_j$ then $M$ is 
positive definite.

Taking logarithms in \eqref{BLfunctional-gen}, we see that $\A$ is a local maximiser for the quantity
$$ (\sum_{j=1}^m p_j \log \mbox{$\det_{H_j}$} A_j) - \log \mbox{$\det_H$} \sum_{j=1}^m p_j B_j^* A_j B_j$$
subject to the constraint $\A \leq \G$.  Now let us fix a $j$. Let
$V_j \subseteq H$ denote the kernel of the positive definite
operator $G_j - A_j$, let $\iota_j: V_j \to H_j$ be the inclusion
map, and let $Q_j: H_j \to H_j$ be an arbitrary self-adjoint
transformation which is negative definite on $V_j$ (i.e.
$\iota_j^* Q_j \iota_j \leq 0$). Then $0 \leq A_j + \eps Q_j \leq
G_j$ for all sufficiently small $\eps > 0$, and hence
$$ \frac{d}{d\eps} p_j \log \mbox{$\det_{H_j}$} (A_j + \eps Q_j) - \log \mbox{$\det_H$} (M + \eps p_j B_j^* Q_j B_j)|_{\eps = 0} \geq 0.$$
Arguing as in Proposition \ref{normal-form} we then conclude that
$$ \tr_{H_j}( (A_j^{-1} - B_j M^{-1} B_j^*) Q_j ) \geq 0 \hbox{ whenever } Q_j \hbox{ is negative definite on } V_j;$$
In particular, this trace is positive whenever $Q_j$ is negative definite on $H_j$, which implies \eqref{bam},
Also, by considering both $Q_j$ and $-Q_j$ we see that
$$ \tr_{H_j}( (A_j^{-1} - B_j M^{-1} B_j^*) Q_j ) = 0 \hbox{ whenever } \iota_j^* Q_j \iota_j = 0, Q_j = Q^*_j$$
which by negative definiteness of $A_j^{-1} - B_j M^{-1} B_j^*$ implies that
$$ A_j^{-1} - B_j M^{-1} B_j^* = \iota_j N_j \iota_j^*$$
for some transformation $N_j: V_j \to V_j$.  In particular, since
$G_j - A_j$ vanishes on $V_j$, so $\iota_j (G_j - A_j) = 0$ and hence by self-adjointness $\iota_j^* (G_j - A_j) = 0$.
Thus we obtain \eqref{caj}.  Applying Corollary \ref{flowform-bl}, and the hypothesis that $\A$ extremises \eqref{BLfunctional-gen},
we conclude
$$ \BL(\B,\p,\G) \leq \BL_\g(\B,\p;\A) = \BL_\g(\B,\p,\G)$$
and the claim follows from \eqref{blg-trivial-gen}.
\end{proof}

Notice that if in addition to the hypotheses of Corollary \ref{lieb-nondeg-gen} we impose \eqref{scaling}, then 
any gaussian extremiser for $(\B,\p,\G)$ is also one for $(\B,\p)$ since \eqref{bam} and invertibility 
of $M$ together imply gaussian-extremisability by Proposition \ref{normal-form}.

Now we remove the condition that an extremiser exists.  The analysis here is in fact somewhat simpler 
than in the non-regularised case,
because it turns out there is a large class of generalised Brascamp--Lieb data, namely the 
\emph{localised} generalised Brascamp--Lieb data, which
are gaussian-extremisable and are in some sense ``dense'' in the space of all generalised Brascamp--Lieb 
data. (The reason for the nomenclature ``localised'' will become clear below.)
 
\begin{definition}[Localised data] A generalised Brascamp--Lieb datum
$(\B,\p,\G)$ is said to be \emph{localised} if there exists $1 \leq j \leq m$ such that $p_j = 1$, $H_j = H$, and $B_j = \id_H$.
\end{definition}

\begin{lemma}\label{localised-lieb} Let $(\B,\p,\G)$ be a localised 
generalised Brascamp--Lieb datum.  Then $(\B,\p,\G)$ is gaussian-extremisable,
and thus (by Corollary \ref{lieb-nondeg-gen}) we have $\BL(\B,\p,\G)$ 
$\hfill$ $= \BL_\g(\B,\p,\G) < \infty$.
\end{lemma}

\begin{proof}  We may assume that $p_j > 0$ for all $j$, since we can drop all indices for which $p_j = 0$.  We can also assume that $m > 1$
since the $m=1$ case is trivial.
Without loss of generality we may assume that $m$ is the localising index, thus $p_m = 1$, $H_m = H$, 
and $B_m = \id_H$.
We can rewrite \eqref{BLfunctional-gen} as
$$
\BL_\g( \B,\p,\G ) = \sup_{\A \leq \G} \mbox{$\det^{-1/2}_H$} (\id_H +
M^{1/2} A_m^{-1} M^{1/2}) \prod_{j=1}^{m-1} (\mbox{$\det_{H_j}$}
A_j)^{p_j/2}$$ where $\A = (A_j)_{1 \leq j\leq m}$ and $M :=
\sum_{j=1}^{m-1} p_j B_j^* A_j B_j$.  Observe that for fixed
$A_1,\ldots,$ $A_{m-1}$, the quantity $\det^{-1/2}_H (\id_H + M^{1/2}
A_m^{-1} M^{1/2})$ with $0 < A_m \leq G_m$ is maximised when $A_m
= G_m$.  Thus we have
$$
\BL_\g( \B,\p,\G ) = \sup_{{0 < A_j \leq G_j, (1 \leq j
\leq m-1)}} \mbox{$\det^{-1/2}_H$}(\id_H + M^{1/2} G_m^{-1} M^{1/2})
\prod_{j=1}^{m-1} (\mbox{$\det_{H_j}$} A_j)^{p_j/2}.$$ Observe
that this expression extends continuously to semi-definite $A_j$,
and thus
$$
\BL_\g( \B,\p,\G ) = \sup_{0 \leq A_j \leq G_j, ( 1 \leq
j \leq m-1)} \mbox{$\det^{-1/2}_H$} (\id_H + M^{1/2} G_m^{-1} M^{1/2})
\prod_{j=1}^{m-1} (\mbox{$\det_{H_j}$} A_j)^{p_j/2}.$$ Thus an
extremiser exists in the range $0 \leq A_j \leq G_j$.  But the
expression in the supremum vanishes when $\det(A_j) = 0$ for some
$j$, and hence at the extremum we have $A_j > 0$ for all $1 \leq j
\leq m-1$.  The claim follows.
\end{proof}

\begin{lemma}[Approximation by localised data]\label{localised-approx}
Let $(\B,\p,\G)$ be a generalised Brascamp--Lieb datum.  Let $(\B,\id_H)$ be 
the $(m+1)$-transformation formed by appending
the identity operator $B_{m+1} = \id_H: H \to H$ to the $m$-transformation $\B$ (thus $H_{m+1}=H$), and for any real number $\lambda > 0$, let
$(\G,\lambda \id_H)$ be the $(m+1)$-type formed by appending the dilation operator $\lambda \id_H: H \to H$ to the $m$-type $\G$.  Then
$$ \BL(\B,\p,\G) = \lim_{\lambda \to 0} \lambda^{-\dim(H)/2} \BL((\B,\id_H), (\p,1), (\G,\lambda \id_H))$$
and
$$ \BL_\g(\B,\p,\G) = \lim_{\lambda \to 0} \lambda^{-\dim(H)/2} \BL_\g((\B,\id_H), (\p,1), (\G,\lambda \id_H)).$$
\end{lemma}

\begin{proof}
We begin with the first equality.  Let $\f = (f_j)_{1 \leq j \leq m}$ be a normalised input of type $\G$.
Then from \eqref{BL-gen} and
\eqref{hilbert} we see that
$$ \int_H \exp(-\lambda \|x\|_H^2) \prod_{j=1}^m (f_j \circ B_j)^{p_j}(x)\ dx \leq \lambda^{-\dim(H)/2} \BL((\B,\id_H), (\p,1), (\G,\lambda \id_H))$$
for all $\lambda > 0$.  Letting $\lambda \to 0$ and using monotone convergence and \eqref{BL-gen} we conclude
\begin{equation}\label{bl-lower}
\BL(\B,\p,\G) \leq \liminf_{\lambda \to 0} \lambda^{-\dim(H)/2} \BL((\B,\id_H), (\p,1), (\G,\lambda \id_H)).
\end{equation}
Conversely, let $\lambda > 0$ and $f_{m+1}: H \to \R^+$ be a function of type $\lambda \id_H$, 
then we can write $$ f_{m+1}(x) = \lambda^{\dim(H)/2} \int_H \exp(-\lambda \|x-y\|_H^2)\ d\mu(y).$$
By the Fubini-Tonelli theorem and \eqref{hilbert} we have $\int_H f_{m+1} = \mu(H)$.  
Now if $(f_j)_{1 \leq j \leq m}$ is of type $\G$ and we set $B_{m+1} = \id_H$ and $p_{m+1} = 1$, 
then by the Fubini-Tonelli theorem again
\begin{eqnarray*}
\begin{aligned}
\int_H \prod_{j=1}^{m+1} (f_j &\circ B_j)^{p_j}(x)\ dx \\
&= \lambda^{\dim(H)/2}  \int _H \int_H \exp(-\lambda \|x-y\|_H^2)\ d\mu(y)  
\prod_{j=1}^m (f_j \circ B_j)^{p_j}(x)\ dx \\ 
&\leq \lambda^{\dim(H)/2}  \int _H (\int_H \prod_{j=1}^m 
(f_j \circ B_j)^{p_j}(x)\ dx) \ d\mu(y) \\
&\leq \lambda^{\dim(H)/2}  \BL(\B,\p,\G) \prod_{j=1}^m 
(\int_{H_j} f_j)^{p_j} \ \mu(H) \\
&= \lambda^{\dim(H)/2}  \BL(\B,\p,\G) \prod_{j=1}^{m+1} 
(\int_{H_j} f_j)^{p_j}.
\end{aligned}
\end{eqnarray*}

We thus conclude that
$$ \lambda^{-\dim(H)/2} \BL((\B,\id_H), (\p,1), (\G,\lambda \id_H)) \leq \BL(\B,\p,\G).$$
Combining this with \eqref{bl-lower} we obtain the first equality of the lemma.  The second equality is proven in exactly the same way
but with all the $f_j$ (and $f_{m+1})$) constrained to be centred gaussians.  
(The fact that $\exp(-\lambda \|x-y\|_H^2)$ is not centred is of no consequence as we are simply 
bounding it by 1).
\end{proof}

Observe that the data $((\B,\id_H), (\p,1), (\G,\lambda \id_H))$ is manifestly localised.  Thus if we combine Lemma \ref{localised-approx}
with Lemma \ref{localised-lieb}, we immediately obtain

\begin{corollary}[Lieb's theorem generalised] \label{qw} Let $(\B,\p,\G)$ 
be a generalised Bras- camp--Lieb datum.  Then
$\BL(\B,\p,\G) = \BL_\g(\B,\p,\G)$.
\end{corollary}

Combining this corollary with \eqref{tlim} we obtain another proof of Lieb's theorem 
(Theorem \ref{lieb-thm}). A different proof of this result (along the lines 
of \cite{Barthe}) has recently been obtained by Valdimarsson, \cite{V}. Another application -- moving now to the setting of 
localised data -- is the following result, which can also be found implicitly in \cite{L}.

\begin{corollary}[Lieb's theorem localised]\label{llt}  Let $(\B,\p)$ be a Brascamp--Lieb datum, and let $G: H \to H$ be 
positive definite.
Then the best constant $0 < K \leq \infty$ in the estimate
\begin{equation}\label{gauss}
\int_H \exp(-\pi \langle Gx, x \rangle_H) \prod_{j=1}^m 
(f_j \circ B_j(x))^{p_j}\ dx 
\leq K \prod_{j=1}^m \left( \int_{H_j} f_j \right)^{p_j}
\end{equation}
is given by
$$ K := \sup_{A_1,\ldots, A_m > 0} \left( \frac{\prod_{j=1}^m (\det_{H_j} A_j)^{p_j}}{\det_H (G + \sum_{j=1}^m p_j B_j^* A_j B_j)} \right)^{1/2}.$$
\end{corollary}

\begin{proof} By testing \eqref{gauss} on centred gaussians we certainly see that we cannot replace $K$ by any better constant.
Thus it will suffice to prove \eqref{gauss} with the specified value of $K$.

Let $\lambda > 0$ be a large parameter.
An inspection of the proof of Lemma \ref{localised-approx} shows that if each $f_j$ is of 
type $\lambda \id_{H_j}$, then we have
$$ \int_H \exp(-\pi \langle Gx, x \rangle_H) 
\prod_{j=1}^m (f_j \circ B_j(x))^{p_j}\ dx 
\leq K_\lambda \prod_{j=1}^m \left(\int_{H_j} f_j \right)^{p_j}$$
where
$$ K_\lambda := \sup_{0 < A_j \leq \lambda \id_{H_j}} \left( \frac{\prod_{j=1}^m (\det_{H_j} A_j)^{p_j}}{\det_H (G + \sum_{j=1}^m p_j B_j^* A_j B_j)} \right)^{1/2}.$$
Taking limits as $\lambda \to \infty$, we obtain the claim.
\end{proof}

We can also obtain a localised version of Theorem \ref{sufbl}. 

\begin{theorem}[Localised necessary and sufficient conditions]\label{sufbl-local} 
Let $(\B,\p)$ be a Brascamp--Lieb datum with $\B$ non-degenerate, and let
$G: H \to H$ be 
positive definite.  Then the estimate \eqref{gauss} holds for some finite constant $K$ if and only if the
inequalities \eqref{dimh}
hold
for all subspaces $0 \subseteq V \subseteq H$.
\end{theorem}

\begin{proof}
Recall that \eqref{dimh} 
asserts that $\dim(H/V) \geq \sum_{j=1}^m p_j \dim(H_j / (B_j V))$.
We can discard those $j$ for which $p_j=0$ or $H_j = \{0\}$ as they have no impact on either claim in the theorem.  By a linear transformation
we may take $G = \id_H$.

The necessity of the conditions \eqref{dimh} can be seen by testing \eqref{gauss} for functions $f_j$ which lie in an $\eps$-neighbourhood
of the unit ball in $B_j V$ and letting $\eps \to 0$; we omit the easy computations.  Now suppose conversely that \eqref{dimh} holds.
Applying Corollary \ref{llt}, we reduce to showing that
$$ \prod_{j=1}^m (\mbox{$\det_{H_j}$} A_j)^{p_j} \leq K^2 \mbox{$\det_H$} (\id_H + M)$$
for all gaussian inputs $\A = (A_j)_{1 \leq j \leq m}$ for some finite
constant $K$, where $M := \sum_{j=1}^m p_j B_j^* A_j B_j$.

We now repeat the proof of Proposition \ref{sufbl-meat}.  By choosing an appropriate orthonormal basis
$e_1,\ldots,e_n \in H$ we may assume that $M = \diag(\lambda_1, \ldots, \lambda_n)$ for some $\lambda_1 \geq \ldots \geq \lambda_n > 0$.
Our task is thus to show that
$$ \prod_{j=1}^m (\mbox{$\det_{H_j}$} A_j)^{p_j} = O( \prod_{i=1}^n (1 + \lambda_i) ).$$
Here we allow the constants in $O( \ )$ to depend on the
Brascamp--Lieb data. Applying Lemma \ref{babygauss}, we can find
$I_j \subseteq \{1,\ldots,n\}$ for each $1 \leq j \leq m$ of
cardinality $|I_j| = \dim(H_j)$ obeying \eqref{p-comb-alt-2} and
\eqref{bjei}.  From the arguments in Proposition \ref{sufbl-meat}
we have
$$ \prod_{j=1}^m (\det A_j)^{p_j} \leq C \prod_{i=1}^n \lambda_i^{\sum_{j=1}^m p_j |I_j \cap \{i\}|}.$$
Writing $\mu_i := 1 + \lambda_i$ and $c_i := \sum_{j=1}^m p_j |I_j \cap \{i\}|$, it then suffices to show that
$$ \prod_{i=1}^n \mu_i^{c_i} \leq \prod_{i=1}^n \mu_i.$$
We can telescope the right-hand side as
$$ (\prod_{i=1}^n \mu_i) \prod_{k=1}^n (\frac{\mu_{k+1}}{\mu_k})^{k - c_1 - \ldots - c_k}$$
where we adopt the convention $\mu_{n+1} := 1$.  But from \eqref{p-comb-alt-2} we have $k - c_1 - \ldots - c_k \geq 0$ for all $k$,
and from construction we have $\frac{\mu_{k+1}}{\mu_k} \leq 1$ for all $k$.  The claim follows.
\end{proof}

\section{Uniqueness of extremals}

In Theorem \ref{extremiser} we gave necessary and sufficient conditions for gaussian extremisers or extremisers to exist.
In this section we address the issue of whether these extremisers are unique.

One trivial source of non-uniqueness is if $p_j = 0$ for some $j$, then $f_j$ can clearly be arbitrary.  However one can simply omit
these indices $j$ to eliminate this source of non-uniqueness, and thus we shall only consider the case where $p_j > 0$ for all $j$.

From Theorem \ref{extremiser} (and Lemma \ref{scaling-lemma}) it now suffices to consider the 
case of geometric Brascamp--Lieb data $(\B,\p)$, which
has the obvious gaussian extremisers $(\lambda \id_{H_j})_{1 \leq j \leq m}$ for any $\lambda > 0$, and similarly the obvious
extremisers $(c_j \exp( - \lambda \| x-B_j x_0 \|_{H_j}^2 ))_{1 \leq j \leq m}$ for any $\lambda, c_1,\ldots, c_m > 0$ and $x_0 \in H$; compare
with Lemma \ref{closures} and Lemma \ref{closures-gauss}.  The natural question is whether any other extremisers exist.

In the gaussian case, the answer is provided by the following proposition, which can be viewed as a continuation of Proposition \ref{normal-form}
in the geometric case.

\begin{proposition}[Characterisation of gaussian extremisers]\label{extremiser-char}  Let $(\B,\p)$ 
be a geometric Brascamp--Lieb datum with $p_j > 0$ for all $1 \leq j \leq m$.  Let $\A = (A_j)_{1 \leq j \leq m}$ be a gaussian input, and let $M: H \to H$ be a 
positive definite transformation.  Then the following seven statements are logically equivalent.
\begin{enumerate}
\item $\A$ is a gaussian extremiser, and $M = \sum_{j=1}^m p_j B_j^* A_j B_j$.
\item $M = \sum_{j=1}^m p_j B_j^* A_j B_j$, and $A_j^{-1} = B_j M^{-1} B_j^*$ for all $1 \leq j \leq m$.
\item  $B_j M = A_j B_j$ and $A_j^{-1} = B_j M^{-1} B_j^*$ for all $1 \leq j \leq m$.
\item  $M B_j^* = B_j^* A_j$ and $A_j^{-1} = B_j M^{-1} B_j^*$ for all $1 \leq j \leq m$.
\item  $M$ leaves the space $B_j^* H_j$ invariant for all $1 \leq j \leq m$, and $A_j^{-1} = B_j M^{-1} B_j^*$ for all $1 \leq j \leq m$.
\item  $M$ leaves the space $\ker(B_j)$ invariant for all $1 \leq j \leq m$, and $A_j^{-1} = B_j M^{-1} B_j^*$ for all $1 \leq j \leq m$.
\item  Each proper eigenspace of $M$ is a critical subspace for $(\B,\p)$, and $A_j^{-1} = B_j M^{-1} B_j^*$ 
for all $1 \leq j \leq m$.
\end{enumerate}
\end{proposition}

We remark that in (g), we use the term eigenspace to denote a {\em{maximal}} 
subspace consisting of eigenvectors of $M$ with a common eigenvalue.

\begin{proof}  For purposes of visualisation, the reader may wish to follow the arguments below in the case when
the $H_j$ are subspaces of $H$ and $B_j$ is an orthogonal projection,
so that $B^*_j$ is simply the inclusion map from $H$ to $H_j$.  (The general case is in fact equivalent to this case after some Euclidean isomorphisms.)
However for notational reasons it is slightly more convenient for us to keep $H_j$ and $H$ separate from each other.

The equivalence of (a) and (b) follows from Proposition \ref{normal-form}.  Now to see that (c) $\implies$ (b), we observe that if
$B_j M = A_j B_j$ then $\sum_{j=1}^m p_j B_j^* A_j B_j = \sum_{j=1}^m p_j B_j^* B_j M = M$ by \eqref{pbb}.

From duality we see that (c) and (d) are equivalent.   Next we prove that (b) $\implies$ (d).  Using the hypothesis
$M = \sum_{j=1}^m p_j B_j^* A_j B_j$ and \eqref{pbb}, we compute
\begin{align*}
\sum_{j=1}^m p_j &\tr_H( M^{-1} (M B_j^* - B_j^* A_j) (M B_j^* - B_j^* A_j)^* ) \\
&= \tr_H( \sum_{j=1}^m p_j M^{-1} (M B_j^* - B_j^* A_j) (B_j M - A_j B_j) ) \\
&=\tr_H( \sum_{j=1}^m p_j B_j^* B_j M - p_j M^{-1} B_j^* A_j B_j M - p_j B_j^* A_j B_j + p_j M^{-1} B_j^* A_j^2 B_j ) \\
&= \tr_H( M - M - M + M^{-1} \sum_{j=1}^m p_j B_j^* A_j^2 B_j ) \\
&= - \sum_{j=1}^m p_j \tr_H( B_j^* A_j B_j ) + \sum_{j=1}^m p_j \tr_H( M^{-1} B_j^* A_j^2 B_j ) \\
&= - \sum_{j=1}^m p_j \tr_{H_j}( A_j B_j B_j^* ) + \sum_{j=1}^m p_j \tr_{H_j}( B_j M^{-1} B_j^* A_j^2 ).
\end{align*}
Since $B_j B_j^* = \id_{H_j}$ and $B_j M^{-1} B_j^* = A_j^{-1}$, we conclude
\begin{align*}
\sum_{j=1}^m p_j &\tr_H( M^{-1} (M B_j^* - B_j^* A_j) (M B_j^* - B_j^* A_j)^* ) \\
&= -\sum_{j=1}^m p_j \tr_{H_j}(A_j) + \sum_{j=1}^m p_j \tr_{H_j}(A_j) \\
&= 0.
\end{align*}
On the other hand, $M^{-1}$ and $(M B_j^* - B_j^* A_j) (M B_j^* - B_j^* A_j)^*$ are positive semi-definite operators, and thus their product
has non-negative trace.  Since $p_j > 0$, we conclude that
$$ \tr_H( M^{-1} (M B_j^* - B_j^* A_j) (M B_j^* - B_j^* A_j)^* ) = 0 \hbox{ for all } j.$$
Since $M^{-1}$ is 
positive definite, we conclude that
$$ (M B_j^* - B_j^* A_j) (M B_j^* - B_j^* A_j)^* = 0 \hbox{ for all } j,$$
and hence $M B_j^* - B_j^* A_j = 0$.  This gives (d).

The implication (d) $\implies$ (e) is immediate, since $M B_j^* H_j = B_j^* A_j H_j = B_j^* H_j$.  Now we show
that (e) $\implies$ (c); thus suppose that
$M$ leaves $B_j^* H_j$ invariant.  Since $B_j^*$ is an isometry, we have $(B_j^* H_j)^\perp = \ker(B_j)$; since $M$ is self-adjoint
we conclude that $M$ also leaves $\ker(B_j)$ invariant.  Since $M$ is invertible, we see that $M^{-1}$ will then also leave the spaces
$B_j^* H_j$ and $\ker(B_j)$ invariant.  From $A_j^{-1} = B_j M^{-1} B^*_j$ we conclude that
$$(A_j B_j - B_j M) (M^{-1} B^*_j) = A_j B_j M^{-1} B^*_j - B_j B^*_j = A_j A_j^{-1} - \id_{H_j} = 0$$
and hence $A_j B_j - B_j M$ vanishes on $M^{-1} B^*_j H_j = B^*_j H_j$.  Also, since $M^{-1}$ preserves $\ker(B_j)$ we see that
$B_j M$ and $A_j B_j$ both annihilate $\ker(B_j)$.  Thus $A_j B_j - B_j M$ vanishes identically, and we conclude (c).

By duality we see that (e) and (f) are equivalent.
Now we show that (e) and (f) together imply (g).
Let $V$ be a proper eigenspace of $M$.  Let $W_j = B^*_j H_j$. Since $M$ 
leaves $W_j$ invariant and $V$ is maximal, we have 
$W_j = (W_j \cap V) \oplus (W_j \cap V^{\perp})$. 
Similarly we have 
$W_j^{\perp} = (W_j^{\perp} \cap V) \oplus (W_j^{\perp} \cap V^{\perp})$.
Thus $H = W_j \oplus W_j^{\perp}  = (W_j \cap V) \oplus 
(W_j \cap V^{\perp}) \oplus (W_j^{\perp} \cap V) \oplus 
(W_j^{\perp} \cap V^{\perp})$. Now let $\Pi: H \to H$ be the orthogonal 
projection onto $V$;  $B_j^* B_j: H \to H$ is the orthogonal projection 
onto $W_j$. Then $\Pi$ and $B_j^* B_j$ have a common orthonormal set of 
eigenvectors and hence they commute. Consequently, $B_j^* B_j \Pi$ is the orthogonal
projection onto $B_j^* B_j V$, and $\dim (B_j^* B_j V) = \tr_{H} (B_j^* B_j \Pi)$. 
Now $B_j^*$, being an isometry, is injective. So
$\dim (B_j V) = \dim (B_j^* B_j V) = \tr_{H} (B_j^* B_j \Pi)$.
Multiplying this by $p_j$ and summing, we see from \eqref{pbb} That
$$ \sum_{j=1}^m p_j \dim(B_j V) = \tr_H( \sum_{j=1}^m p_j B_j^* B_j \Pi ) 
= \tr_H (\Pi) = \dim(V).$$
Thus $V$ is a critical subspace as desired.

Finally, we show that (g) implies (e) and (f).  Decomposing $M$ as a direct sum of scalar dilation maps on 
the eigenspaces, we see that it suffices
to show that each eigenspace $V$ splits as the direct sum of a subspace of $\ker(B_j)$ and a subspace 
of the orthogonal complement
$B^*_j H_j$.  But by hypothesis the eigenspace $V$ is a critical subspace, and hence by Lemma 
\ref{critpaircrit-converse}, $( V, \ V^\perp \ )$ is a critical pair.  The claim follows.
\end{proof}

The equivalence of (a) and (g) provides a means to construct all the gaussian extremisers 
$\A$ of a geometric Brascamp--Lieb datum $(\B,\p)$. Namely,
we first decompose $\B$ into the direct sum of indecomposable mutually orthogonal spaces, and then let 
$M$ be the direct sum of arbitrary positive scalar dilations on these indecomposable components.  
One then defines the $A_j$ by the formula $A_j^{-1} = B_j M^{-1} B_j^*$.  As a corollary
we obtain

\begin{corollary}[Uniqueness of gaussian extremisers]\label{gauss-unique}  Let $(\B,\p)$ be a 
gaussian-extremisable Brascamp--Lieb datum with $p_j > 0$ for all $j$.  Then the following three statements are equivalent.
\begin{itemize}
\item[(a)] The gaussian extremisers $\A$ of $(\B,\p)$ are unique (up to scaling $\A \mapsto \lambda \A$).
\item[(b)] $\B$ is indecomposable.
\item[(c)] $(\B,\p)$ is simple.
\end{itemize}
\end{corollary}

We remark that this result was obtained in the rank-one case by Barthe \cite{Barthe}.

\begin{proof}  The implication (b) $\implies$ (c) follows from Theorem \ref{extremiser}
(or Lemma \ref{critpaircrit-converse}).

Now we prove (c) $\implies$ (a).  If $(\B,\p)$ is simple, then by statement (g) of 
Proposition \ref{extremiser-char} there are no proper eigenspaces 
of $M$. Consequently the unique eigenspace of $M$ must be the whole of $H$, i.e., $M$ is a scalar 
multiple of the identity. Proposition \ref{extremiser-char} now gives the desired uniqueness. 

Now we prove (a) $\implies$ (b) Conversely, if $\B$ is decomposable, then by Lemma \ref{geom-crit} 
we have a proper critical pair $(V, W)$, which we can choose to be orthogonal complements, 
from Lemma \ref{critpaircrit-converse}.
One can now define a two-parameter family of positive-definite operators $M$ with eigenspaces 
$V,W$, and each one generates a gaussian extremiser by Proposition \ref{extremiser-char}.  The claim 
follows.
\end{proof}

Now we address the issue of uniqueness of extremisers which are not gaussian.  Here the situation appears to be significantly more complicated,
especially in the decomposable case.  For instance, for the Loomis--Whitney inequality \eqref{lw-eq} the extremisers take the form
$$ f(y,z) = a G(y) H(z); \quad g(x,z) = b F(x) H(z); \quad h(x,y) = c F(x) G(y)$$
for arbitrary scalars $a,b,c > 0$ and non-negative integrable functions $F,G,H$ of one variable.  Even in the
indecomposable case there can be plenty of extremisers.  For instance, as is well known, for H\"older's inequality (Example \ref{holder-ex}) one has an
extremiser whenever the functions $f_j$ are scalar multiples of each other, even in the one-dimensional case which is irreducible and thus has unique
gaussian extremisers up to scaling.  However, in the rank one case it is known that for any irreducible 
Brascamp--Lieb datum whose domain $H$ has
dimension strictly larger than one, the only extremisers are the gaussian ones; see \cite[Theorem 4]{Barthe}.  Further analysis of the extremisers in the rank one case was conducted to \cite{CLL}.  We cannot fully generalise this analysis
to the higher rank case.  We do however have the following partial result.

\begin{theorem}[Uniqueness of extremisers]\label{extreme-unique}  Let $(\B,\p)$ be an extremisable Brascamp--Lieb datum with $0 < p_j < 1$ for all $j$.
Suppose also that the spaces $B_j^* H_j$ are all disjoint except
at 0, thus $B_i^* H_i \cap B_j^* H_j = \{0\}$ whenever $1 \leq i <
j \leq m$. Then if $\f = (f_j)$ is an extremising input, then all the
$f_j$ are gaussians, thus there exist real numbers $c_j > 0$,
positive definite transformations $A_j: H_j \to H_j$, and points
$x_j \in H_j$ such that $f_j(x) = c_j \exp( - \pi \langle A_j
(x-x_j), (x-x_j) \rangle_{H_j} )$.
\end{theorem}

\begin{remark}
Remark \ref{bar} tells us that the permitted $x_j$'s in the 
conclusion of the theorem are precisely of the form $B_j \overline{w}$ as $\overline{w}$ varies 
over $H$. 
\end{remark}

\begin{remark} The condition $p_j < 1$ is automatic if one also assumes that $(\B,\p)$ is simple, $\dim(H) > 1$, and that $\dim(H_j) > 0$ for all $j$, as can be seen by testing \eqref{dimension} on one-dimensional subspaces of $H$ not in the kernel of $B_j$.  In the rank one case (Example \ref{bbl}), the condition that the $B_j^* H_j$ are all disjoint amounts to the assertion that no two of the vectors $v_1,\ldots,v_m$ are scalar multiples of each other.  But
in the case when two of the $v_i$ are multiples of each other they
can easily be concatenated using H\"older's inequality; if the
concatenated extremiser is gaussian, then by the converse H\"older
inequality we see that the original extremisers are also gaussian.
Using this observation we can recover the result of \cite[Theorem
4]{Barthe} that in the simple rank one case in dimensions two or
greater, the only extremisers are given by gaussians.  We also
recover the well-known fact that for the sharp Young's inequality
(Example \ref{bl-ex}) with $p_1, p_2, p_3$ strictly less than 1,
the only extremisers are gaussians.
\end{remark}

\begin{proof} Our arguments here are based on a more careful analysis of the heat flow monotonicity argument, as in \cite{CLL}.
By Theorem \ref{extremiser} it suffices to consider the case when $(\B,\p)$ is geometric.
In particular we already have
$\f_0 = (\exp(-\pi \|x\|_{H_j}^2))_{1 \leq j \leq m}$ as an extremising input.  From Lemma \ref{closures} we see that
$\f$ is any extremising input for $(\B,\p)$, then so is $(\f * \f_0) \f_0$.  But since $\f_0$ is Schwartz and $\f$ is integrable,
we see that $(\f * \f_0) \f_0$ is Schwartz and strictly positive (in fact it is of type $\id_{H_j}$).  Also one can easily verify (using the Fourier transform) that $(\f * \f_0) \f_0$ consists of gaussians if and only if $\f$ consists of gaussians.  Thus to prove the claim, it suffices to do so for inputs $\f$ which are Schwartz strictly positive functions of type $\id_{H_j}$.

We now review the proof of Proposition \ref{gbl-prop}, which among other things proves $\BL(\B,\p;\f) \leq 1 = \BL(\B,\p)$.  But since $\f$ is extremising, we must have equality at every stage of the proof of this Proposition.  In particular, since this proof needed \eqref{vvj}
to be non-negative, we now see that in fact \eqref{vvj} needs to be zero everywhere:
$$ \sum_{j=1}^m p_j \langle B_j^* B_j(\vec v - \vec v_j), (\vec v - \vec v_j) \rangle_H = 0.$$
Since $p_j > 0$ for all $j$, we conclude that
$$ \| B_j(\vec v - \vec v_j) \|_H = 0 \hbox{ for all } j.$$
Since $\vec v = \sum_{i=1}^m p_i \vec v_i$, we thus have
$$ B_j \sum_{i=1}^m p_i \vec v_i = B_j \vec v_j \hbox{ for all } j.$$
Since $p_j$ is strictly less than 1, and $\vec v_j$ lies in the range of $B_j^* B_j$, we conclude that
$$ \vec v_j = B_j \sum_{i \neq j} \frac{p_i}{1-p_j} B_j^* B_j \vec v_i.$$
Next, observe from the chain rule that $\vec v_i$ takes values in $B_i^* H_i \subseteq H$, and
we can write $\vec v_i = w_i \circ B_i$ for some function $w_i: H_i \to B_i^* H_i$, thus
\begin{equation}\label{wb}
w_j \circ B_j = B_j \sum_{i \neq j} \frac{p_i}{1-p_j} B_j^* B_j w_i \circ B_i.
\end{equation}
The disjointness of the spaces $B_i^* H_i$ will now force a lot of structure on these functions $w_i$; this is easiest to establish using
Fourier analysis and the theory of distributions\footnote{An alternate approach, which exploits the smoothness of the $w_i$, is to differentiate \eqref{wb} in various directions, designed to eliminate the terms on the right-hand side but not on the left.}.  First observe that since
$\f$ is Schwartz and strictly positive, it is not hard to see from the definition of $\vec v$
that the $w_i$ are continuous and grow at most linearly
(thus we have $\|w_i(x)\|_H \leq C_{\f,i} (1+\|x\|_{H_j})$ for all $x \in H_j$
for some constant $C_{\f,i}$ depending on the input $\f$ and the index $i$).  In particular $w_i$ is a tempered distribution.
We can now take distributional Fourier transforms of \eqref{wb} in the Euclidean space $H$.
The left-hand side is a tempered distribution supported in $B_j^* H_j$, while the right-hand side is supported in $\bigcup_{i \neq j} B_i^* H_i$.
Thus by hypothesis, both distributions are in fact supported on $\{0\}$.  Inverting the Fourier transform, we conclude in particular
that $w_j \circ B_j$ (and hence $w_j$) is a polynomial.  But since $w_j$ has at most linear growth, we conclude that $w_j$ is a linear polynomial.
We now specialise to the limiting time case $t=0$, in which $v_j = - \nabla \log f_j$, taking advantage of the fact that $f_j$ is of type $\id_{H_j}$, to conclude that $\nabla \log f_j$ is a linear polynomial.  Integrating this we see that $\log f_j$ must be a quadratic polynomial; since $f_j$ is integrable, the leading term must be strictly negative definite.  The claim follows by completing the square of the quadratic polynomial.
\end{proof}

\section{Sliding kernels}\label{sliding-sec}

We now recast the monotonicity formulae obtained earlier as a monotonicity property of sliding gaussians.  This leads naturally to the
question of whether such monotonicity also holds for other kernels than gaussians; in the linear case, we will show that this is true
for log-concave kernels.

Let us first return to Proposition \ref{gbl-prop}.  A key component of the proof of that proposition was the claim that the quantity
$$ Q(t) := \int_H \prod_{j=1}^m u_j^{p_j}(t,x)\ dx$$
was monotone non-decreasing in time for $t > 0$, where
$$u_j(t,x) = \frac{1}{(4\pi t)^{\dim(H_j)/2}} \int_{H_j} e^{-\| B_j x - z \|_{H_j}^2/4t} f_j(z)\ dz$$
and $(f_j)_{1 \leq j \leq m}$ was an arbitrary input.  Making the change of variables $s := (4\pi t)^{-1/2}$, $y := sx$, and $v := z$, and $d\mu_j(v) = f_j(v) d\mu(v)$, we thus see that the quantity
$$ \int_H \prod_{j=1}^m \left(\int_{H_j} \exp(- \pi \| B_j y - v s \|_{H_j}^2 )\ d\mu_j(v)\right)^{p_j}\ dy$$
is monotone non-increasing in $s$ for $s > 0$.  We view $s$ as a time parameter, $y$ as a position variable, $v$ as a velocity variable,
and $\mu_j$ as a velocity distribution.  
Each function $\exp(- \pi \| B_j y - v s \|_{H_j}^2 )$ then becomes a ``sliding gaussian''
which equals $\exp(-\pi \|B_j y \|_{H_j}^2)$ at time $s=0$ and then moves with velocity $v$ (where we embed $H_j$ into $H$ using the isometry
$B_j^*$).  The above monotonicity then represents the intuitively plausible fact that the multilinear $L^p$ norm of these sliding gaussians
is maximised at time $s=0$, at which time all the gaussians are centred at the origin.

In light of Proposition \ref{normal-form}, there should be an analogous 
monotonicity of sliding gaussians for any Brascamp--Lieb datum 
$(\B,\p)$ for which one can locate a gaussian input $\A$ obeying 
\eqref{scaling} and the inequalities $B_j^* A_j B_j \leq \sum_{i=1}^m p_i B_i^* A_i B_i$ for all
$j$.  Indeed, in light of Proposition \ref{flowform} (in the special case $\G=\A$), the scaling condition can be dropped:

\begin{proposition}\label{multislide}  Let $(\B,\p)$ be a Brascamp--Lieb datum with all the $B_j$ surjective,
and let $\A = (A_j)_{1 \leq j \leq m}$ be a gaussian input such that $B_j^* A_j B_j \leq 
\sum_{i=1}^m p_i B_i^* A_i B_i$ for all
$j$.  Then for any positive finite Borel measures $d\mu_j$ on $H_j$, the quantity
$$ \int_H \prod_{j=1}^m \left(\int_{H_j} \exp(- \pi \langle A_j (B_j y - v s), (B_j y - v s) \rangle_{H_j} )\ d\mu_j(v)\right)^{p_j}\ dy$$
is monotone non-increasing in $s$ for $s > 0$.
\end{proposition}

\begin{proof} An easy scaling argument shows that it suffices to prove this when $0 < s < (4\pi)^{-1/2}$.  The claim then follows from
Proposition \ref{flowform} in the case $\G=\A$, after the change of variables $s := (4\pi t)^{-1/2}$, $y := sx$, and $v := z$
as in the preceding discussion.
\end{proof}

Now we turn to log-concave kernels.  We begin with a divergence estimate, which is a weak analogue to Corollary \ref{divag}.

\begin{proposition}[Divergence estimate for log-concave kernels]\label{log-concave}  Let $\psi: \R^n \to \R^+$ be a smooth, strictly positive, absolutely integrable
function which is log-concave (thus $\nabla^2 \psi(x) \geq 0$ for all $x$).  Let $\mu$ be a non-zero positive finite Borel measure,
and let $\vec y: \R^n \to \R^n$ be
the centre-of-mass vector field
$$ \vec y(x) := \frac{ \int_{{\bf{R}}^n} \psi(x-y)\ y\ d\mu(y) }{ \int_{{\bf{R}}^n} \psi(x-y)\ d\mu(y)}.$$
Then $\div \vec y \geq 0$, with equality if and only if $\mu$ is a point mass.
\end{proposition}

\begin{proof}  From the definition of $\vec y$ we have
\begin{equation}\label{vec-ident}
\int_{{\bf{R}}^n} \psi(x-y) (\vec y(x) - y)\ d\mu(y) = 0 \hbox{ for all } x \in \R^n.
\end{equation}
Taking divergences of both sides and writing $\nabla \psi = \psi \nabla \log \psi$, we obtain
$$ \int_{{\bf{R}}^n} \psi(x-y) [ \langle \nabla \log \psi(x-y), \vec y(x) - y \rangle + \div(\vec y)(x) ]\ d\mu(y) = 0.$$
On the other hand, from \eqref{vec-ident} again we have
$$ \int_{{\bf{R}}^n} \psi(x-y) \langle \nabla \log \psi(x-\vec y), \vec y(x) - y \rangle\ d\mu(y) = 0.$$
Combining these two equations we obtain
$$ \div(\vec y)(x) = \frac{ \int_{{\bf{R}}^n} \psi(x-y) \langle \nabla \log \psi(x-\vec y) - \nabla \log \psi(x-y), \vec y - y \rangle\ d\mu(y) }
{ \int_{{\bf{R}}^n} \psi(x-y) d\mu(y) }.$$
Since $\nabla^2 \log \psi \geq 0$, we see from the mean-value theorem that
$\langle \nabla \log \psi(x-\vec y) - \nabla \log \psi(x-y), \vec y - y \rangle$ is non-negative, and the claim follows.
\end{proof}

\begin{remark} In the gaussian case $\psi(x) = \exp(-\pi \langle G x, x \rangle)$ the above proposition follows from the
$A=G^{-1}$ case of Corollary \ref{divag} after some simple algebraic manipulations which we omit here.
\end{remark}

As a consequence of this divergence estimate we obtain the following monotonicity formula.

\begin{lemma}[$L^p$ monotonicity for log-concave kernels]\label{logconcave}  Let $\psi: \R^n \to \R^+$ be a strictly positive log-concave function which vanishes at infinity.
Then for any positive finite Borel measure $\mu$ on $\R^n$ and any $p \geq 1$ the quantity
$$ Q(t) = \int_{{\bf{R}}^n} (\int_{{\bf{R}}^n} \psi(x-vt)\ d\mu(v))^p\ dx$$
is non-increasing in time for $t \geq 0$.
\end{lemma}

\begin{remark} The intuition here is that the travelling 
waves $\psi(x-vt)$ are diverging from each other as $t > 0$ increases, and the total mass of
$\int_{{\bf{R}}^n} \psi(x-vt)\ d\mu(v)$ is constant, so the higher $L^p$ norms should decrease.
\end{remark}

\begin{proof}  We may assume of course that $\mu$ is not identically zero.  By a limiting argument we can also assume that $\mu$ is compactly
supported and $\psi$ is rapidly decreasing.  We set
$$ u(t,x) := \int_{{\bf{R}}^n} \psi(x-vt)\ d\mu(v)$$
and
$$ \vec v(t,x) := \frac{ \int_{{\bf{R}}^n} 
v \psi(x-vt)\ d\mu(v)}{ \int_{{\bf{R}}^n} \psi(x-vt)\ d\mu(v) }$$
(note that the denominator is strictly positive by our assumptions on $\psi, d\mu$).  A simple calculation then shows that we have the transport
equation
$$ \partial_t u + \div( \vec v u ) = 0.$$
From Proposition \ref{log-concave} we also have
$$ \div( \vec v ) \geq 0.$$
Also by our assumptions on $\psi, \mu$ we see that $\vec v u^p$ is rapidly decreasing in space.
The claim now follows from Lemma \ref{multitransport} with $m=1$, $\alpha = 0$, $u_1 := u$, and $\vec v_1 := \vec v$, and with the reversed signs.
\end{proof}

\begin{remarks} A similar argument shows that when $0 < p < 1$ the quantity $Q(t)$ is non-increasing in time.  Note that $Q$ is constant in the
boundary case $p=1$.  One can easily prove a similar result when $f$ is merely a positive finite measure rather than a Schwartz function by a standard
limiting argument which we omit here.  In the case when $\psi$ is strictly log-concave, a more refined analysis in the $p \neq 1$ case
also shows that $Q$ is strictly monotone unless $f$ is a point mass; we again omit the details.
\end{remarks}

\begin{remark}  In the gaussian case $\psi(x) = \exp(-\pi \langle A x, x \rangle)$ this proposition is a special case of
Proposition \ref{multislide}.  Indeed, the above arguments can be easily modified to give a direct proof of Proposition \ref{multislide}
from Lemma \ref{multitransport}, using Corollary \ref{divag} instead of Proposition \ref{log-concave}; we omit the details.
\end{remark}

We can apply this lemma to concrete log-concave kernels such as the one-dimensional heat kernel $\psi(x) := e^{-\pi x^2}$.

\begin{corollary}\label{harm-cor}  Let $f: \R \to \R^+$ be a positive finite measure 
and $p \geq 1$.  Similarly,
if $u: \R^+ \times \R \to \R^+$ denotes the heat extension
$$ u(t,x) := \frac{1}{\sqrt{4\pi t}} \int_{{\bf{R}}} e^{-|x-y|^2/4t} f(y)\ dy$$
then $t^{1/2p'} \| \phi(\cdot,t)\|_{L^p({\bf{R}})}$ is non-decreasing in $t$.  If $f$ is not a point mass, then these quantities are strictly increasing in $t$.  Here $p'$ is the dual exponent of $p$, defined by $1/p + 1/p' = 1$.  
\end{corollary}

This innocuous-sounding result does not appear to be previously in the literature; it can be derived by explicit computation when $p \geq 1$
is an integer but is not trivial to prove for other values of $p$; one can also establish these results directly from
Lemma \ref{multitransport} of course.  It is an interesting question to ask whether an analogous result holds for the harmonic
extension $\phi(t,x)$ of $u$ (with $t^{1/2p'}$ replaced by $t^{1/p'}$); this corresponds to setting $\psi(x) = \frac{1}{\pi(1+x^2)}$, which is not log-concave.  In this regard, it is a classical result, essentially due to Hardy and Littlewood \cite{HL} (see also \cite{Duren}, Theorem 5.9) 
that if $p \geq 1$, then $t^{1/p'} \|\phi( \cdot, t) \|_p$ is bounded as a function of $s$, and in fact it tends to zero as $s$
decreases to zero.   A similar claim can also be proven easily for the heat extension.  It is also amusing to note that there is a dyadic analogue of
these monotonicity formulae: if for each integer $k$ we let $E_k(f)$ denote the orthogonal 
projection of $f$ onto functions which are constant on dyadic intervals of 
length $2^k$, then it is easy to verify that the quantity $2^{k/p'} \| E_k(f) \|_{L^p({\bf{R}})}$ is non-decreasing in $k$.

\end{document}